\documentclass[nonblindrev]{informs1_noinforms}              
\usepackage{pdfsync}
\usepackage{enumerate,float}

\newtheorem{observation}{Observation}
\newtheorem{interpretation}{Interpretation}

\newcommand{\E}{{\Bbb{E}}}

\newcommand{\X}{{\cal X}}

\def\supp {{\rm supp}}

\def\bbr{{\Bbb{R}}} 
\def\bbe{{\Bbb{E}}} 

\def\bbp{{\Bbb{P}}}

\newcommand{\pr}{\mathbb{P}}

\usepackage{natbib}
 \NatBibNumeric
 \def\bibfont{\small}%
 \bibpunct[, ]{[}{]}{,}{n}{}{,}%

\usepackage[colorlinks=true,breaklinks=true,bookmarks=true,urlcolor=blue,
     citecolor=blue,linkcolor=blue,bookmarksopen=false,draft=false]{hyperref}

\def\EMAIL#1{\href{mailto:#1}{#1}}
\def\URL#1{\href{#1}{#1}}         

\TheoremsNumberedThrough     

\begin{document}


\RUNAUTHOR{Xin and Goldberg}

\RUNTITLE{Distributionally robust inventory control when demand is a martingale}

\TITLE{Distributionally robust inventory control when demand is a martingale}

\ARTICLEAUTHORS{%
\AUTHOR{Linwei Xin}
\AFF{University of Chicago, \EMAIL{linwei.xin@chicagobooth.edu}, \URL{http://faculty.chicagobooth.edu/linwei.xin/}}
\AUTHOR{David A. Goldberg}
\AFF{Cornell University, \EMAIL{dag369@cornell.edu}, \URL{https://people.orie.cornell.edu/dag369/}}
} 

\ABSTRACT{Demand forecasting plays an important role in many inventory control problems.  To mitigate the potential harms of model misspecification in this context, various forms of distributionally robust optimization have been applied.  Although many of these methodologies suffer from the problem of time-inconsistency, the work of \citet{KSLS} established a general time-consistent framework for such problems by connecting to the literature on robust Markov decision processes.

Motivated by the fact that many forecasting models exhibit very special structure, as well as a desire to understand the impact of positing different dependency structures, in this paper we formulate and solve a time-consistent distributionally robust multi-stage newsvendor model which naturally unifies and robustifies several inventory models with demand forecasting.  In particular, many simple models of demand forecasting have the feature that demand evolves as a martingale (i.e. expected demand tomorrow equals realized demand today).  We consider a robust variant of such models, in which the sequence of future demands may be any martingale with given mean and support.  Under such a model, past realizations of demand are naturally incorporated into the structure of the uncertainty set going forwards.

We explicitly compute the minimax optimal policy (and worst-case distribution) in closed form, by combining ideas from convex analysis, probability, and dynamic programming.  We prove that at optimality the worst-case demand distribution corresponds to the setting in which inventory may become obsolete at a random time, a scenario of practical interest.  To gain further insight, we prove weak convergence (as the time horizon grows large) to a simple and intuitive process.  We also compare to the analogous setting in which demand is independent across periods (analyzed previously in \citet{S-12}), and identify interesting differences between these models, in the spirit of the price of correlations studied in \citet{ADSY}.  Finally, we complement our theoretical results by providing a targeted and concise numerical experiment further demonstrating the benefits of our model.
}



\maketitle

%


%
%
%

\setcounter{equation}{0}

\section{Introduction and literature review}\label{sec-intro}
\subsection{Introduction}
The fundamental problem of managing an inventory over time in the presence of stochastic demand is one of the core problems of Operations Research, and related models have been studied since at least the seminal work of \citet{edgeworth1888mathematical}.  In many practical settings of interest, demands are correlated over time (cf. \citet{scarf1959bayes,scarf1960some,iglehart1962optimal}).  As a result, there is a vast literature investigating inventory models with correlated demand, including: studies of the so-called bull-whip effect (cf. \citet{ryan1998analysis,chen2000quantifying,lee2004information}); models with Markov-modulated demand (cf. \citet{feldman1978continuous,iglehart1962optimal,kalymon1971stochastic}); models with forecasting, including models in which demand follows an auto-regressive / moving average (ARMA) or exponentially smoothed process (cf. \citet{JT,Lovejoy,blinder1983inventories,pindyck1982adjustment,miller1986scarf,badinelli1990inventory,gallego2001integrating,levi2008lost,CGSZ2015,Boute18}); and models obeying the Martingale Model of Forecast Evolution (MMFE) and its many generalizations (cf. \citet{HJ,graves1986two,toktay2001analysis,erkip1990optimal,IZ,LSR,MK}).  Although several of these works offer insights into the qualitative impact of correlations on the optimal policy (and associated costs) when managing an inventory over time, these results are typically proven under very particular distributional assumptions, which assume perfect knowledge of all relevant distributions.  
This is potentially a significant problem, since various authors have previously noted that model misspecification when demand is correlated can lead to very sub-optimal policies (cf. \citet{badinelli1990inventory}).  Indeed, the use of such time series and forecasting models in Operations Research practice is well-documented, and concerns over the practical impact of model misspecification have been raised repeatedly in the forecasting and Operations Research literature (cf. \citet{Fildes1995,Fildes2008}).
\\\indent One approach taken in the literature to correcting for such model uncertainty is so-called distributionally robust optimization.  In this framework, one assumes that the joint distribution (over time) of the sequence of future demands belongs to some set of joint distributions, and solves the minimax problem of computing the control policy which is optimal against a worst-case distribution belonging to this set.  Such a distributionally robust approach is motivated by the reality that perfect knowledge of the exact distribution of a given random process is rarely available (cf. \citet{D-87,D-01,Prekopa,Za}). A typical minimax formulation is as below:
\[
   \min_{x\in \X}\max_{Q\in \mathfrak{M}} \bbe_Q\left[f(x,\xi)\right],
\]
where $\X$ is a set of decisions and $\mathfrak{M}$ is a set of probability measures. The objective is to pick a decision $x$ that minimizes the average cost of $f$ under a worst-case distribution.
\\\indent The application of such distributionally robust approaches to the class of inventory control problems was pioneered in \citet{Scarf58}, where it was assumed that $\mathfrak{M}$ contains all probability measures whose associated distributions satisfy certain moment constraints.  Such an approach has been taken to many variants of the single-stage model since then (cf. \citet{GM,gallego1998new,PR,moon1994distribution,yue2006expected, hanasusanto2012distributionally,zhu2013newsvendor}), and the single-stage distributionally robust model is quite broadly understood.
\\\indent The analogous questions become more subtle in the multi-stage setting, due to questions regarding the specification of uncertainty in the underlying joint distribution.  There have been two approaches taken in the literature, depending on whether the underlying optimization model is \emph{static} or \emph{dynamic} in nature.  In a static formulation, one specifies a family of joint distributions for demand over time, typically by e.g. fixing various moments and supports, or positing that the distribution belongs to some appropriately defined ball, and then solves an associated global minimax optimization problem (cf. \citet{DY,SS,Popescu, bertsimas2010models, D-87, D-01, BGK14}).  Such static formulations generally cannot be decomposed and solved by dynamic programming (DP), because the distributional constraints do not contain sufficient information about how the distribution behaves under conditioning.  Put another way, such models generally do not allow for the incorporation of realized demand information into model uncertainty / robustness going forwards (i.e. re-optimization in real time), and are generally referred to as time-inconsistent in the literature (cf. \citet{ACS,IPS,XGS,ES}).  This inability to incorporate realized demand information may make such approaches undesirable for analyzing models which explicitly consider the forecasting of future demand (cf. \citet{KSLS,IPS,S-09}).  We note that although some of these static formulations have been able to model settings in which information is revealed over time and/or temporal dependencies manifest through time series (cf. \citet{SS,lam2011,Gupta2015}), the fundamental inability to incorporate realized demand (in the sense of time in-consistency) remains.
\\\indent Alternatively, in a dynamic formulation, the underlying family of potential joint distributions must admit a certain decomposition which ensures a desirable behavior under conditioning, and a resolution by DP.  The existence of such a decomposition is generally referred to as the \emph{rectangularity property} (cf. \citet{Iyengar,IPS,ES,Shapiro15}).  For example, the set of all joint distributions for the vector of demands $(D_1,D_2)$ such that (s.t.): $\E[D_1] = 1, \E[D_2 | D_1] = D_1$, and $(D_1,D_2)$ has  support on the non-negative integers is rectangular, since the feasible set of joint distributions may be decomposed as follows.  To each possible realized value $d$ for $D_1$ (i.e. each non-negative integer), we may associate a fixed collection $S(d)$ of possible conditional distributions for $D_2$ (i.e. those distributions with mean $d$ and support on the non-negative integers).  Furthermore, every feasible distribution for the vector $(D_1,D_2)$ may be constructed by first selecting a feasible distribution ${\mathcal D}$ for $D_1$ (i.e. any distribution with mean one and support on the non-negative integers), and for each element $d$ in the support of ${\mathcal D}$, setting the conditional distribution of $D_2$ (given $\lbrace D_1 = d \rbrace$) to be any fixed element of $S(d)$.  Moreover, the set of joint distributions constructible in this manner is precisely the set of feasible distributions for the vector $(D_1,D_2)$.  Alternatively, if we had instead required that $\E[D_1] = \E[D_2] = 1$, and $\E[D_1 D_2] = 2$ (with the same support constraint), the corresponding set of feasible joint distributions would \emph{not} be rectangular, as it may be verified that such a decomposition is no longer possible.  For a formal definition and more complete / precise description of the rectangularity property, we refer the reader to cf. \citet{Iyengar, IPS, ES, XGS, NE, Bielecki2014, Shapiro15}, and note that since various communities have studied several closely related notions at different times, a complete consensus on a common rigorous definition has not yet been reached.  As discussed in \citet{NE}, we note that such dynamic formulations are closely related to the theory of stochastic games (cf. \citet{Altman2000}), and have also been studied within the mathematical finance community (cf. \citet{Bay13,BayYao13,Glasserman14,CR,FL,NO,Riedel,TK}).
\\\indent There are several works which formulate DP approaches to distributionally robust inventory models (cf. \citet{gallego2001minimax, ACS, choi2008risk, S-12, XGS}).  More generally, such dynamic problems can typically be formulated as so-called robust Markov decision processes (MDP) (cf. \citet{Iyengar,NE,WKR}), akin to robust formulations considered previously within the control community (cf. \citet{hansensargent2001}).  To the best of our knowledge, the only such work which considers 
applications to correlated demands and forecasting models within Operations Research is that of \citet{KSLS}.  In that work, the authors make the connection between dynamic distributionally robust inventory models and robust MDP, in a very general sense, and prove e.g. the optimality of state-dependent $(s,S)$ policies, where the state incorporates the entire history of demand.  The authors also consider a more specialized model in which the demand in period $t$ may come from any distribution within a certain ball of the empirical histogram of all past demands.  We note that many simple forecasting models common in the Operations Research and statistics literature (e.g. linear time series) have considerable special structure, and to the best of our knowledge precisely how this structure would manifest within the general framework of \citet{KSLS} (and associated general framework of robust MDP) remains an open question.
 \\\indent On a related note, to the best of our knowledge, there seems to have been no systematic study of the impact of positing different joint dependency structures in such multi-stage distributionally robust inventory control problems, e.g. comparing the properties and cost of a minimax optimal policy.  The quest to develop such an understanding in the broader context of stochastic optimization (not specifically inventory control) was initiated in \citet{ADSY}, where the authors define the so-called \emph{price of correlations} as the ratio between the optimal minimax value when all associated random variables (r.v.) are independent, and the setting where they may take any joint distribution belonging to an allowed family.  Although the authors do not look specifically at any inventory problems, they stress the general importance of understanding how positing different joint distribution uncertainty impacts the underlying stochastic optimization.  Also, although several authors have demonstrated that the family of non-stationary history-dependent s-S policies is optimal for a broad class of minimax inventory problems (cf. \citet{ACS,Notzon70,KSLS}), the impact of the posited dependency structure on the corresponding s-S policy remains poorly understood.  We note that although these questions remain open, the general theory of sensitivity analysis for (distributionally robust) optimization problems is very well developed (cf. \citet{BS2013}), and that different yet related questions are the subject of several recent investigations (cf. \citet{Gupta2015}).
\\\indent Combining the above, we are led to the following question.
\\\begin{question}\label{mainquestion1}
Can we construct effective dynamic distributionally robust variants of the structured time series and forecasting models common in the Operations Research and statistics literature?  Furthermore, can we develop a theory of how positing different dependency structures impacts the optimal policy and cost for such models?
\end{question}
\subsection{Additional literature review}
Highly relevant to any discussion of distributionally robust optimization (either static or dynamic) is the vast literature on classical (i.e. deterministic) robust optimization, in which the only constraints made are on the supports of the associated r.v.s (cf. \citet{BEN,BBC}).  In this setting, the worst-case distribution is always degenerate, namely a point-mass on a particular worst-case trajectory.  Such models often lead to tractable global optimization problems for fairly complex models, and have been successfully applied to several inventory settings (cf. \citet{kasugai1961note,aharon2009robust,BGNV,BIP,BT,BHWMR,carrizosaa2014robust,Sol.15,Sun17,Mamani17}).  Indeed, these models tend to be more tractable than corresponding distributionally robust formulations, where the question of tractability (under both static and dynamic formulations) may be more sensitive to the particular modeling assumptions (cf. \citet{SS,BSZ2014,WKR,DY}).  In spite of their potential computational advantages, one drawback of such classical robust approaches is that they may be overly conservative, and unable to capture the stochastic nature of many real-world problems (cf. \citet{SS}).  We note that the precise relationship between classical robust and distributionally robust approaches remains an intriguing open question.  On a related note, questions of robust time series models and their applications to inventory models, similar in spirit to those we will consider in the present work, were recently studied in \citet{carrizosaa2014robust}, at least for single-stage models.  There the authors formulated a general notion of robust time series, and prove that the corresponding single-stage minimax newsvendor problem can be formulated as a convex optimization problem by using a budget-of-uncertainty approach similar to that studied in \cite{bandi12}.  Such an approach has also been used to study robust versions of multi-period inventory models in which demand is correlated (cf. \cite{Mamani17}).  As regards connections to our earlier discussion of static and dynamic formulations in the distributionally robust setting, we note that the potential modeling benefits of being able to dynamically update one's uncertainty sets has also been considered in the deterministic robust setting (cf. \citet{BV15,Wag.17}).  Indeed, an approach to robust multi-period inventory problems based on deterministic robust optimization / the budget-of-uncertainty approach and optimal control was recently presented in \cite{Wag.17}, where we note that their approach and findings are in a sense incomparable to our own.  
\\\indent Also relevant to questions involving distributionally robust optimization is the search for ``optimal" probability bounds.  In particular, some of the most celebrated results in probability theory are of the form ``any random vector with property A has probability at most p of belonging to set S" for appropriate choices of A,p,S.  That obtaining tight versions of these bounds can be phrased as an appropriate distributionally robust optimization problem is by now well-known within the Operations Research community, especially since this connection was clarified in the seminal work of \citet{BP2005}.  Such problems, e.g. those considered in \citet{BP2005}, typically involve formulating a static optimization problem in which one specifies knowledge or bounds on certain moments, covariances, supports, etc., and frames finding a ``worst-case" distribution satisfying these conditions as a convex optimization problem.  Recent work of \citet{beigblock2014} has brought to light the fact that certain probability inequalities, specifically well-known inequalities for martingales such as the celebrated Burkholder inequality, can be similarly formulated using distributionally robust optimization, but are more amenable to a dynamic formulation, due to the fundamentally conditional structure of the martingale property.  Another related line of work is that on so-called optimal transport problems, especially that subset focusing on so-called optimal martingale transport, and we refer the interested reader to \citet{Dolinsky14} and \citet{beigblock2015} for details and additional references.
\\\indent We close this overview by noting that there is a vast literature on time series in the statistics community (cf. \citet{BJ2011}), and we make no attempt to survey that literature here.  These investigations include considerable work on robust approaches (cf. \citet{Mar2006}), where different lines of research use different notions of robustness (cf. \citet{Stock87}).  Most relevant to our own investigations are the many such works related to distributional misspecification, which has been investigated within the statistics community for over 50 years (cf. \citet{Tukey60,Huber64,Stock87,MR15}), and includes the analysis of minimax notions of robustness (cf. \citet{Franke1985,Luz2014}) and so-called robust Kalman filters (cf. \citet{Ruck2010}).  In spite of this vast literature, it is worth noting that even recently, there have been several calls for the development of more tools and methodologies for the robust analysis of time series data.  For example, in \citet{GGH06}, it is noted that a persistent problem in the application of time series models is \textbf{``the frequent misuse of methods based on linear models with Gaussian i.i.d. distributed errors."}  A similar commentary was given very recently in \citet{Guerrier14}, in which the authors note that 
\textbf{``the robust estimation of time series parameters is still a widely open topic in statistics for various reasons."}  Such calls for the development of new methodologies further serve to underscore the potential impact of combining ideas from robust optimization with time series analysis, e.g. through Question\ \ref{mainquestion1}.
\subsection{Our contributions}
In this paper we take a step towards answering Question\ \ref{mainquestion1} by formulating and solving a dynamic distributionally robust multi-stage newsvendor model which naturally unifies the analysis of several inventory models with demand forecasting.  
In particular, we consider a dynamic distributionally robust multi-stage inventory control problem (with backlogged demand) in which the sequence of future demands may be any martingale with given mean and support (assumed the same in every period).  Note that the martingale demand assumption has been widely used in many demand forecasting models in the literature, including the additive-MMFE in which the demand process is given by $D_t-D_{t-1}=\epsilon_t$ where $\{\epsilon_t\}_{t\geq 1}$ is a sequence of i.i.d. r.v.s with mean zero; and the multiplicative-MMFE in which the demand process is given by $\frac{D_t}{D_{t-1}}=\exp\left(\eta_t\right)$ where $\{\eta_t\}_{t\geq 1}$ is a sequence of i.i.d. r.v.s with $\bbe\left[\exp\left(\eta_1\right)\right]=1.$  Here we refer the interested reader to \citet{WAK} for further discussion of such models and their study in the inventory literature.  Our contributions are four-fold.  First, we explicitly compute the minimax optimal policy (and associated worst-case distribution) in closed form.  Our main proof technique involves a non-trivial induction, combining ideas from convex analysis, probability, and DP.  Second, we prove that at optimality the worst-case demand distribution corresponds to the setting in which the inventory may become obsolete at a random time, a scenario of practical interest which has been studied previously in the literature (cf. \citet{Pierskalla69,Brown1964,JoglekarLee93,SZ96,Cerag2010}).  Third, to gain further insight into our explicit solution, we prove weak convergence (as the time horizon grows large) to a simple and intuitive process.  Fourth, in the spirit of \citet{ADSY}, we compare to the analogous setting in which demand is independent across periods (analyzed previously in \citet{S-12}), and identify interesting differences between these two models. In particular, we prove that when the initial inventory level is zero, the minimax optimal cost is maximized under the dependency structure of the independent-demand model.  More precisely, the minimax optimal cost when nature is restricted to selecting a joint distribution (for demand) under which demand is independent over time (and satisfies the given support constraints) is equivalent to the minimax optimal cost when nature need only satisfy the support constraints (and may select any joint distribution whatsoever).  Furthermore, we compute the limiting ratio (as the time horizon grows large) of the minimax optimal value under the martingale-demand model to that under the independent-demand model, which has a simple closed form and evaluates to exactly $\frac{1}{2}$ in the perfectly symmetric case.  Interestingly, we find that under different initial conditions, the situation may in fact be reversed, showing that such comparative statements may be quite subtle.  Finally, we complement our theoretical results by providing a targeted and concise numerical experiment further demonstrating the benefits of our model.  Indeed, our experiments find that when the true demand distribution has the dependency structure of an additive-MMFE model, the robust minimax-optimal policy we devise in this work significantly outperforms the previously studied robust policy which is minimax-optimal under an independence assumption.

\subsection{Outline of paper}\label{outlinesec}
Our paper proceeds as follows.  We begin by reviewing the model with independent demand (analyzed in \citet{S-12}) in Section\ \ref{indiesection}.  We then introduce our distributionally robust model with martingale dependency structure in Section\ \ref{martpresubsec}, and prove the validity of a certain DP formulation (very similar to the robust MDP framework of \citet{Iyengar}).  We present our main results in Section\ \ref{mainsubsec}, which include the explicit solution to our distributionally robust martingale model, the weak convergence, and a comparison between our martingale-demand model and the independent-demand model. The proofs of our main results regarding the explicit solution of the associated DP and several related structural properties are given in Section\ \ref{mainmart1}.  In Section\ \ref{sec-mainproof2} we prove our weak convergence results.  In Section\ \ref{sec-nume}, we conduct a targeted and concise numerical experiment.  We provide a summary of our results, directions for future research, and concluding remarks in Section\ \ref{sec-conclusion}.  We also provide a technical appendix, which contains the proofs of several results used throughout the paper, as well as some additional supporting results, in Section\ \ref{rm-sec-appendix}.

\section{Main results}\label{mainsec}
Before stating our main results, we briefly review the inventory model analyzed in \citet{S-12}, which we refer to as the independent-demand model.
\subsection{Independent-demand model}\label{indiesection}

Consider the following distributionally robust inventory control problem with backlogging, finite time horizon $T$, strictly positive linear backorder per-unit cost $b$, and a per-unit holding cost of 1, where we note that assuming a holding cost equal to 1 is without loss of generality (w.l.o.g.) due to a simple scaling argument.  Let $D_t$ be the demand in period $t$, and $x_t$ be the inventory level at period $t$ after placing a non-negative order, $t \in [1,T]$, where we note that the order must be placed at the start of period $t$ before $D_t$ is known.  For a $t$-dimensional vector $\mathbf{x}$ (potentially of real numbers, r.v.s, functions, etc.), let $\mathbf{x}_{[s]} = (x_1,\ldots,x_s)$, and $\mathbf{x}_{[t_1,t_2]} = (x_{t_1},\ldots,x_{t_2})$.  In that case, we require that $x_t$ is a measurable function of $D_{[t-1]}$, 
$x_t \geq x_{t-1} - D_{t-1}$ for $t \in [2,T]$, and that $x_1$ is a real constant, where $x_0 \in [0,U]$ denotes the initial inventory level and $x_1 \geq x_0$.

For $t \in [1,T]$, the cost incurred in period $t$ equals
$$C_t \stackrel{\Delta}{=} b[D_t-x_t]_+ + [x_t-D_t]_+,$$

where $z_+ \stackrel{\Delta}{=} \max(z,0), z_- \stackrel{\Delta}{=} \max(-z,0), z \in \bbr$.  As a notational convenience, we also define $C(x,d) \stackrel{\Delta}{=} b(d- x)_+ + (x- d)_+$.

As in \citet{Iyengar}, to avoid later possible concerns about measurability, throughout we make the minor technical assumption that all relevant demand distributions have rational support, and all policies $\pi$ considered are restricted to order rational quantities (i.e. not irrational quantities such as $\sqrt{2}$), and the relevant problem primitives $x_0, \mu, U, b$ (corresponding to initial inventory level, mean demand, bound on support of demand, and backorder penalty respectively) also belong to ${\mathcal Q}^+$, with $x_0,\mu \in [0,U]$, where ${\mathcal Q}^+$ is the set of all non-negative rational numbers.  Under this assumption, it may be easily verified (due to the finiteness and non-negativity of all relevant cost functions) that all relevant functions arising in the various DP which we shall later define are sufficiently regular, e.g. measurable, integrable, and finite-valued.  For clarity of exposition, and in light of this assumption, we do not dwell further on related questions of measurability, and instead refer the interested reader to the excellent reference \citet{Puterman14} for further discussion.

We now formally and more precisely define the relevant notion of admissible policy in our setting.  As our analysis will be restricted to the setting in which there is a uniform upper bound on the support of demand (denoted $U$), for simplicity we will incorporate $U$ directly into several of our definitions.  Let $[0,U]_{{\mathcal Q}} \stackrel{\Delta}{=} [0,U] \bigcap {\mathcal Q}$, where ${\mathcal Q}$ is the set of all rational numbers.  Then we define an admissible policy $\pi$ to be a $T$-dimensional vector of functions $\lbrace x^{\pi}_t, t \in [1,T] \rbrace$, s.t. $x^{\pi}_1$ is a constant satisfying $x^{\pi}_1 \geq x^{\pi}_0 \stackrel{\Delta}{=} y^{\pi}_1 \stackrel{\Delta}{=} x_0$; and for $t \in [2,T]$, $x^{\pi}_t$ is a deterministic map from $[0,U]_{{\mathcal Q}}^{t-1}$ to ${\mathcal Q}$, s.t. for $t \in [1,T-1]$ and $\mathbf{d} = (d_1,\ldots,d_t) \in [0,U]_{{\mathcal Q}}^t$, 
$x^{\pi}_{t+1}(\mathbf{d}) \geq y^{\pi}_{t+1}(\mathbf{d}) \stackrel{\Delta}{=} x^{\pi}_t(\mathbf{d}_{[t-1]}) - d_t$, which is equivalent to requiring that a non-negative amount of inventory is ordered in each period.  Sometimes as a notational convenience, we will denote $x^{\pi}_1$ by $x^{\pi}_1(\mathbf{d}_{[0]})$, or more generally as $x^{\pi}_1(q)$ for $q \in [0,U]_{{\mathcal Q}}$.  

Let $\hat{\Pi}$ denote the family of all such admissible policies.  Note that once a particular policy $\pi \in \hat{\Pi}$ is specified, the associated costs $\lbrace C^{\pi}_t, t \in [1,T] \rbrace$ are explicit functions of the vector of demands, and we sometimes make this dependence explicit with the notation $C^{\pi}_t(D_{[t]})$, where $C^{\pi}_t(D_{[t]}) = C(x^{\pi}_t(D_{[t-1]}),D_t)$.

Since we will exclusively consider settings in which the demand in every period belongs to $[0,U]$, and the initial inventory belongs to $[0,U]$, it will be convenient to prove that in such a setting one may w.l.o.g. restrict to policies which always order up to a level in $[0,U]$.  In particular, let $\Pi$ denote that subset of $\hat{\Pi}$ consisting of policies $\pi$ s.t. $x^{\pi}_t(\mathbf{d}_{[t-1]}) \in [0,U]$ for all $t \in [1,T]$ and $\mathbf{d} \in [0,U]_{{\mathcal Q}}^T$ whenever $x_0 \in [0,U]_{{\mathcal Q}}$ (which again is to be assumed throughout).  Then the following lemma provides the desired result, where we defer the relevant proof to the Technical Appendix in Section\ \ref{rm-sec-appendix}.

\begin{lemma}\label{lemma0u}
For every $\hat{\pi} \in \hat{\Pi}$, there exists $\pi \in \Pi$ s.t. for all $\mathbf{d} = (d_1,\ldots,d_T) \in [0,U]_{{\mathcal Q}}^T$ and $t \in [1,T]$, $C^{\pi}_t(\mathbf{d}_{[t]}) \leq C^{\hat{\pi}}_t(\mathbf{d}_{[t]})$.
\end{lemma}

Let us say that a probability measure $Q$ on an appropriate measurable space is rationally supported if it places probability 1 on $t$-dimensional rational vectors for some finite $t$, where we note that all measures which we consider will have this property.  For such a measure $Q$, let $\supp(Q)$ (i.e. the support of $Q$) denote the (possibly infinite) set of rational vectors which $Q$ assigns strictly positive probability, and $D^Q$ denote a r.v. (or vector if $\supp(Q)$ is not one-dimensional) distributed as $Q$.  For $\mu \in [0,U]_{{\mathcal Q}}$, let $\mathfrak{M}(\mu)$ be the collection of all rationally supported probability measures $Q$ s.t. $\supp(Q) \subseteq [0,U]_{{\mathcal Q}}$, and $E[D^Q] = \mu$.  Furthermore, let $\textbf{IND}_T$ denote the collection of all $T$-dimensional product measures s.t. all $T$ marginal distributions belong to $\mathfrak{M}(\mu)$.  When there is no ambiguity, we will suppress the dependence on $T$, simply writing $\textbf{IND}$.  In words, the joint distribution of demand belongs to $\textbf{IND}$ iff the demand is independent across time periods, and the demand in each period has support a subset of $[0,U]_{{\mathcal Q}}$ and mean $\mu$.  Then the particular optimization problem considered in \citet{S-12} Example 5.1.2 is
$\inf_{\pi \in \hat{\Pi}} \sup_{Q \in \textbf{IND}} \bbe_Q[ \sum_{t=1}^T C^{\pi}_t ],$ which by Lemma\ \ref{lemma0u} is equivalent to
\begin{equation}\label{dynamic0-indep}
\inf_{\pi \in \Pi} \sup_{Q \in \textbf{IND}} \bbe_Q[ \sum_{t=1}^T C^{\pi}_t ].
\end{equation}
We note that although the exact problem considered in \citet{S-12} does not explicitly restrict the family of policies and demand distributions to ${\mathcal Q}$, it may be easily verified that the same results go through in this setting (under our assumptions on the rationality of the relevant problem parameters).  In Example 5.1.2 of \citet{S-12}, the author proves that due to certain structural properties, Problem (\ref{dynamic0-indep}) can be reformulated as a DP, as we now describe.  We note that although the family of measures $\textbf{IND}$ does not possess the so-called rectangularity property, here the decomposition follows from the fact that convexity and independence imply that irregardless of the inventory level and past realizations of demand, a worst-case demand distribution always has support $\lbrace 0, U \rbrace$, where the precise distribution is then dictated by the requirement of having expectation $\mu$.  To formally define the relevant DP, we now define a sequence of functions $\lbrace V^t, t \geq 1 \rbrace, \lbrace f^t, t \geq 1 \rbrace$, $\lbrace g^t, t \geq 1 \rbrace$.  In general, such a DP would be phrased in terms of the so-called ``cost-to-go" functions, with the $t$-th such function representing the remaining cost incurred by an optimal policy during periods $T-t+1$ through $T$ (i.e., if there are $t$ periods to go), subject to the given state at time $t$.  
Here and throughout, we use the fact that the backorder and holding costs are the same in every period to simplify the relevant concepts and notations.  In particular, due to this symmetry and the associated self-reducibility which it induces, for all $T \geq 1$, it will suffice to define the aforementioned cost-to-go function for the first time period only.  Indeed, in the following function definitions, $V^T(y,\mu)$ will coincide with the optimal value of Problem (\ref{dynamic0-indep}) when the initial inventory level is $y$ and the associated mean is $\mu$ (we leave the dependence on $U$ implicit), where $f^T,g^T$ have analogous interpretations.  We note that to prevent later potential confusion and ambiguity regarding the order of operations such as supremum and expectation, we explicitly write out certain expectations as relevant sums over appropriate domains.  Then the relevant DP is as follows.

\begin{equation}\label{dynamic1-indep}
\begin{aligned}
&\ f^1(x_1,d_1) \stackrel{\Delta}{=} C(x_1,d_1)\ \ \ \ \ \ \ \ \ \ \ \ \ \ \ ,\ \ \ g^1(x_1,\mu) \stackrel{\Delta}{=} \sup_{Q_1 \in \mathfrak{M}(\mu)} \sum_{q_1 \in [0,U]_{{\mathcal Q}}} f^1(x_1,q_1) Q_1(q_1),\\
&\ V^1(y_1,\mu) \stackrel{\Delta}{=} \inf_{\substack{z_1 \geq y_1\\z_1 \in [0,U]_{{\mathcal Q}}}} g^1(z_1,\mu)\ \ \ \ \ \ ,\ \ \ Q^1(x_1,\mu) \stackrel{\Delta}{=} \argmax_{Q_1 \in \mathfrak{M}(\mu)} \sum_{q_1 \in [0,U]_{{\mathcal Q}}} f^1(x_1,q_1) Q_1(q_1),\\
&\ \Phi^1(y_1,\mu) \stackrel{\Delta}{=} \argmin_{\substack{z_1 \geq y_1\\z_1 \in [0,U]_{{\mathcal Q}}}} g^1(z_1,\mu);
\end{aligned}
\end{equation}
and for $T \geq 2$,
\begin{equation}\label{dynamic2-indep}
\begin{aligned}
&\ f^T(x_T,d_T) \stackrel{\Delta}{=} C(x_T,d_T) +  V^{T-1}(x_T - d_T,\mu)\ \ \ ,\ \ \ g^T(x_T,\mu) \stackrel{\Delta}{=} \sup_{ Q_T \in \mathfrak{M}(\mu) }  \sum_{q_T \in [0,U]_{{\mathcal Q}}} f^T(x_T,q_T) Q_T(q_T),\\
&\ V^T(y_T,\mu) \stackrel{\Delta}{=} \inf_{\substack{z_T \geq y_T\\z_T \in [0, U]_{{\mathcal Q}}}} g^T(z_T,\mu)\ \ \ \ \ \ \ \ \ \ \ \ \ \ \ \ \ \ \ \ ,\ \ \ Q^T(x_T,\mu) \stackrel{\Delta}{=} \argmax_{Q_T \in \mathfrak{M}(\mu)} \sum_{q_T \in [0,U]_{{\mathcal Q}}} f^T(x_T,q_T) Q_T(q_T),\\
&\ \Phi^T(y_T,\mu) \stackrel{\Delta}{=} \argmin_{\substack{z_T \geq y_T\\z_T \in [0,U]_{{\mathcal Q}}}} g^T(z_T,\mu).
\end{aligned}
\end{equation}
We have given the argument of the relevant functions with superscript $T$ the subscript $T$ (e.g. $x_T,d_T$ are the arguments for $f^T$) to prevent any possible ambiguities when expanding the relevant recursions, although this is of course not strictly necessary.  We do note that conceptually, if the length of the time horizon is fixed to $T$, then the functions with superscript $T$ (e.g. $V^T$) actually correspond to the cost-to-go in period 1 (not period $T$), and the corresponding arguments (e.g. $x_T,d_T$) correspond to the inventory level and demand in period 1.  More generally, for $s \in [0,T-1]$, $V^{T-s}$ will correspond to the cost-to-go in period $s + 1$, and the corresponding arguments  (e.g. $x_{T-s},d_{T-s}$) correspond to the inventory level and demand in period $s + 1$.  We have chosen our notational scheme to maximize clarity of exposition, and again note that several of the resulting simplifications do rely on the symmetries / time-homogeneities of the specific problem considered. 
\\\\To illustrate how the relevant functions evolve, note that $V^2(y_2,\mu)$ equals 
$$\inf_{\substack{z_2 \geq y_2\\z_2 \in [0, U]_{{\mathcal Q}}}} \sup_{ Q_2 \in \mathfrak{M}(\mu) }  \sum_{q_2 \in [0,U]_{{\mathcal Q}}} \big( C(z_2,q_2) + V^1(z_2 - q_2, \mu) \big) Q_2(q_2),$$
which itself equals
$$\inf_{\substack{z_2 \geq y_2\\z_2 \in [0, U]_{{\mathcal Q}}}} \sup_{ Q_2 \in \mathfrak{M}(\mu) }  \sum_{q_2 \in [0,U]_{{\mathcal Q}}} \bigg( C(z_2,q_2) 
+ \inf_{\substack{z_1 \geq z_2 - q_2\\z_1 \in [0,U]_{{\mathcal Q}}}} \sup_{Q_1 \in \mathfrak{M}(\mu)} \sum_{q_1 \in [0,U]_{{\mathcal Q}}} C(z_1,q_1) Q_1(q_1)
\bigg) Q_2(q_2).$$
\\\\We now review the results of Lemma 5.1 and Example 5.1.2 of \citet{S-12}, which characterizes the optimal policy, value, and worst-case distribution (against an optimal policy) for Problem\ \ref{dynamic0-indep}, and relates the solution to DP formulation (\ref{dynamic1-indep}) - (\ref{dynamic2-indep}).  Recall that a policy $\pi \in \Pi$ is said to be a \emph{base-stock policy} if there exist constants $\lbrace x^*_t, t \in [1,T] \rbrace$, s.t. $x^{\pi}_t = \max\left\{y^{\pi}_t ,\ x^{*}_t\right\}$ for all $t \in [1,T]$.  Let us define
$$
\begin{aligned}
   & \chi_{\textbf{IND}}(\mu,U,b) \stackrel{\Delta}{=} \begin{cases}
   0& \text{if}\ \mu \leq \frac{U}{b+1},\\
   U& \text{if}\ \mu > \frac{U}{b+1};
   \end{cases}
\end{aligned}
$$
$$
\begin{aligned}
   & \text{Opt}^T_{\textbf{IND}}(\mu,U,b)  \stackrel{\Delta}{=} \begin{cases}
   T b \mu & \text{if}\ \mu \leq \frac{U}{b+1},\\
   T( U - \mu) & \text{if}\ \mu > \frac{U}{b+1}.
   \end{cases}
\end{aligned}
$$

Let ${\mathcal D}_{\mu} \in \mathfrak{M}(\mu)$ be the probability measure s.t. ${\mathcal D}_{\mu}(0) = 1 - \frac{\mu}{U}, {\mathcal D}_{\mu}(U) = \frac{\mu}{U}$.  Then the following is proven in \citet{S-12}.  
\begin{theorem}[\cite{S-12}]\label{thm-4.1}
For all strictly positive $U,b \in {\mathcal Q}^+, T \geq 1$, and $\mu,x_0 \in [0,U]_{{\mathcal Q}}$, Problem\ \ref{dynamic0-indep} always has an optimal base-stock policy $\pi^*$, in which the associated base-stock constants $\lbrace x^{*}_t, t \in [1,T] \rbrace$ satisfy $x^{*}_t = \chi_{\textbf{IND}}(\mu,U,b)$ for all $t \in [1,T]$.  If $x_0 \leq \chi_{\textbf{IND}}(\mu,U,b)$, the optimal value of Problem\ \ref{dynamic0-indep} equals $\text{Opt}^T_{\textbf{IND}}(\mu,U,b)$; and for all $x_0 \in [0,U]$, the product measure s.t. all marginals are distributed according to law ${\mathcal D}_{\mu}$ belongs to $\argmax_{Q \in \textbf{IND}} \bbe_Q[ \sum_{t=1}^T C^{\pi^*}_t ]$.  Furthermore, the DP formulation (\ref{dynamic1-indep}) - (\ref{dynamic2-indep}) can be used to compute these optimal policies, distributions, and values.  In particular, for all $x_0,x,\mu \in [0,U]_{{\mathcal Q}}$, $V^T(x_0,\mu)$ is the optimal value of Problem\ \ref{dynamic0-indep}, ${\mathcal D}_{\mu} \in Q^T(x,\mu)$, and $\max\big(x,\chi_{\textbf{IND}}(\mu,U,b)\big) \in 
\Phi^T(x,\mu)$.
\end{theorem}

\subsection{Martingale-demand model}\label{martpresubsec}
In this subsection, we formally define our distributionally robust martingale-demand model, and establish some preliminary results.  Let $T,b,\Pi,C^{\pi}_t,\mathfrak{M}(\mu)$ be exactly as defined for the independent-demand model in Subsection\ \ref{indiesection}, i.e. the time horizon, backorder and holding costs, set of admissible policies, cost incurred in period $t$ under policy $\pi$, and collection of all probability measures with support $[0,U]_{{\mathcal Q}}$ and mean $\mu$, respectively.
Furthermore, let $\textbf{MAR}_T$ denote the collection of all rationally supported probability measures corresponding to discrete-time martingale sequences $(D_1,\ldots,D_T)$ s.t. for all $t$, the marginal distribution of $D_t$ belongs to $\mathfrak{M}(\mu)$.  In words, the joint distribution of demand belongs to $\textbf{MAR}_T$ iff the sequence of demands is a martingale, and the demand in each period has support a subset of $[0,U]_{{\mathcal Q}}$ and mean $\mu$.  When there is no ambiguity, we will suppress the dependence on $T$, simply writing $\textbf{MAR}$.
Note that by the martingale property, the \emph{conditional} distribution of $D_t$, conditioned on $D_{[t-1]}$, belongs to $\mathfrak{M}(D_{t-1})$.  Then in analogy with Problem\ \ref{dynamic0-indep}, the optimization problem of interest is
$\inf_{\pi \in \hat{\Pi}} \sup_{Q \in \textbf{MAR}} \bbe_Q[ \sum_{t=1}^T C^{\pi}_t ],$ which by Lemma\ \ref{lemma0u} is equivalent to
\begin{equation}\label{dynamic0-mar}
\inf_{\pi \in \Pi} \sup_{Q \in \textbf{MAR}} \bbe_Q[ \sum_{t=1}^T C^{\pi}_t ].
\end{equation}
In analogy with the DP formulation (\ref{dynamic1-indep}) - (\ref{dynamic2-indep}) given for the independent-demand model, we now define an analogous formulation for the martingale-demand setting, and will later prove that this DP formulation indeed corresponds (in an appropriate sense) to Problem\ \ref{dynamic0-mar}.

\begin{equation}\label{dynamic1-mar}
\begin{aligned}
&\ \hat{f}^1(x_1,d_1) \stackrel{\Delta}{=} C(x_1,d_1)\ \ \ \ \ \ \ \ \ \ \ \ \ \ ,\ \ \ \hat{g}^1(x_1,\mu_1) \stackrel{\Delta}{=} \sup_{Q_1 \in \mathfrak{M}(\mu_1)} \sum_{q_1 \in [0,U]_{{\mathcal Q}}} \hat{f}^1(x_1,q_1) Q_1(q_1),\\
&\ \hat{V}^1(y_1,\mu_1) \stackrel{\Delta}{=} \inf_{\substack{z_1 \geq y_1\\z_1 \in [0,U]_{{\mathcal Q}}}} \hat{g}^1(z_1,\mu_1)\ \ \ ,\ \ \ 
\hat{Q}^1(x_1,\mu_1) \stackrel{\Delta}{=} \argmax_{Q_1 \in \mathfrak{M}(\mu_1)} \sum_{q_1 \in [0,U]_{{\mathcal Q}}} \hat{f}^1(x_1,q_1) Q_1(q_1),\\
&\ \hat{\Phi}^1(y_1,\mu_1) \stackrel{\Delta}{=} \argmin_{\substack{z_1 \geq y_1\\z_1 \in [0,U]_{{\mathcal Q}}}} \hat{g}^1(z_1,\mu_1);
\end{aligned}
\end{equation}
and for $T \geq 2$,
\begin{equation}\label{dynamic2-mar}
\begin{aligned}
&\ \hat{f}^T(x_T,d_T) \stackrel{\Delta}{=} C(x_T,d_T) +  \hat{V}^{T-1}(x_T - d_T,d_T)\ \ \ \ \ ,\ \ \ \hat{g}^T(x_T,\mu_T) \stackrel{\Delta}{=} \sup_{ Q_T \in \mathfrak{M}(\mu_T) }  \sum_{q_T \in [0,U]_{{\mathcal Q}}} \hat{f}^T(x_T,q_T) Q_T(q_T),\\
&\ \hat{V}^T(y_T,\mu_T) \stackrel{\Delta}{=} \inf_{\substack{z_T \geq y_T\\z_T \in [0, U]_{{\mathcal Q}}}} \hat{g}^T(z_T,\mu_T)\ \ \ \ \ \ \ \ \ \ \ \ \ \ \ \ \ \ \ \ ,\ \ \ \hat{Q}^T(x_T,\mu_T) \stackrel{\Delta}{=} \argmax_{Q_T \in \mathfrak{M}(\mu_T)} \sum_{q_T \in [0,U]_{{\mathcal Q}}} \hat{f}^T(x_T,q_T) Q_T(q_T),\\
&\ \hat{\Phi}^T(y_T,\mu_T) \stackrel{\Delta}{=} \argmin_{\substack{z_T \geq y_T\\z_T \in [0,U]_{{\mathcal Q}}}} \hat{g}^T(z_T,\mu_T).
\end{aligned}
\end{equation}

Note that in the definition of ${f}^T(x_T,d_T)$ in \eqref{dynamic1-indep}-\eqref{dynamic2-indep}, the second argument in the function $ {V}^{T-1}$ is $\mu$, which is a constant. By contrast, in the definition of $\hat{f}^T(x_T,d_T)$ above, the second argument in the function $ \hat{V}^{T-1}$ is $d_T$, which represents the demand just realized in the past period and is not a constant.

\ \\To again illustrate how the relevant functions evolve, note that $\hat{V}^2(y_2,\mu_2)$ equals 
$$
\inf_{\substack{z_2 \geq y_2\\z_2 \in [0, U]_{{\mathcal Q}}}} \sup_{ Q_2 \in \mathfrak{M}(\mu_2) }  \sum_{q_2 \in [0,U]_{{\mathcal Q}}} \big( C(z_2,q_2) + \hat{V}^1(z_2 - q_2, q_2) \big) Q_2(q_2),$$
which itself equals
$$\inf_{\substack{z_2 \geq y_2\\z_2 \in [0, U]_{{\mathcal Q}}}} \sup_{ Q_2 \in \mathfrak{M}(\mu_2) }  \sum_{q_2 \in [0,U]_{{\mathcal Q}}} \bigg( C(z_2,q_2) 
+ \inf_{\substack{z_1 \geq z_2 - q_2\\z_1 \in [0,U]_{{\mathcal Q}}}} \sup_{Q_1 \in \mathfrak{M}(q_2)} \sum_{q_1 \in [0,U]_{{\mathcal Q}}} C(z_1,q_1) Q_1(q_1)
\bigg) Q_2(q_2).$$

We now prove that the DP\ \ref{dynamic1-mar}\ -\ \ref{dynamic2-mar} indeed corresponds to Problem\ \ref{dynamic0-mar}, and begin with several additional definitions.   Recall that given a probability measure $Q \in \textbf{MAR}_T$, $D^Q = (D^Q_1,\ldots,D^Q_T)$ denotes a random vector distributed as $Q$.  Given $t_1 \leq t_2 \in [1,T]$, we let $Q_{[t_1,t_2]}$ denote the probability measure induced by $Q$ on $(D^Q_{t_1},\ldots,D^Q_{t_2})$, letting $Q_{[t_2]}$ denote $Q_{[1,t_2]}$, and $Q_{t_2}$ denote $Q_{[t_2,t_2]}$.  As a notational convenience, we define $Q_0$ to be the probability measure which assigns $\mu$ probability one.  Thus $\supp(Q_0) = \lbrace \mu \rbrace$.  Also, we sometimes let $\mu_0$ denote $\mu$.  Given two finite vectors $\mathbf{x} \in \bbr^n ,\mathbf{y} \in \bbr^m$, we let $\mathbf{x}:\mathbf{y} \in \bbr^{n + m}$ denote the concatenation of $\mathbf{x}$ and $\mathbf{y}$.  Namely, $(x:y)_k = x_k$ for $k \in [1,n]$, and $(x:y)_k = y _{k - n}$ for $k \in [n+1,n+m]$.  Given $Q \in \textbf{MAR}_T$, $t \in [1,T]$ and $\mathbf{q} \in \supp(Q_{[t-1]})$, let $Q_{t | \mathbf{q}}$ denote the unique probability measure s.t. for all $d \in [0,U]_{{\mathcal Q}}$, $Q_{t | \mathbf{q}}(d) = \frac{Q_{[t]}(\mathbf{q}:d)}{Q_{[t-1]}(\mathbf{q})}$.  
Equivalently, $Q_{t | \mathbf{q}}(d) = \pr\big(D^Q_t = d\ \big|\ D^Q_{[t-1]} = \mathbf{q}\big)$.  Note that by the martingale property, $Q_{t | \mathbf{q}} \in \mathfrak{M}(q_{t-1})$.
\\\\Before precisely stating the desired result, we will need some preliminary results regarding the DP\ \ref{dynamic1-mar}\ -\ \ref{dynamic2-mar}, which assert the non-emptiness of various sets, as well as a certain axiom-of-choice type result asserting the ability to appropriately select elements from various sets.
\begin{lemma}\label{nonempty1}
For all $t \in [1,T]$, and $x,d \in [0,U]_{{\mathcal Q}}$, $\hat{Q}^t(x , d) \neq \emptyset, \hat{\Phi}^t(x, d) \neq \emptyset$.  Furthermore, it is possible to assign to every $t \in [1,T]$, and $x,d \in [0,U]_{{\mathcal Q}}$, a measure $\overline{Q}^t_{x,d} \in \hat{Q}^t(x , d)$; and rational number $\overline{\Phi}^t_{x,d} \in \hat{\Phi}^t(x, d)$.
\end{lemma}

Although we temporarily defer the proof of Lemma\ \ref{nonempty1}, as the lemma will follow from the proof of our main result in which we explicitly construct such $\overline{Q}^t_{x,d}$, $\overline{\Phi}^t_{x,d}$, suppose for now that such an assignment $\overline{Q}^t_{x,d}$, $\overline{\Phi}^t_{x,d}$ has been made in an arbitrary yet fixed manner.  We now use these quantities to define a corresponding policy and distribution, which will correspond to the optimal policy, and worst-case distribution at optimality, for Problem\ \ref{dynamic0-mar}.  Let $\hat{\pi}^* = (\hat{x}^*_1,\ldots,\hat{x}^*_T) \in \Pi$ denote the following policy, which we define inductively.  $\hat{x}^*_1 = \overline{\Phi}^T_{x_0,\mu}$.  Now, suppose we have completely specified $\hat{x}^*_s$ for all $s \in [1,t]$ for some $t \in [1,T-1]$, and thus have also specified $\hat{y}^*_s \stackrel{\Delta}{=} y^{\hat{\pi}^*}_s$ for all $s \in [1,t+1]$.  Then we define $\hat{x}^*_{t+1}$ as follows.  For all $\mathbf{d} \in [0,U]_{{\mathcal Q}}^{t}$, $\hat{x}^*_{t+1}(\mathbf{d}) = \overline{\Phi}^{T-t}_{\hat{y}^*_{t+1}(\mathbf{d}),d_t}$.  It may be easily verified from the definition of $\hat{\Phi}$ that this defines a valid policy belonging to $\Pi$.  Similarly, let $\hat{Q}^* \in \textbf{MAR}_T$ denote the following probability measure, which we define inductively through its conditional distributions.  $\hat{Q}^*_1 = \overline{Q}^T_{x_0,\mu}$.   Now, suppose we have completely specified $\hat{Q}^*_{[t]}$.  Then we specify $\hat{Q}^*_{[t+1]}$ as follows.  For all $\mathbf{d} \in \supp(\hat{Q}^*_{[t]})$, $\hat{Q}^*_{t+1 | \mathbf{d}} = 
\overline{Q}^{T-t}_{x^{\hat{\pi}^*}_{t+1}(\mathbf{d}),d_t}$.  Due to its countable support, it may be easily verified that this uniquely and completely specifies a measure $\hat{Q}^* \in \textbf{MAR}_T$.  As a notational convenience, we let 
$\hat{\textbf{D}}^* = (\hat{D}^*_1,\ldots,\hat{D}^*_T)$ denote a vector distributed as $\hat{Q}^*$.
\\\\Then we have the following theorem, in analogy with Theorem\ \ref{thm-4.1}, whose proof we defer to the Technical Appendix in Section\ \ref{rm-sec-appendix}.
\begin{theorem}\label{thm-4.2}
For all strictly positive $U,b \in {\mathcal Q}^+, T \geq 1$, and $\mu,x_0 \in [0,U]_{{\mathcal Q}}$, $\hat{\pi}^*$ is an optimal policy for Problem\ \ref{dynamic0-mar}, and $\hat{Q}^* \in \argmax_{Q \in \textbf{MAR}} \bbe_Q[ \sum_{t=1}^T C^{\hat{\pi}^*}_t ]$.  Furthermore, 
$\hat{V}^T(x_0,\mu)$ is the optimal value of Problem\ \ref{dynamic0-mar}.  
\end{theorem}
\ \\Note that in light of Lemma\ \ref{nonempty1}, Theorem\ \ref{thm-4.2} implies that computing an optimal policy (and associated worst-case distribution at optimality) for Problem\ \ref{dynamic0-mar} is equivalent to solving the DP\ \ref{dynamic1-mar}\ -\ \ref{dynamic2-mar}, and computing the relevant minimizers and maximizers.  Furthermore, in light of Theorem\ \ref{thm-4.2}, solving Problem\ \ref{dynamic0-mar} is thus reduced from a minimax problem to purely a maximization problem (of the expected performance of the policy $\hat{\pi}^*$), over a restricted family of martingales.  This puts the problem within the general framework of \citet{beigblock2014}, in which the authors develop a general methodology, based on concave envelopes, to compute optimal martingale inequalities.  Indeed, our resolution of the DP\ \ref{dynamic1-mar}\ -\ \ref{dynamic2-mar} can be viewed as an explicit execution of this methodology, which exploits the special structure of $\hat{\pi}^*$.  It is also interesting to note that Theorem\ \ref{thm-4.2} asserts that one can w.l.o.g. restrict to policies and distributions which are Markovian when one considers the three-dimensional state-space consisting of both the post-ordering inventory level and demand in the previous period, as well as a dimension tracking the current time period, even though in principle both the policy and distribution could have exhibited a more complicated form of history dependence.

\subsection{Main results}\label{mainsubsec}
We now state our main results.  We begin by introducing some additional definitions and notations.  Whenever possible, we will suppress dependence on the parameters $\mu,U,b$ for simplicity.  As a notational convenience, let us define all empty products to equal unity, and all empty sums to equal zero.   For an event $A$, let $I(A)$ denote the corresponding indicator.  For $j \in [-1,T]$, let
\begin{align*}
A^T_j &\stackrel{\Delta}{=} \begin{cases}
U \prod_{k = j + 1}^{T - 1} \frac{k}{b + k} & \textrm{if}\  j \leq T - 1,\\
\frac{b + T}{T} U & \textrm{if}\ j = T;
\end{cases}
\end{align*}
and $$B^T_j \stackrel{\Delta}{=} \frac{j}{b + T} A^T_j.$$  Note that $A^T_{-1} = B^T_{-1} = B^T_0 = 0$, $A^T_{T-1} = B^T_T = U$, and both $\lbrace A^T_j, j \in [-1,T] \rbrace$ and $\lbrace B^T_j, j \in [-1,T] \rbrace$ are increasing in $j$, and decreasing in $T$.  For $x,\mu \in [0,U]$, let ${\mathfrak q}^T_{x,\mu}$ be the unique (up to sets of measure zero) probability measure s.t.
\begin{align}\label{last-hope}
\begin{cases}
{\mathfrak q}^T_{x,\mu}(A^T_j) = \frac{A^T_{j+1} - \mu}{A^T_{j+1} - A^T_j}\ ,\ {\mathfrak q}^T_{x,\mu}(A^T_{j+1}) = \frac{\mu - A^T_j}{A^T_{j+1} - A^T_j} & \textrm{if}\ \mu \in (A^T_j, A^T_{j+1}]\ ,\ x \in [0, B^T_{j+1});\\
{\mathfrak q}^T_{x,\mu}(0) = 1 - \frac{\mu}{A^T_k}\ ,\ {\mathfrak q}^T_{x,\mu}(A^T_k) = \frac{\mu}{A^T_k}& \textrm{if}\ \mu \in (A^T_j, A^T_{j+1}]\ ,\ x \in [B^T_{k}, B^T_{k+1})\ ,\ k \geq j + 1;\\
{\mathfrak q}^T_{x,\mu}(0) = 1 - \frac{\mu}{U}\ ,\ {\mathfrak q}^T_{x,\mu}(U) = \frac{\mu}{U}& \textrm{if}\ \mu = 0\ \textrm{or}\ x = U.
\end{cases}
\end{align}
Also, for $j \in [0,T]$, let
$$G^T_j(x,\mu) \stackrel{\Delta}{=} (T - \frac{b+T}{A^T_j} \mu) x + (T - j) b \mu;$$
in which case we define
\begin{align*}
\Gamma^T_{\mu} &\stackrel{\Delta}{=} \begin{cases}
0 & \textrm{if}\ \mu = 0,\\
j + 1 & \textrm{if}\  \mu \in (A^{T+1}_j, A^{T+1}_{j+1}]\ ,\ j \in [-1,T-1];
\end{cases}
\end{align*}
and note that $\mu \in (A^T_{\Gamma^{T-1}_{\mu} - 1}, A^T_{\Gamma^{T-1}_{\mu}}]$ for all $T \geq 2$ and $\mu \in (0,U]$, while $A^T_{\Gamma^{T-1}_0 - 1} = 0$.  We also define
$$
\chi^T_{\textbf{MAR}}(\mu,U,b) \stackrel{\Delta}{=} \beta^T_{\mu} \stackrel{\Delta}{=} B^T_{\Gamma^T_{\mu}}\ \ \ ,\ \ \ \text{Opt}^T_{\textbf{MAR}}(\mu,U,b) \stackrel{\Delta}{=} G^T_{\Gamma^T_{\mu}}(\beta^T_{\mu},\mu).$$

Then an explicit solution to the DP\ \ref{dynamic1-mar}\ -\ \ref{dynamic2-mar}, and thus to Problem\ \ref{dynamic0-mar}, is given as follows.

\begin{theorem}\label{mainmart1}
For all $t \in [1,T]$, and all $x,d \in [0,U]_{{\mathcal Q}}$, $\max\big(x, \chi^t_{\textbf{MAR}}(d,U,b) \big) \in \hat{\Phi}^t(x, d)$, and we may set $\overline{\Phi}^t_{x,d}  = \max\big(x, \chi^t_{\textbf{MAR}}(d,U,b) \big)$.
Thus an optimal policy $\hat{\pi}^*$ for Problem\ \ref{dynamic0-mar} is given as follows.  $\hat{x}^*_1 = \max\big(x_0, \chi^T_{\textbf{MAR}}(\mu,U,b)\big)$.  For $t \in [1,T-1]$, 
$$\hat{x}^*_{t+1}(\mathbf{d}_{[t]}) = \max\big(\hat{x}^*_t(\mathbf{d}_{[t-1]}) - d_t, \chi^{T-t}_{\textbf{MAR}}(d_t,U,b) \big).$$
Similarly, ${\mathfrak q}^t_{x,d} \in \hat{Q}^t(x,d)$, and we may set $\overline{Q}^t_{x,d} = {\mathfrak q}^t_{x,d}$.  Thus a worst-case (at optimality) distribution $\hat{Q}^* \in \argmax_{Q \in \textbf{MAR}} \bbe_Q[ \sum_{t=1}^T C^{\hat{\pi}^*}_t ]$ may be given as follows.  $\hat{Q}^*_1 = {\mathfrak q}^T_{\hat{x}^*_1,\mu}$.  For $t \in [1,T-1]$ and $\mathbf{d} \in \supp(\hat{Q}^*_{[t]})$, 
$$\hat{Q}^*_{t+1 | \mathbf{d}} = {\mathfrak q}^{T-t}_{\hat{x}^*_{t+1}(\mathbf{d}),d_t}.$$
Furthermore, for all $x_0 \leq \beta^T_{\mu}$, the optimal value of Problem\ \ref{dynamic0-mar} equals $\text{Opt}^T_{\textbf{MAR}}(\mu,U,b)$.
\end{theorem}

We note that this optimal solution is not necessarily unique, but do not attempt to characterize the family of all possible solutions here.  Furthermore, as the optimal policy is essentially given in closed form, it may be computed very efficiently.  
\\\indent It is possible to infer many interesting qualitative features of the optimal policy $\hat{\pi}^*$ and worst-case (at optimality) measure $\hat{Q}^*$ from Theorem\ \ref{mainmart1}.  First, we have the following observation, whose proof we defer to the Technical Appendix in Section\ \ref{rm-sec-appendix}.
\begin{observation}.  \label{thereiszero}
For any initial inventory level $x_0 \in [0,U]$, the (optimal) policy $\hat{\pi}^* = (\hat{x}^*_1,\ldots,\hat{x}^*_T)$, and the random vector $\hat{\mathbf{D}}^*$ (distributed as the worst-case measure $\hat{Q}^*$) satisfy the following relation.  For all $t \in [1,T]$, w.p.1, either $\hat{D}^{*}_t \geq \hat{x}^*_t(\hat{D}^*_{[t-1]})$ (i.e. the demand in period $t$ clears the inventory), or $\hat{D}^{*}_t = 0$ (i.e. the demand in period $t$ equals 0, all future demands also equal 0, and the policy-maker is left holding her inventory for the remainder of the horizon).  Furthermore, for all $T \geq 1$, $\mu \in (0,U)$, and $x \in \big[ \chi^{T}_{\textbf{MAR}}(\mu,U,b) , U \big], {\mathfrak q}^T_{x,\mu}$ has support on two points, one of which is zero, and the other of which is at least x.  It also follows that under these assumptions, the first case in the definition of ${\mathfrak q}^T_{x,\mu}$ cannot occur.
\end{observation}
We next further formalize this observation under the additional assumption that one's initial inventory is sufficiently small, i.e.\ $x_0 \leq \chi^T_{\textbf{MAR}}(\mu,U,b)$.  In this case the stochastic inventory and demand process induced by the policy $\hat{\pi}^*$, and sequence of demands distributed as $\hat{Q}^*$, takes a simpler form which we now describe explicitly.  We also relate this simpler form, which results from the fact that the (conditional) distribution of demand always has 2-point support with positive probability at 0 (as noted in Observation\ \ref{thereiszero}), to the concept of inventory obsolescence (see Interpretation\ \ref{obsoleteinterp}).  Let $\Lambda^T \stackrel{\Delta}{=} \max(T -  \Gamma^T_{\mu},1), D_0^T\stackrel{\Delta}{=} \mu, X_0^T\stackrel{\Delta}{=} \frac{\Gamma^T_\mu}{b+T+1} \mu$, and for $t \in [1, \Lambda^T]$,
\begin{equation}\label{rm-eq-def}
D^T_t \stackrel{\Delta}{=} A^{T + 1 - t}_{\min(\Gamma^T_{\mu},T-1)}\ \ \ ,\ \ \ X^T_t \stackrel{\Delta}{=} B^{T + 1 - t}_{\Gamma^T_{\mu}}.
\end{equation}
Note that $\Lambda^T$ represents the first time that $D^T_t$ reaches $U$.

We now formally relate the corresponding inventory and demand dynamics to an appropriate Markov chain.  Consider the following discrete time Markov chain $\lbrace M^T(t) \stackrel{\Delta}{=} \big( {\mathcal X}^{T}_t, {\mathcal D}^{T}_t \big), t \geq 1 \rbrace$, with randomized initial conditions.
\begin{align*}
\big(  {\mathcal X}^{T}_1, {\mathcal D}^{T}_1 \big)
& = \begin{cases}
\big( X^T_1, D^T_1 \big) & \textrm{w.p.}\  \frac{\mu}{D^T_1},\\
\big( X^T_1, 0 \big) & \textrm{w.p.}\  1 - \frac{\mu}{D^T_1};
\end{cases}
\end{align*}
for $t \in [2, \Lambda^T - 1]$,
\begin{align*}
\big( {\mathcal X}^{T}_t, {\mathcal D}^{T}_t \big) & = \begin{cases}
\big( {\mathcal X}^{T}_{t-1}, 0 \big) & \textrm{w.p.}\ 1,\ \textrm{if}\ {\mathcal D}^{T}_{t-1} = 0,\\
\big( X^T_t, D^T_t \big) & \textrm{w.p.}\ \frac{D^T_{t-1}}{D^T_t},\ \textrm{if}\ {\mathcal D}^{T}_{t-1} \neq 0,\\
\big( X^T_t, 0 \big) & \textrm{w.p.}\ 1 - \frac{D^T_{t-1}}{D^T_t},\ \textrm{if}\ {\mathcal D}^{T}_{t-1} \neq 0;
\end{cases}
\end{align*}
for $t = \Lambda^T$ (if $\Lambda^T \neq 1$),
\begin{align*}
\big( {\mathcal X}^{T}_t, {\mathcal D}^{T}_t  \big)
& = \begin{cases}
\big( {\mathcal X}^{T}_{t-1}, 0 \big) & \textrm{w.p.}\ 1,\ \textrm{if}\ {\mathcal D}^{T}_{t-1} = 0,\\
\big( X^T_t, U \big) & \textrm{w.p.}\ \frac{D^T_{t-1}}{U},\ \textrm{if}\ {\mathcal D}^{T}_{t-1} \neq 0,\\
\big( X^T_t, 0 \big) & \textrm{w.p.}\ 1 - \frac{D^T_{t-1}}{U},\ \textrm{if}\ {\mathcal D}^{T}_{t-1} \neq 0;
\end{cases}
\end{align*}
for $t \geq \min\big(\Lambda^T + 1, T\big)$,
\begin{align*}
\big( {\mathcal X}^{T}_t, {\mathcal D}^{T}_t  \big) & = \begin{cases}
\big( {\mathcal X}^{T}_{\Lambda^T}, 0 \big)& \textrm{w.p.}\ 1,\ \textrm{if}\ {\mathcal D}^{T}_{\Lambda^T} = 0;\\
\big( U, U \big)& \textrm{w.p.}\ 1,\ \textrm{if}\ {\mathcal D}^{T}_{\Lambda^T} = U.\\
\end{cases}
\end{align*}

The state transitions of the above discrete time Markov chain $\left\{M^T(t), t\geq 1\right\}$ can be  described as follows. In period $t=1$, the inventory level is set to be $X_1^T$ and the demand is realized as either $D_1^T$ or zero with certain probabilities (note that the probabilities depend on $\mu$). Going forward, if the demand realized in the past period $t-1$ is zero, then the demand in period $t$ remains zero due to the martingale property, and the inventory level does not change from period $t-1$ to period $t$.  Between periods $2$ and $\Lambda^T$, if the demand has not yet hit zero, the post-ordering inventory level increases (in period t) from $X_{t-1}^T$ to $X_t^T,$ and the demand (in period t) is realized as either $D_t^T$ or zero with certain probabilities (the probabilities depend on the demand realized in the past period). In particular, $D_{t}^T$ reaches $U$ when $t=\Lambda^T$, where $\lbrace D^T_t, t \in [1, \Lambda^T] \rbrace$ is an increasing sequence.  Starting from $t=\Lambda^T+1$, the demand is either zero or $U$ in all remaining periods. If the demand is zero, then the inventory level does not change, stuck at its value before the first zero demand.  If the demand is $U$, then the post-ordering inventory level is set to $U$ in all remaining periods.

\begin{corollary}\label{mainmartcorr1}
If $x_0 \in [0, \chi^T_{\textbf{MAR}}(\mu,U,b)]$, then one can construct $\hat{\mathbf{D}}^*$ on a common probability space with $M^T$ s.t.
$\lbrace \big( \hat{x}^*_t(\hat{D}^*_{[t-1]}), \hat{D}^*_t \big), t \in [1,T] \rbrace$ equals $\lbrace \big( {\mathcal X}^T_t, {\mathcal D}^T_t \big), t \in [1,T] \rbrace$ w.p.1.
Furthermore, $\lbrace D^T_t , t \in [1,\Lambda^T] \rbrace$ and $\lbrace X^T_t, t \in [1,\Lambda^T] \rbrace$ are both monotone increasing, with $X^T_t \leq D^T_t$ for all $t \in [1,\Lambda^T]$ and $D^T_{\Lambda^T} = U$.
\end{corollary}

We now give an alternate description of the dynamics described in Corollary\ \ref{mainmartcorr1}, explicitly describing the random amount of time until the demand is either $0$ or $U$.  Let $Z^T$ denote the r.v., with support on the integers belonging to $[1, \Lambda^T]$, whose corresponding probability measure ${\mathcal Z}^T$ satisfies
\begin{align*}
{\mathcal Z}^T(t) & = \begin{cases}
\big( 1 - \frac{D^T_{t-1}}{D^T_t} \big) \frac{\mu}{D^T_{t-1}}& \textrm{if}\ t \in [1, \Lambda^T - 1];\\
\frac{\mu}{D^T_{t - 1}} & \textrm{if}\ t = \Lambda^T.
\end{cases}
\end{align*}
Also, let $Y^T$ denote a r.v., with support on $\lbrace 0,U \rbrace$, independent of $Z^T$, whose corresponding probability measure ${\mathcal Y}^T$ satisfies
\begin{align*}
{\mathcal Y}^T(x) & = \begin{cases}
\frac{D^T_{\Lambda^T - 1}}{U} & \textrm{if}\ x = U;\\
1 - \frac{D^T_{\Lambda^T - 1}}{U} & \textrm{if}\ x = 0.
\end{cases}
\end{align*}

\begin{corollary}\label{mainmartcorr2}
Under the same assumptions, and on the same probability space, as described in Corollary \ref{mainmartcorr1}, one can also construct $Z^T, Y^T$ s.t. all of the following hold w.p.1.  $\lbrace \big( {\mathcal X}^T_t, {\mathcal D}^T_t \big), t \in [1,Z^T - 1] \rbrace$ equals $\lbrace \big( X^T_t, D^T_t \big), t \in [1,Z^T - 1] \rbrace$.
If $Z^T \leq \Lambda^T - 1$, then $\big( {\mathcal X}^T_t, {\mathcal D}^T_t \big) =  \big( X^T_{Z^T}, 0 \big)$ for all $t \in [Z^T, T]$.  If $Z^T = \Lambda^T$, and $Y^T = U$, then $\big( {\mathcal X}^T_{\Lambda^T}, {\mathcal D}^T_{\Lambda^T} \big) =  \big( X^T_{\Lambda^T}, U \big)$, and $\big( {\mathcal X}^T_t, {\mathcal D}^T_t \big) = \big(U,U\big)$ for all $t \in [\Lambda^T + 1, T]$.  If $Z^T = \Lambda^T$, and $Y^T = 0$, then $\big( {\mathcal X}^T_t, {\mathcal D}^T_t \big) = \big(X^T_{\Lambda^T}, 0\big)$ for all $t \in [\Lambda^T, T]$.  Furthermore, $\hat{D}^*_t = D^T_t I(Z^T > t), t \in [1, \Lambda^T - 1]$; and $\hat{D}^*_t = U I(Z^T = \Lambda^T, Y^T = U), t \in [\Lambda^T, T]$.
\end{corollary}

Note that Corollary\ \ref{mainmartcorr2} implies that $\hat{\mathbf{D}}^*$ increases along a sequence of pre-determined constants (i.e. $D^T_t$) until a random time at which it either drops to 0 or jumps to $U$, and stays there.
Similar martingales have played an important role in studying extremal martingales in other contexts, e.g. certain worst-case distributions which arise in the study of the celebrated prophet inequalities (cf. \citet{HK83}).  In light of our results, the appearance of such martingales has an appealing interpretation.

\begin{interpretation}.\label{obsoleteinterp}
Under our model of demand uncertainty, a robust inventory manager should be most concerned about the possibility that her product becomes obsolete, i.e. demand for her product drops to zero at a random time.
\end{interpretation}

We again note that the scenario of product obsolescence has been recognized as practically relevant and analyzed extensively in the inventory literature.  In contrast to much previous literature on this topic, here our model \emph{did not assume that products would become obsolete} - instead the risks of obsolescence arose naturally from our worst-case analysis.  Another interesting feature of our analysis is that even in a worst-case setting, the system may enter a state s.t. for all subsequent time periods, demand becomes perfectly predictable, and can be perfectly matched (at the level $U$).  This is again intuitively appealing, as it coincides with the notion that one's product has reached such a level of popularity that demand predictably attains its maximum.  Furthermore, it suggests a certain dichotomy - namely, that either one's product eventually becomes obsolete, or eventually reaches a level of popularity for which demand predictably attains its maximum.  

We now present an intuitive description of these dynamics, explained as a game between the ``inventory manager" (selector of the policy) and ``nature" (selector of the worst-case martingale).  Nature reasons that, if her demand in some period $t$ is ever zero, the martingale property ensures that she will order zero in all subsequent periods.  This will leave the inventory manager stuck holding all the inventory which she held at the start of period $t$ for the entire remainder of the time horizon, which will incur a large cost.  Such a situation is thus desirable from nature's perspective, and one might think that to maximize the probability of this happening, nature will thus maximize the probability that the demand is 0, which would be achieved by putting all probability on $0$ and $U$.  However, as Corollaries\ \ref{mainmartcorr1}\ -\ \ref{mainmartcorr2} show, this is not the case.  In fact, the adversary does not put any probability mass at $U$ until time $\Lambda^T$, i.e. the last time period in which there are any meaningful dynamics.
The reason is that in the martingale-demand model, there is an additional hidden cost for nature associated with putting probability on or near $U$.  In particular, if a demand of $U$ ever occurs, then the martingale property ensures that all future demands must be $U$ as well.  However, it is ``easy" for an inventory manager to perform well against an adversary that always orders $U$ - in particular, she can simply order up to $U$ in every period, incurring zero cost.  In other words, the aforementioned hidden cost to nature comes in the form of a ``loss of randomness", making nature perfectly predictable going forwards.  Corollaries\ \ref{mainmartcorr1}\ -\ \ref{mainmartcorr2} indicate that at optimality, this tradeoff manifests by having nature always put some probability at 0, and some probability on a different quantity $D^T_t$, which increases monotonically to $U$ as $t \uparrow \Lambda^T$.  This ``ramping up" can be explained by observing that although higher values for this second point of the support result in a greater ``loss of randomness", the cost of such a loss becomes smaller over time, as the associated window of time during which the adversary is perfectly predictable shrinks accordingly.  
\\\indent Note that Corollary\ \ref{mainmartcorr1} implies that (in a worst-case setting) the inventory manager continues to ramp up production until either the demand drops to zero, or time $\Lambda^T$ is reached, where the precise manner in which production is ramped up results from an attempt to balance the cost of potential obsoletion (i.e. being stuck holding one's inventory) with the cost of under-ordering.  To better understand this ramping up, and more generally the dynamics of an optimal policy against a worst-case distribution for Problem\ \ref{dynamic0-mar}, we now prove that the Markov chain $M^T$, properly rescaled, converges weakly (as $T \rightarrow \infty$) to a simple and intuitive limiting process.   For a review of the formal definition of weak convergence, and the relevant topological spaces and metrics, we refer the reader to the excellent texts \citet{Billingsley} and \citet{Whitt}.

Let us define
\[
\gamma \stackrel{\Delta}{=} \frac{\mu}{U},\ \ \ \ \ \ \ \Lambda^{\infty} \stackrel{\Delta}{=}
  1 - \gamma^{\frac{1}{b}}.
\]
Let $Z^{\infty}$ denote the mixed (i.e. both continuous and discrete components) r.v., with continuous support on $[0, \Lambda^{\infty} \big)$, and discrete support on the singleton $\Lambda^{\infty}$, whose corresponding probability measure ${\mathcal Z}^{\infty}$ satisfies
\begin{align*}
{\mathcal Z}^{\infty}(\alpha) & = \begin{cases}
 b(1 - \alpha)^{b-1} d\alpha
& \textrm{if}\ \alpha \in [0, \Lambda^{\infty}\big);\\
\gamma  & \textrm{if}\ \alpha = \Lambda^{\infty}.
\end{cases}
\end{align*}
For $\alpha \in \left[0, \Lambda^{\infty} \right]$, where $\alpha$ is again a continuous parameter, let us also define
$$D^{\infty}_\alpha \stackrel{\Delta}{=} \mu (1 - \alpha)^{- b} \ \ \ ,\ \ \
X^{\infty}_\alpha \stackrel{\Delta}{=} \mu \gamma^{\frac{1}{b}}(1 - \alpha)^{-(b+1)}.$$
Note that $D^{\infty}_{\Lambda^{\infty}}= X^{\infty}_{\Lambda^{\infty}}= U$. We now define an appropriate limiting process.  Let ${\mathcal M}^{\infty}(\alpha)_{0 \leq \alpha \leq 1}$ denote the following two-dimensional process, constructed on the same probability space as $Z^{\infty}$.
\begin{align*}
{\mathcal M}^{\infty}(\alpha) & = \begin{cases}
\big( X^{\infty}_\alpha, D^{\infty}_\alpha \big) & \textrm{if}\ \alpha \in [0, Z^{\infty});\\
\big( X^{\infty}_{Z^{\infty}} , 0 \big) & \textrm{if}\ \alpha \geq Z^{\infty}\ \textrm{and}\ Z^{\infty} <  \Lambda^{\infty};\\
\big( U , U \big) & \textrm{if}\ \alpha \geq Z^{\infty}\ \textrm{and}\ Z^{\infty} =  \Lambda^{\infty}.\\
\end{cases}
\end{align*}

Then we have the following weak convergence result.  For $\alpha \in [0,1]$, let ${\mathcal M}^T(\alpha) \stackrel{\Delta}{=} M^T\big(\max(\lfloor \alpha T \rfloor,1)\big)$. 

\begin{theorem}\label{weakconvergethm1}
For all $U,b > 0$ and $\mu \in (0,U)$, the sequence of stochastic processes $\lbrace {\mathcal M}^T(\alpha)_{0 \leq \alpha \leq 1}, T \geq 1 \rbrace$ converges weakly to the process ${\mathcal M}^{\infty}(\alpha)_{0 \leq \alpha \leq 1}$ on the space $D\left([0,1],\ \bbr^2\right)$ under the $J_1$ topology.
\end{theorem}

We note that Theorem\ \ref{weakconvergethm1} reveals several qualitative insights into the corresponding dynamics.  For example, it provides simple analytical formulas for how (in a worst-case setting) the production and demand both ramp up over time, and demonstrates that the limiting probability that obsolescence \emph{ever} occurs (again in such a worst-case setting) equals $1 - \gamma$.  
\\\indent A different lens through which one can gain further insight into Theorem\ \ref{mainmart1} follows by observing that for certain values of $b$ and $\mu$, several quantities appearing in our main results simplify considerably.
\begin{observation}.
When $b = 1$ and $\mu = \frac{k U}{T+1}$ for some $k \in [1,T]$, $A^T_j = \frac{U (j + 1)}{T}$ for $j \in [-1,T]$, and $\Lambda^T = T + 1 - k$.  For $t \in [1,\Lambda^T]$, $D^T_t = \frac{k U}{T + 1 - t}$, 
$$X^T_t = \big( (1 - \frac{t}{T+2})(1 - \frac{t}{T+1}) \big)^{-1} \frac{k}{T+2} \frac{k-1}{T+1} U\ \ ,\ \ \text{Opt}^T_{\textbf{MAR}}(\mu,U,b) = \frac{k}{T} (1 - \frac{k}{T+1}) U T.$$
Furthermore,
\begin{align*}
{\mathcal Z}^T(t) & = \begin{cases}
\frac{1}{T+1}& \textrm{if}\ t \in [1, T - k];\\
\frac{k+1}{T+1} & \textrm{if}\ t = T + 1 - k;
\end{cases}
\end{align*}
and 
\begin{align*}
{\mathcal Y}^T(x) & = \begin{cases}
1 - \frac{1}{k + 1} & \textrm{if}\ x = U;\\
\frac{1}{k+1} & \textrm{if}\ x = 0.
\end{cases}
\end{align*}
Namely, $Z^T$ is equally likely to take any value in the interval $[1, T - k]$.
\end{observation}
Of course, another avenue for more concretely understanding Theorem\ \ref{mainmart1} is to explicitly evaluate all relevant quantities for small values of $T$, as we now do for the case $T = 2$.

\begin{observation}.\label{explicitt2}
\end{observation}
\ \\\\
\begin{tabular}{l l l l l l}
$A^1_{-1} = 0,$ & $A^1_0 = U,$ & $A^1_1 = (b + 1) U,$ & $B^1_{-1} = 0,$ & $B^1_0 = 0,$ & $B^1_1 = U,$\\
$A^2_{-1} = 0,$ & $A^2_0 = \frac{U}{b+1},$ & $A^2_1 = U,$ & $A^2_2 = (1 + \frac{b}{2}) U,$ & $B^2_{-1} = 0,$ & $B^2_0 = 0,$\\
$B^2_1 = \frac{U}{b + 2},$ & $B^2_2 = U,$ & $A^3_{-1} = 0,$ & $A^3_0 = \frac{2 U}{(b+1)(b+2)},$ & $A^3_1 = \frac{2 U}{b+2},$ & $A^3_2 = U,$\\ 
$A^3_3 = (1 + \frac{b}{3})U,$ & $B^3_{-1} = 0,$ & $B^3_0 = 0,$ & $B^3_1 = \frac{2 U}{(b+2)(b+3)},$ &  $B^3_2 = \frac{2 U}{b+3},$ & $B^3_3 = U$.
\end{tabular}

\begin{align*}
\begin{cases}
{\mathfrak q}^2_{x,\mu}(\frac{U}{b+1}) = \frac{U - \mu}{U(1 - \frac{1}{b+1})}
\ ,\ {\mathfrak q}^2_{x,\mu}(U) = \frac{\mu - \frac{U}{b+1}}{U(1 - \frac{1}{b+1})} & \textrm{if}\ \mu \in (\frac{U}{b+1}, U]\ ,\ x \in [0, \frac{U}{b+2});\\
{\mathfrak q}^2_{x,\mu}(0) = 1 - (b+1)\frac{\mu}{U}\ ,\ {\mathfrak q}^2_{x,\mu}(\frac{U}{b+1}) = (b+1)\frac{\mu}{U} & \textrm{if}\ \mu \in (0 , \frac{U}{b+1}]\ ,\ x \in [0, \frac{U}{b+2});\\
{\mathfrak q}^2_{x,\mu}(0) = 1 - \frac{\mu}{U}\ ,\ {\mathfrak q}^2_{x,\mu}(U) = \frac{\mu}{U} & \textrm{if}\ \mu \in (0 , \frac{U}{b+1}]\ ,\ x \in [\frac{U}{b+2}, U);\\
{\mathfrak q}^2_{x,\mu}(0) = 1 - \frac{\mu}{U}\ ,\ {\mathfrak q}^2_{x,\mu}(U) = \frac{\mu}{U} & \textrm{if}\ \mu \in (\frac{U}{b+1} , U]\ ,\ x \in [\frac{U}{b+2}, U);\\
{\mathfrak q}^2_{x,\mu}(0) = 1 - \frac{\mu}{U}\ ,\ {\mathfrak q}^2_{x,\mu}(U) = \frac{\mu}{U}& \textrm{if}\ \mu = 0\ \textrm{or}\ x = U.
\end{cases}
\end{align*}
$$G^2_0(x,\mu) = \big(2 - (b+1)(b+2)\frac{\mu}{U}\big)x + 2 b \mu\ \ \ ,\ \ \ 
G^2_1(x,\mu) = \big(2 - (b+2) \frac{\mu}{U}\big)x + b \mu\ \ \ ,\ \ \ 
G^2_2(x,\mu) = 2(1 - \frac{\mu}{U})x.$$
\[
\begin{aligned}
\Gamma^2_{\mu} &= \begin{cases}
0 & \textrm{if}\  \mu \in [0,  \frac{2 U}{(b+1)(b+2)}],\\
1 & \textrm{if}\  \mu \in (\frac{2 U}{(b+1)(b+2)}, \frac{2 U}{b+2}],\\
2 & \textrm{if}\  \mu \in (\frac{2 U}{b+2}, U].
\end{cases}
\end{aligned}
\begin{aligned}
\Lambda^2 &= \begin{cases}
2 & \textrm{if}\  \mu \in [0,  \frac{2 U}{(b+1)(b+2)}],\\
1 & \textrm{if}\  \mu \in (\frac{2 U}{(b+1)(b+2)}, U].
\end{cases}
\end{aligned}
\]
\[
\begin{aligned}
\chi^2_{\textbf{MAR}} &= \begin{cases}
0 & \textrm{if}\  \mu \in [0,  \frac{2 U}{(b+1)(b+2)}],\\
\frac{U}{b+2} & \textrm{if}\  \mu \in (\frac{2 U}{(b+1)(b+2)}, \frac{2 U}{b+2}],\\
U & \textrm{if}\  \mu \in (\frac{2 U}{b+2}, U].
\end{cases}
\end{aligned}
\begin{aligned}
\text{Opt}^2_{\textbf{MAR}} &= \begin{cases}
2 b \mu & \textrm{if}\  \mu \in [0,  \frac{2 U}{(b+1)(b+2)}],\\
 \frac{2 U}{b+2} + (b - 1) \mu & \textrm{if}\  \mu \in (\frac{2 U}{(b+1)(b+2)}, \frac{2 U}{b+2}],\\
2 (U - \mu) & \textrm{if}\  \mu \in (\frac{2 U}{b+2}, U].
\end{cases}
\end{aligned}
\]
\[
\begin{aligned}
D^2_1 &= \begin{cases}
\frac{U}{b+1} & \textrm{if}\  \mu \in [0,  \frac{2 U}{(b+1)(b+2)}],\\
U & \textrm{if}\  \mu \in (\frac{2 U}{(b+1)(b+2)}, U].
\end{cases}
\end{aligned}
\begin{aligned}
D^2_2 &= \begin{cases}
U & \textrm{if}\  \mu \in [0,  \frac{2 U}{(b+1)(b+2)}].
\end{cases}
\end{aligned}
\]

\[
\begin{aligned}
X^2_1 &= \begin{cases}
0 & \textrm{if}\  \mu \in [0,  \frac{2 U}{(b+1)(b+2)}],\\
\frac{U}{b + 2} & \textrm{if}\  \mu \in (\frac{2 U}{(b+1)(b+2)}, \frac{2 U}{b+2}],\\
U & \textrm{if}\  \mu \in (\frac{2 U}{b+2}, U].
\end{cases}
\end{aligned}
\begin{aligned}
X^2_2 &= \begin{cases}
0 & \textrm{if}\  \mu \in [0,  \frac{2 U}{(b+1)(b+2)}].
\end{cases}
\end{aligned}
\]

\[
\begin{aligned}
\begin{cases}
{\mathcal Z}^2(1) = 1 - (b+1)\frac{\mu}{U}\ ,\ {\mathcal Z}^2(2) = (b+1)\frac{\mu}{U} & \textrm{if}\ \mu \in [0,  \frac{2 U}{(b+1)(b+2)}];\\
{\mathcal Z}^2(1) = 1 & \textrm{if}\ \mu \in (\frac{2 U}{(b+1)(b+2)}, U].
\end{cases}
\end{aligned}
\]
\[
\begin{aligned}
\begin{cases}
{\mathcal Y}^2(U) = \frac{1}{b+1}\ ,\ {\mathcal Y}^2(0) = 1 - \frac{1}{b+1} & \textrm{if}\ \mu \in [0,  \frac{2 U}{(b+1)(b+2)}];\\
{\mathcal Y}^2(U) = \frac{\mu}{U}\ ,\ {\mathcal Y}^2(0) = 1 - \frac{\mu}{U} & \textrm{if}\ \mu \in (\frac{2 U}{(b+1)(b+2)}, U].
\end{cases}
\end{aligned}
\]

To further help interpret Observation\ \ref{explicitt2} and the explicit forms for the various quantities, we show (in Section\ \ref{furtherobssec} of the Technical Appendix)
how (for the case $T = 2$) one may derive the same quantities from a certain ``heuristic relaxation" of our original problem.  In doing so, we see precisely where and why certain quantities arise, and we defer the relevant formulation and discussion to Section\ \ref{furtherobssec} of the Technical Appendix.  We note that Observation\ \ref{explicitt2} also brings another interesting feature of Theorem\ \ref{mainmart1} to light, namely the fact that $\chi^T_{\textbf{MAR}}(\mu,U,b)$ is not monotone increasing in $b$.  This is surprising, since one would expect that as the backlogging penalty increases, one would wish to stock higher inventory levels.  We also note that for fixed $T,\mu,U$, $\chi^T_{\textbf{IND}}(\mu,U,b)$ is monotone increasing in $b$, i.e. such a non-monotonicity does not manifest in the independent-demand model.  Furthermore, Theorem\ \ref{weakconvergethm1}, and the monotonicity (in $b$) of $X^{\infty}_0 = \mu \gamma^{\frac{1}{b}}$, implies that this non-monotonicity does not manifest ``in the limit".

We conclude by presenting several comparative results between the independent-demand and martingale-demand models.  First, we prove that the expected cost incurred by a minimax optimal policy is highest under the independent-demand model, among all possible dependency structures, whenever $x_0 = 0$.  In particular, in this case the expected cost incurred under the independent-demand model is \emph{at least} the expected cost incurred under the martingale-demand model. 
Let $\textbf{GEN}_T$ denote the family of all $T$-dimensional measures $Q$ s.t. $Q_t \in \mathfrak{M}(\mu)$ for all $t$, and $\text{Opt}^T_{\textbf{GEN}}(\mu,U,b)\stackrel{\Delta}{=}\ \inf_{\pi \in \Pi} \sup_{Q \in \textbf{GEN}_T} \bbe_Q[ \sum_{t=1}^T C^{\pi}_t ]$ under the initial condition $x_0 = 0$.  Then we prove the following comparative result, and defer the proof to the Technical Appendix in Section\ \ref{rm-sec-appendix}.

\begin{theorem}\label{compare1theorem}
For all strictly positive $U,b \in {\mathcal Q}^+, T \geq 1$, and $\mu \in [0,U]_{{\mathcal Q}}$, 
$\text{Opt}^T_{\textbf{GEN}}(\mu,U,b) = \text{Opt}^T_{\textbf{IND}}(\mu,U,b)$. 
\end{theorem}
Although Theorem\ \ref{compare1theorem} implies that $\frac{\text{Opt}^T_{\textbf{MAR}}(\mu,U,b)}{\text{Opt}^T_{\textbf{IND}}(\mu,U,b)} \leq 1$, it is interesting to further understand the behavior of this ratio, in the spirit of the so-called price of correlations introduced in \citet{ADSY}.  Indeed, we prove the following, which provides a simple analytical form for this ratio as $T \rightarrow \infty$, again deferring the proof to the Technical Appendix in Section\ \ref{rm-sec-appendix}.

\begin{theorem}\label{compare2theorem}
For all strictly positive $U,b \in {\mathcal Q}^+$, and $\mu \in (0,U)_{{\mathcal Q}}$,  
$$
\begin{aligned}
   & \lim_{T \rightarrow \infty} \frac{\text{Opt}^T_{\textbf{MAR}}(\mu,U,b)}{\text{Opt}^T_{\textbf{IND}}(\mu,U,b)}  = \begin{cases}
   1 - \gamma^{\frac{1}{b}} & \text{if}\ \mu \leq \frac{U}{b+1},\\
   \frac{(1 - \gamma^{\frac{1}{b}}) b \mu}{U - \mu} & \text{if}\ \mu > \frac{U}{b+1}.
   \end{cases}
\end{aligned}
$$
\end{theorem}

We also note the following corollary, which follows immediately from Theorem\ \ref{compare2theorem}.
\begin{corollary}\label{compare2corr}
Suppose $U \in {\mathcal Q}^+, \mu = \frac{U}{2}$, and $b = 1$.  Then
$$\lim_{T \rightarrow \infty} \frac{\text{Opt}^T_{\textbf{MAR}}(\mu,U,b)}{\text{Opt}^T_{\textbf{IND}}(\mu,U,b)} = \frac{1}{2}.$$
\end{corollary}

One setting which arises in various inventory applications is that in which $b >> h$, which motivates considering what our results say about the regime in which $b \rightarrow \infty$.  In particular, the following may be easily verified by combining our main result Theorem\ \ref{mainmart1} with a straightforward calculation.
\begin{corollary}\label{largeb1}
For all strictly positive $U \in {\mathcal Q}^+, \mu \in [0,U]_{{\mathcal Q}}$, $T \geq 1$, and $b \in {\mathcal Q}^+$ s.t. $b > (\frac{U}{\mu} - 1) T$, it holds that $\chi^T_{\textbf{MAR}}(\mu,U,b) = \chi_{\textbf{IND}}(\mu,U,b) = U$, and $\frac{\text{Opt}^T_{\textbf{MAR}}(\mu,U,b)}{\text{Opt}^T_{\textbf{IND}}(\mu,U,b)} = 1$.  
\end{corollary}
Intuitively, Corollary\ \ref{largeb1} indicates that with all other parameters fixed, if $b$ is sufficiently large, then in both scenarios the inventory manager orders up to $U$ in the first period out of fear of incurring large backlogging costs.  Furthermore, although the demand dynamics are different in the two models in this case, they both result in the same minimax expected cost.  Interestingly, when one considers the case of large $b$ in the asymptotic setting of Theorem\ \ref{compare2theorem} (in which $T$ is also very large), a different behavior emerges.  Namely, by letting $b \rightarrow \infty$ in Theorem\ \ref{compare2theorem}, one may conclude the following after a straightforward calculation.
\begin{corollary}\label{largeb2}
For all strictly positive $U \in {\mathcal Q}^+$, and $\mu \in (0,U)_{{\mathcal Q}}$,  
$$
\lim_{b \rightarrow \infty} \lim_{T \rightarrow \infty} \frac{\text{Opt}^T_{\textbf{MAR}}(\mu,U,b)}{\text{Opt}^T_{\textbf{IND}}(\mu,U,b)}  = \frac{\log(\gamma^{-1})}{\gamma^{-1} - 1}.$$
\end{corollary}
Combining Corollaries\ \ref{largeb1} and \ref{largeb2}, we conclude that in this case an interchange of limits does not hold.  Namely, 
as Corollary\ \ref{largeb1} implies that $\lim_{T \rightarrow \infty} \lim_{b \rightarrow \infty} \frac{\text{Opt}^T_{\textbf{MAR}}(\mu,U,b)}{\text{Opt}^T_{\textbf{IND}}(\mu,U,b)} = 1$, here the order in which one lets $b,T$ get large makes a fundamental difference.  Indeed, as it is easily verified that $\lim_{b \rightarrow \infty} (b \Lambda^{\infty}) = \log(\gamma^{-1})$, while $D^{\infty}_0 = \mu$ irregardless of $b$, what happens in this second setting as $b$ gets large is that $D^{\infty}_{\alpha}$ ramps up very quickly from $\mu$ to $U$, where the time required to ramp up behaves (asymptotically) as $\frac{\log(\gamma^{-1})}{b} T$ (remember that as large as $b$ may be, $T$ is always much larger due to the order of limits here).  However, over this relatively short time horizon (short as a fraction of $T$, still large in an absolute sense as $T$ is very large), a non-trivial backlogging penalty may be incurred since $b$ is so large.  It can be shown that these effects result in a non-trivial limiting dynamics under a proper rescaling, although we do not present a complete analysis here for the sake of brevity.

As a final observation, we note that the comparative results of Theorem\ \ref{compare1theorem} are sensitive to the initial conditions, and in fact for different values of $x_0$ the situation may be reversed.  For example, we have the following result, whose proof we similarly defer to the Technical Appendix in Section\ \ref{rm-sec-appendix}.
\begin{observation}.\label{compare1theorembb}
For all strictly positive $U,b \in {\mathcal Q}^+, T \geq 1$, and $\mu \in [0,U]_{{\mathcal Q}}$, if $x_0 = U$ and $\frac{\mu}{U} < \frac{1}{b+1}$, then 
$\lim_{T \rightarrow \infty} \frac{\inf_{\pi \in \Pi} \sup_{Q \in \textbf{MAR}} \bbe_Q[ \sum_{t=1}^T C^{\pi}_t ]}{\inf_{\pi \in \Pi} \sup_{Q \in \textbf{IND}} \bbe_Q[ \sum_{t=1}^T C^{\pi}_t ]} = \frac{U - \mu}{b \mu} > 1$.
\end{observation}

\indent We also would like to emphasize that many of the insights of our analysis, especially the ``obsolescence phenomena", i.e. the property that conditional on the past demand realizations and the current inventory level (under an optimal policy) there always exists a worst-case distribution that assigns a strictly positive probability to zero and some level which clears the inventory, is quite sensitive to the particular assumptions of our model.  Indeed, although the fact that such a worst-case (conditional) distribution assigns probability to at most 2 points will hold under a broad range of modeling assumptions (and follows essentially from the sparsity of basic feasible solutions to linear programs, which would extend to some other small number of support points under additional constraints), the special properties of putting probability at 0 and clearing the inventory are quite fragile.  In Section\ \ref{nozerohere} of the Technical Appendix, we provide several examples showing that this feature need not hold if one relaxes our modeling assumptions.  Although we defer an explicit description of those examples and findings to Section\ \ref{nozerohere}, we summarize those findings here.  We consider three examples, which collectively demonstrate that : 1. if one allows for time-dependent costs, then a worst-case distribution may not put positive probability at zero; 2. if one removes the upper bound on the support, then a worst-case distribution may not even exist; and 3. if one imposes a lower bound on the support, then a worst-case distribution may not put positive probability on this lower bound, and furthermore may not clear the inventory.  Collectively, these findings suggest that extending our framework to more complex models will require several fundamentally new ideas, as much of the structure which allowed for our explicit analysis may no longer hold.  
\\\indent Of course, the resulting dynamic programs (although more complex) still remain quite structured even in these more general settings, and we believe that searching for new (possibly more complex and subtle) phenomena may lead to the resolution of those settings as well.  For example, in the setting that a lower bound is imposed, although our analysis in Section\ \ref{nozerohere} (specifically the third example) shows that under an optimal policy a worst-case distribution need not clear the inventory, it does not rule out the possibility that in such a setting a worst-case distribution always takes the inventory below the lower bound.  Also, although our analysis in Section\ \ref{nozerohere} (specifically the first example) shows that if the holding costs are held constant and the backlogging costs are monotone increasing (over time) then the ``obsolescence phenomena" need not hold, it may be that imposing other monotonicities and/or appropriate restrictions on the sequence of holding and backlogging costs allows one to recover such obsolescence properties even when costs vary over time.  Indeed, in Section\ \ref{nozerohere} (specifically the fourth example) we do show the following positive result: the obsolescence phenomena again manifests even if holding costs can vary arbitrarily over time, if one enforces (the admittedly strong assumption) that all backlogging costs are 0.  In this  case under a worst-case measure either all demands are 0 or all demands are U.  Although computing the optimal policy is trivial in that setting (as it is to order nothing), the fact that the worst-case measure utilizes obsolescence (in the strongest sense possible) is more subtle, and further demonstrates the worst-case nature of the correlations (over time) induced by obsolescence (especially when holding costs are a primary concern).  Formulating and analyzing such structural properties (which naturally generalize those considered in our analysis), and applying them to develop solution methodologies for more general modeling frameworks, remains an interesting direction for future research.

\section{Proof of Theorem\ \ref{mainmart1} and Corollaries\ \ref{mainmartcorr1} - \ref{mainmartcorr2}}\label{sec-mainproof1}
In this section, we complete the proof of our main result, Theorem\ \ref{mainmart1}, which yields an explicit solution to Problem\ \ref{dynamic0-mar}, as well as Lemma\ \ref{nonempty1}, and 
Corollaries\ \ref{mainmartcorr1} - \ref{mainmartcorr2}.  We proceed by induction, and note that our analysis combines ideas from convex analysis with the theory of martingales.  We begin by making several additional definitions, and proving several auxiliary results.

For $x,\mu \in [0,U]$, and $j \in [-1,T-1]$ , let
$$F^T_j(x,\mu) \stackrel{\Delta}{=} - b x + (b + T) B^T_{j+1} + \big( T b  - (b + 1)(j + 1) \big) \mu.$$
Also, let us define
\begin{align*}
{\mathfrak g}^T(x,\mu)  &\stackrel{\Delta}{=} \begin{cases}
F^T_j(x,\mu) & \textrm{if}\ \mu \in (A^T_j, A^T_{j+1}]\ ,\ x \in [0, B^T_{j+1});\\
G^T_k(x,\mu) & \textrm{if}\ \mu \in (A^T_j, A^T_{j+1}]\ ,\ x \in [B^T_{k}, B^T_{k+1})\ ,\ k \geq j + 1;\\
G^T_k(x,0) & \textrm{if}\ \mu = 0 ,\ x \in [B^T_{k}, B^T_{k+1});\\
G^T_T(U,\mu) & \textrm{if}\ x = U.
\end{cases}
\end{align*}
For later proofs, it will be convenient to note that ${\mathfrak g}$ may be equivalently expressed as follows.
Let
\begin{align*}
\Upsilon^T_x &\stackrel{\Delta}{=} \begin{cases}
T & \textrm{if}\ x = U,\\
j & \textrm{if}\  x \in [B^{T+1}_j, B^{T+1}_{j+1});
\end{cases}
\end{align*}
and note that $x \in [B^T_{\Upsilon^{T-1}_x}, B^T_{\Upsilon^{T-1}_x + 1})$ for all $T \geq 2$ and $x \in [0,U)$, while $U = B^T_{\Upsilon^{T-1}_U+1}$.  Noting that $G^T_T(U,\mu)= G^T_{T-1}(U,\mu)$ implies
\begin{align*}
{\mathfrak g}^T(x,\mu)  &= \begin{cases}
F^T_{\Gamma^{T-1}_{\mu} - 1}(x,\mu) & \textrm{if}\ x < B^T_{\Gamma^{T-1}_{\mu}};\\
G^T_{\Upsilon^{T-1}_x}(x,\mu) & \textrm{if}\ x \geq B^T_{\Gamma^{T-1}_{\mu}}.
\end{cases}
\end{align*}

We now prove several properties of ${\mathfrak g}^T$.  It will first be useful to review a well-known sufficient condition for convexity of non-differentiable functions.
Recall that for a one-dimensional function $f(x)$, the right-derivative of $f$ evaluated at $x_0$, which we denote by $\partial^+_x f(x_0)$, equals $\lim_{h \downarrow 0} \frac{f(x_0 + h) - f(x_0)}{h}$.  When this limit exists, we say that $f$ is right-differentiable at $x_0$.
Then the following sufficient condition for convexity is stated in \citet{royden1988real} Section 5, Proposition 18.
\begin{lemma}[\citet{royden1988real}]\label{suffcon1}
A one-dimensional function $f(x)$, which is continuous and right-differentiable on an open interval $(a,b)$ with non-decreasing right-derivative on $(a,b)$, is convex on $(a,b)$.
\end{lemma}
Then we may derive the following properties of ${\mathfrak g}^T$ and other relevant quantities, whose proofs we defer to the Technical Appendix in Section\ \ref{rm-sec-appendix}.

\begin{lemma}\label{proveme1}
${\mathfrak g}^T(x,d)$ is a continuous and convex function of $x$ on $(0,U)$, and a right (left) continuous function of $x$ at $0\ (U)$.
\end{lemma}
\begin{lemma}\label{interlace1}
$A^T_j < A^{T+1}_{j+1} < A^T_{j+1}$ for $j \in [0,T-2]$.  It follows that for all $T \geq 2$ and $d \in [0,U]$, $\Gamma^T_d \in \lbrace \Gamma^{T-1}_d, \Gamma^{T-1}_d + 1 \rbrace$.
\end{lemma}
\begin{lemma}\label{proveme2}
$\beta^T_d \in \argmin_{z \in [0,U]} {\mathfrak g}^T(z,d)$.
\end{lemma}

To simplify various notations and concepts before completing the proof of Theorem\ \ref{mainmart1}, it will be useful to make several additional definitions, and prove a few more relevant preliminary results.  For $x,d \in [0,U]$, and $j \in [0,T]$, let
$$\overline{F}^T_j(x,d) \stackrel{\Delta}{=} T x + \big( (b-1) T - b j - \frac{b+T}{A^T_j} x \big) d + \frac{b+T}{A^T_j} d^2;$$
and for $j \in [-1,T-1]$, let
$$\overline{G}^T_j(d) \stackrel{\Delta}{=} T B^T_{j+1} + \big( T b - (b + 1)(j+ 1) \big) d.$$
Note that $\beta^T_d$ is monotone increasing in $d$, with $\beta^T_0 = 0$, and $\beta^T_U = U$.  It follows that for all $x \in [0,U]$,
$$z^T_x \stackrel{\Delta}{=} \inf \lbrace d \geq 0\ \ \ \textrm{s.t.}\ \ \ \beta^T_d \geq x - d \rbrace,$$
is well defined,
\begin{equation}\label{ztcompx}
z^T_x \leq x\ \ \ \textrm{for}\ \ \ x \in [0,U],
\end{equation}
and $z^T_x = 0$ iff $x = 0$.
For $x,d \in [0,U]$, let us define
\begin{align*}
\overline{{\mathfrak g}}^T(x,d)  &\stackrel{\Delta}{=} \begin{cases}
\overline{F}^T_j(x,d) & \textrm{if}\ d \in [0, z^T_x\big) \bigcap (x - B^T_{j+1}, x - B^T_j]\ ,\ j \in [0,T-1];\\
\overline{G}^T_j(d) & \textrm{if}\ d \in [z^T_x, U] \bigcap (A^{T+1}_j, A^{T+1}_{j+1}]\ ,\ j \in [-1,T-1];\\
\overline{G}^T_{-1}(0) & \textrm{if}\ d = 0\ , x = 0;\\
\overline{F}^T_{T-1}(U,0) & \textrm{if}\ d = 0\ , x = U.
\end{cases}
\end{align*}
For later proofs, it will be convenient to note that $\overline{{\mathfrak g}}$ may be equivalently expressed as follows.
\begin{align*}
\overline{{\mathfrak g}}^T(x,d)  &= \begin{cases}
\overline{F}^T_{\Upsilon^{T-1}_{x-d}}(x,d) & \textrm{if}\ d < z^T_x;\\
\overline{G}^T_{\Gamma^T_d - 1}(d) & \textrm{if}\ d \geq z^T_x.
\end{cases}
\end{align*}

In that case, we may derive the following properties of $\overline{{\mathfrak g}}^T(x,d)$, whose proofs we defer to the Technical Appendix in Section\ \ref{rm-sec-appendix}.
\begin{lemma}\label{bareqbar2}
For all $T \geq 1$, $x,d \in [0,U]$, $\overline{{\mathfrak g}}^T(x,d) = {\mathfrak g}^T\bigg( \max\big( \beta^T_d , x - d \big) , d \bigg)$.
\end{lemma}

\begin{lemma}\label{convexconcave1}
For each fixed $x \in [0,U]$,
$\overline{{\mathfrak g}}^T(x,d)$ is a continuous function of $d$ on $(0,U)$, and a right (left) - continuous function of $d$ at $0\ ( U )$.  Also, $\overline{{\mathfrak g}}^T(x,d)$ is a convex function of $d$ on $(0, z^T_x)$, and a concave function of $d$ on $(z^T_x, U)$.
\end{lemma}

Also, for $x,d \in [0,U]$, let us define $\eta_d \stackrel{\Delta}{=} x - d$, $\alpha^T_x \stackrel{\Delta}{=} A^T_{\Gamma^{T-1}_x}, \zeta^T_x \stackrel{\Delta}{=} \Gamma^T_{z^T_x}, {\mathcal A}^T_x \stackrel{\Delta}{=} A^T_{\zeta^{T-1}_x}$, $\aleph^T_x \stackrel{\Delta}{=} A^T_{\Upsilon^{T-1}_x}$,
$${\mathfrak f}^T(x,d)  \stackrel{\Delta}{=} b (d - x)_+ + (x - d)_+ +\overline{{\mathfrak g}}^{T-1}(x,d).$$
Note that for all $j \in [-1, T-1]$,
\begin{eqnarray}
B^T_{j+1} &=& \frac{j+1}{b+T} A^T_{j+1} \nonumber
\\&\leq& \frac{j+1}{b+j+1} A^T_{j+1}\ \ \ =\ \ \ A^T_j. \label{BAcompare}
\end{eqnarray}
Combining with (\ref{ztcompx}) and a straightforward contradiction argument, we conclude that for all $x \in [0,U]$ and $T \geq 2$,
\begin{equation}\label{BAcompare2}
\Upsilon^T_x \geq \Gamma^T_x \geq \zeta^T_x\ \ \ ,\ \ \ \textrm{and}\ \ \ \aleph^T_x \geq \alpha^T_x \geq {\mathcal A}^T_x.
\end{equation}

In that case, we may derive the following structural properties of ${\mathfrak f}^T$, whose proofs we defer to the Technical Appendix in Section\ \ref{rm-sec-appendix}.

\begin{lemma}\label{piececon1}
For each fixed $x \in [0,U]$, ${\mathfrak f}^T(x,d)$ is a continuous function of $d$ on $(0,U)$, and a right (left) - continuous function of $d$ at $0\ ( U )$.  Also, for all $j \in [-1, T-2]$, ${\mathfrak f}^T(x,d)$ is a convex function of $d$ on $(A^T_j, A^T_{j+1})$.  Furthermore, ${\mathfrak f}^T(x,d)$ is a convex function of $d$ on $(0, {\mathcal A}^T_x)$, and a concave function of $d$ on $(\alpha^T_x , U)$.
\end{lemma}

We now introduce one last preliminary result, which will provide a mechanism for certifying a distribution as the solution to a given distributionally robust optimization problem, whose proof we again defer to the Technical Appendix in Section\ \ref{rm-sec-appendix}.

\begin{lemma}\label{linedom1}
Suppose $f:[0,U] \rightarrow \bbr$ is any bounded real-valued function with domain $[0,U]$.  Suppose $0 \leq L \leq \mu < R \leq U$, and the linear function $\eta(d)\ \stackrel{\Delta}{=} \frac{ f(R) - f(L) }{R - L} d + \frac{ R f(L) - L f(R) }{R - L}$, i.e. the line intersecting $f$ at the points $L$ and $R$, satisfies $\eta(d) \geq f(d)$ for all $d \in [0,U]$, i.e. lies above $f$ on $[0,U]$.  Then the measure $q$ s.t. $q(L) = \frac{R - \mu}{R - L}, q(R) = \frac{\mu - L}{R - L}$, belongs to $\argmax_{Q \in \mathfrak{M}(\mu)} \bbe_Q[f(D)]$.
\end{lemma}

With all of the above notations, definitions, and preliminary results in place, we now proceed with the proof of Theorem\ \ref{mainmart1}, which will follow from the following result.

\begin{theorem}\label{centralmart1}
For all $x,d \in [0,U]_{{\mathcal Q}}$ and $T \geq 1$, ${\mathfrak g}^T(x,d) = \hat{g}^T(x,d)$, and ${\mathfrak q}^T_{x,d} \in \hat{Q}^T(x,d)$.
\end{theorem}

\proof{Proof:} We proceed by induction, beginning with the base case $T = 1$.  It follows from Observation\ \ref{explicitt2} that for all $x,\mu \in [0,U]_{{\mathcal Q}}$, $\Gamma^0_{\mu} = \Upsilon^0_x = B^1_{\Gamma^0_{\mu}} = 0, \Gamma^1_{\mu} = I(\mu > \frac{U}{b+1}), \beta^1_\mu =  U \times I(\mu > \frac{U}{b+1}), {\mathfrak g}^1(x,\mu) = G^1_0(x,\mu) = (1 - \frac{b+1}{U} \mu) x + b \mu,$ and ${\mathfrak q}^1_{x,\mu}(0) = 1 - \frac{\mu}{U}\ ,\ {\mathfrak q}^1_{x,\mu}(U) = \frac{\mu}{U}$.  The desired result then follows directly from the results of \citet{S-12}, and we omit the details.
\\\\Now, suppose that the result is true for all $s \in [1,T-1]$ for some $T \geq 2$.  We now prove that the results hold also for $T$.  It follows from the induction hypothesis, and Lemmas\ \ref{proveme1},\ \ref{proveme2},\ and\ \ref{bareqbar2} that for all $x,d \in [0,U]_{{\mathcal Q}}$, 
\begin{equation}\label{viswhat}
\hat{V}^{T-1}(x,d) = {\mathfrak g}^{T-1}\bigg( \max\big( \beta^{T-1}_d , x \big) , d \bigg)\ \ \ ,\ \ \ \hat{V}^{T-1}(x - d, d) = \overline{{\mathfrak g}}^{T-1}(x,d)\ \ \ ,\ \ \ \textrm{and}\ \ \ \hat{f}^T(x,d) = {\mathfrak f}^T(x,d).
\end{equation}
From definitions, to prove the desired induction, it thus suffices to demonstrate that 
\begin{equation}\label{showgf1}
{\mathfrak g}^T(x,d) = \sup_{Q \in \mathfrak{M}(d)} \bbe_Q[ {\mathfrak f}^T(x,D) ]\ \ \ ,\ \ \ \textrm{and}\ \ \ {\mathfrak q}^T_{x,d} \in \argmax_{Q \in \mathfrak{M}(d)} \bbe_Q[ {\mathfrak f}^T(x,D) ].
\end{equation}
We now use Lemma\ \ref{piececon1} to explicitly construct (for each fixed value of $x$) a family ${\mathcal F}$ of lines $\lbrace L_i \rbrace$, s.t. each line $L_i$ lies above ${\mathfrak f}^T(x,d)$ for all $d \in [0,U]$, each line $L_i$ intersects
${\mathfrak f}^T(x,d)$ at exactly two points $p^1_i, p^2_i$, and for each $\mu \in [0,U]$ there exists $L_i \in {\mathcal F}$ s.t. $\mu \in [p^1_i,p^2_i]$.  Our construction will ultimately allow us to apply Lemma\ \ref{linedom1}, explicitly solve the distributionally robust optimization problem $\sup_{Q \in \mathfrak{M}(d)} \bbe_Q[ {\mathfrak f}^T(x,D) ]$, and complete the proof.  We begin by explicitly constructing the family of lines ${\mathcal F}$.  For $d \in \bbr$, let us define
$${\mathcal K}^T(x,d) \stackrel{\Delta}{=} \frac{ {\mathfrak f}^T(x,\aleph^T_x ) - T x}{\aleph^T_x} d + T x;$$
and for $j \in [-1,T-2]$,
$${\mathcal L}^T_j(x,d) \stackrel{\Delta}{=}  \frac{ {\mathfrak f}^T(x, A^T_{j+1}) - {\mathfrak f}^T(x, A^T_j)}{A^T_{j+1} - A^T_j} d + \frac{ A^T_{j+1} {\mathfrak f}^T(x, A^T_j) - A^T_j {\mathfrak f}^T(x, A^T_{j+1})}{A^T_{j+1} - A^T_j}.$$
Noting that $\overline{F}^{T-1}_j(x,0) = (T - 1) x$ for all $j$, it follows from Lemma\ \ref{bareqbar2} that
\begin{equation}\label{f00}
{\mathfrak f}^T(x,0) = T x.
\end{equation}
It may be easily verified, using (\ref{f00}), that ${\mathcal K}^T$ defines the unique line passing through the $(x,y)$ co-ordinates
$\big(0, {\mathfrak f}^T(x,0) \big)$ and $\big(\aleph^T_x, {\mathfrak f}^T(x, \aleph^T_x) \big)$; and
${\mathcal L}^T_j$ defines the unique line passing through the $(x,y)$ co-ordinates $\big( A^T_j, {\mathfrak f}^T(x,A^T_j) \big)$ and $\big( A^T_{j+1}, {\mathfrak f}^T(x, A^T_{j+1}) \big).$  Then we have the following result, showing that ${\mathcal K}^T$ lies above ${\mathfrak f}^T$, and that ${\mathcal L}^T_\ell$ lies above ${\mathfrak f}^T$ for certain values of $\ell$, whose proof we defer to the Technical Appendix in Section\ \ref{rm-sec-appendix}.
\begin{lemma}\label{wearefamily1}
For each fixed $x \in [0,U]$, ${\mathcal K}^T(x,d) \geq {\mathfrak f}^T(x,d)$ for all $d \in [0,U]$.  Also, for all $\ell \in [\Upsilon^{T-1}_x, T-2]$, ${\mathcal L}^T_\ell(x,d) \geq {\mathfrak f}^T(x,d)$ for all $d \in [0,U]$.
\end{lemma}

That ${\mathfrak q}^T_{x,d} \in \argmax_{Q \in \mathfrak{M}(d)} \bbe_Q[ {\mathfrak f}^T(x,D) ]$ follows from Lemmas\ \ref{linedom1} and \ref{wearefamily1}, definitions, and a straightforward case analysis, the details of which we omit.  We now prove that ${\mathfrak g}^T(x,d) = \sup_{Q \in \mathfrak{M}(d)} \bbe_Q[ {\mathfrak f}^T(x,D) ]$, and proceed by a case analysis.  Let $j = \Gamma^{T-1}_d$ and $k = \Upsilon^{T-1}_x$.  First, suppose $d \in [0, A^T_k]$.  In light of Lemmas\ \ref{linedom1}\ and\ \ref{wearefamily1}, in this case it suffices to demonstrate that
\begin{equation}\label{showgend1}
{\mathfrak g}^T(x,d) = {\mathcal K}^T(x,d).
\end{equation}
In this case, the left-hand side of (\ref{showgend1}) equals
\begin{equation}\label{doornumberg1}
G^T_k(x,d) =  (T - \frac{b+T}{A^T_k} d) x + (T - k) b d.
\end{equation}
Alternatively, from \eqref{zetaiswhat1}, the right-hand side of (\ref{showgend1}) equals
\begin{equation}\label{doornumberg1b}
   \big( b (T - k) - \frac{ (b + T) x }{A^T_k} \big) d + T x.
\end{equation}
Noting that (\ref{doornumberg1}) equals (\ref{doornumberg1b}) completes the proof in this case.
\\\indent Alternatively, suppose $d > A^T_k$.  In this case it suffices to demonstrate that
\begin{equation}\label{showgend2}
{\mathfrak g}^T(x,d) =  {\mathcal L}^T_{j-1}(x,d).
\end{equation}
Note that the left-hand side of (\ref{showgend2}) equals
\begin{equation}\label{doornumberg2}
F^T_{j-1}(x,d) =  - b x + (b + T) B^T_j + \big( T b  - (b + 1) j \big) d.
\end{equation}
Alternatively, it follows from (\ref{ztcompx}) and \eqref{BAcompare} that the right-hand side of (\ref{showgend2}) equals
\begin{eqnarray}
&\ &\ \frac{ {\mathfrak f}^T(x, A^T_j) - {\mathfrak f}^T(x, A^T_{j-1})}{A^T_j - A^T_{j-1}} d +
\frac{ A^T_j {\mathfrak f}^T(x, A^T_{j-1}) - A^T_{j-1} {\mathfrak f}^T(x, A^T_j)}{A^T_j - A^T_{j-1}} \nonumber
\\&\ &\ \ \ =\ \ \ b (d - x)  + \frac{ \overline{G}^{T-1}_{j-1}(A^T_{j}) - \overline{G}^{T-1}_{j-2}(A^T_{j-1})}{A^T_{j} - A^T_{j-1}} d  + \frac{ A^T_{j} \overline{G}^{T-1}_{j-2}(A^T_{j-1}) - A^T_{j-1} \overline{G}^{T-1}_{j-1}(A^T_{j})}{A^T_{j} - A^T_{j-1}}. \label{ee02}
\end{eqnarray}
We now simplify (\ref{ee02}).  Note that $\overline{G}^{T-1}_{j-1}(A^T_{j}) - \overline{G}^{T-1}_{j-2}(A^T_{j-1})$ equals
\begin{eqnarray}
\ &\ &\ (T -1)(B^{T-1}_{j} - B^{T-1}_{j-1}) +
\bigg( (T - 1) b - (b + 1) j \bigg) (A^T_{j} - A^T_{j-1}) - (b + 1) A^T_{j-1}  \nonumber
\\&\ &\ \ \ =\bigg(  (T - 1) b - (b + 1) j \bigg) \left(A^T_{j} - A^T_{j-1}\right)
+  (T -1)\left(\frac{j}{T - 1}A^T_{j}  - \frac{j-1}{T - 1}A^T_{j-1}\right) - (b + 1) A^T_{j-1} \nonumber
\\&\ &\ \ \ =\bigg( (T - 1) b - (b + 1) j \bigg)  \left(A^T_{j} - A^T_{j-1}\right), \label{firstsimplejack1}
\end{eqnarray}
and $A^T_{j} \overline{G}^{T-1}_{j-2}(A^T_{j-1}) - A^T_{j-1} \overline{G}^{T-1}_{j-1}(A^T_{j})$ equals
\begin{eqnarray}
\ &\ &\ \ \  (T - 1) \left(  A^T_jB^{T-1}_{j-1}- A^T_{j-1}B^{T-1}_{j} \right)+   (b+1) A^T_{j-1} A^T_{j} \nonumber
\\&\ &\ \ \ \ \ = (T - 1) \left(  A^T_j A^{T}_{j-1}\frac{j-1}{T-1}- A^T_{j-1}A^{T}_{j} \frac{j}{T-1} \right)+   (b+1) A^T_{j-1} A^T_{j} \nonumber
\\&\ &\ \ \ \ \ = b A^T_{j-1} A^T_{j} \ \ \ =\ \ \ j A^T_{j} \left(A^T_{j} - A^T_{j-1}\right). \label{firstsimplejack2}
\end{eqnarray}
Plugging (\ref{firstsimplejack1}) and (\ref{firstsimplejack2}) into (\ref{ee02}), and comparing to (\ref{doornumberg2}), completes the proof of (\ref{showgf1}), and the desired induction.  $\Halmos$
\endproof
With Theorem\ \ref{centralmart1} in hand, we now complete the proof of our main result Theorem\ \ref{mainmart1}, as well as Corollaries\ \ref{mainmartcorr1} and\  \ref{mainmartcorr2}. 
\proof{Proof:}[Proof of Theorem\ \ref{mainmart1}]
Lemma\ \ref{nonempty1} follows from Theorem\ \ref{centralmart1} and Lemmas\ \ref{proveme1} and\ \ref{proveme2}, and it follows from definitions and a straightforward case analysis (the details of which we omit) that $\text{Opt}^T_{\textbf{MAR}}(\mu,U,b) = {\mathfrak g}^T(\beta^T_{\mu},\mu)$.  Theorem\ \ref{mainmart1} then follows by combining Lemmas\ \ref{proveme1}\ and\ \ref{proveme2} with Theorem\ \ref{centralmart1}, (\ref{viswhat}), Lemma\ \ref{nonempty1}, and Theorem\ \ref{thm-4.2}.  $\Halmos$
\endproof
\proof{Proof:}[Proof of Corollaries\ \ref{mainmartcorr1}\ and\ \ref{mainmartcorr2}]
First, note that if $\Gamma^T_{\mu} = T$, then $\chi^T_{\textbf{MAR}}\big(\mu, U, b\big) = B^T_T = A^T_{T-1} = X^T_1 = D^T_1 = U$, proving the desired result in this case.  Next, suppose $\Gamma^T_{\mu} \leq T - 1$.  In light of the definitions of $X^T_t$ and $D^T_t$, i.e.\ (\ref{rm-eq-def}), the already established monotonicities of $\lbrace A^T_j, j \in [-1,T-1] \rbrace$ and $\lbrace B^T_j, j \in [-1,T] \rbrace$, and the martingale property, it suffices to establish that for all $t \in [1,\Lambda^T]$: 1. $\chi^{T+1-t}_{\textbf{MAR}}\big(D^T_{t-1}, U, b\big) = X^T_t$, 2. $\exists j \in [-1,T-1], k \geq j + 1$ s.t. $D^T_{t-1} \in (A^{T+1-t}_j, A^{T+1-t}_{j+1}], X^T_t \in [B^{T+1-t}_k, B^{T+1-t}_{k+1}), A^{T+1-t}_k = D^T_t$.  
\\\\First, suppose $t = 1$.
In this case, $D^T_{t-1} = \mu \in (A^T_{\Gamma^{T-1}_{\mu} - 1}, A^T_{\Gamma^{T-1}_{\mu}}],$
$$\chi^{T+1-t}_{\textbf{MAR}}\big(D^T_{t-1}, U, b\big) = \chi^{T}_{\textbf{MAR}}\big(\mu, U, b\big) = B^T_{\Gamma^{T}_{\mu}} = X^T_t,$$
and $D^T_t = A^T_{\Gamma^{T}_{\mu}}$.  Note that since $A^{\tau}_j$ is monotone decreasing in $\tau$, it follows that $\Gamma^{\tau}_{\mu}$ is monotone increasing in $\tau$, and thus $\Gamma^T_{\mu} \geq \Gamma^{T-1}_{\mu}$.  Thus we find that 1. holds, and 2. holds with $j = \Gamma^{T-1}_{\mu} - 1, k = \Gamma^T_{\mu}$.
\\\\Next, suppose $t \in [2, \Lambda^T]$.  In this case, 
$D^T_{t-1} = z \stackrel{\Delta}{=} A^{T+2-t}_{\Gamma^T_{\mu}}$.  Note that 
\begin{equation}\label{afact99}
z \in (A^{T+2-t}_{\Gamma^T_{\mu}-1}, A^{T+2-t}_{\Gamma^T_{\mu}}],
\end{equation}
 and thus 
$\Gamma^{T+1-t}_z = \Gamma^T_{\mu}$.  Furthermore, 
$$\chi^{T+1-t}_{\textbf{MAR}}\big(D^T_{t-1}, U, b\big) = B^{T+1-t}_{\Gamma^{T+1-t}_z} = B^{T+1-t}_{\Gamma^T_{\mu}} = X^T_t.$$  
In addition, the fact that $A^{\tau}_j$ is monotone decreasing in $\tau$ and monotone increasing in j, combined with (\ref{afact99}), implies that $z \in (A^{T+1-t}_{j'-1}, A^{T+1-t}_{j'}]$ for some $j' \leq \Gamma^T_{\mu}$.  Thus we find that 1. holds, and 2. holds with $j = j' - 1, k = \Gamma^T_{\mu}$, completing the proof.  $\Halmos$
\endproof

\section{Asymptotic analysis and proof of Theorem\ \ref{weakconvergethm1}}\label{sec-mainproof2}
In this section, we complete the proof of Theorem\ \ref{weakconvergethm1}, and begin with a lemma which asserts the (uniform) convergence of several one-dimensional functions to appropriate limits.  As all of these convergences follow from repeated application of elementary bounds such as Taylor's Inequality for the exponential function, and the continuous differentiability of all relevant functions over appropriate domains, we omit the details.  For $\alpha \in [0,1)$, let $\mathtt{f}(\alpha)\stackrel{\Delta}{=}  b \left( 1-\alpha\right)^{b-1}$, where we note that $\mathtt{f}(\alpha) = Z^{\infty}(\alpha)$ for $\alpha \in [0, \Lambda^{\infty})$.
\begin{lemma}\label{rm-lemma1}
For all $\mu\in (0,U)$, 
$$\lim_{T\to \infty} \big( T^{-1} \Gamma^{T}_{\mu} \big) = \gamma^{\frac{1}{b}}\ \ ,\ \ \lim_{T\to \infty} A^{T+1}_{\Gamma^{T}_{\mu}}= \mu\ \ ,\ \ \lim_{T\to \infty} \left( T^{-1} \Lambda^T \right)=  \Lambda^{\infty};$$
$$
\limsup_{T \rightarrow \infty} \sup_{t \in [1, \Lambda^T]} |X^T_t - X^{\infty}_{\frac{t}{T}}|\ \ =\ \ 
\limsup_{T \rightarrow \infty} \sup_{t \in [1, \Lambda^T]} |D^T_t - D^{\infty}_{\frac{t}{T}}|\ \ =\ \ \limsup_{T \rightarrow \infty} \sup_{t \in [1, \Lambda^T - 1]} |T {\mathcal Z}^T(t) - \mathtt{f}(\frac{t}{T})|\ \ =\ \ 0.
$$
\end{lemma}

For all $T \geq 1, \alpha \geq 0$, let $m_{\alpha}(T) \stackrel{\Delta}{=} \min(\lfloor\alpha T\rfloor, \Lambda^T-1)$, and $m_{\alpha} \stackrel{\Delta}{=} \min(\alpha,\Lambda^{\infty})$.  For $\alpha_1, \alpha_2\geq 0$ s.t. $\alpha_1< m_{\alpha_2}$, let us define
$X^T_{\alpha_1,\alpha_2}\stackrel{\Delta}{=} X^{T}_{Z^T}I\big(\lfloor\alpha_1 T\rfloor+1\leq Z^T\leq m_{\alpha_2}(T)\big)$, and 
$X^\infty_{\alpha_1,\alpha_2}\stackrel{\Delta}{=} X^{\infty}_{Z^\infty}I\big(\alpha_1\leq Z^\infty< m_{\alpha_2}\big).$
Also, for $\mu \in (0,U)$, note that $\mathtt{f}$ is continuous on $[0,\Lambda^{\infty}]$, and let $\overline{\mathtt{f}} \stackrel{\Delta}{=} \sup_{\alpha \in [0,\Lambda^{\infty}]} \mathtt{f}(\alpha) < \infty$.  In light of Lemma\ \ref{rm-lemma1}, it follows that 
\begin{equation}\label{boundtz}
\limsup_{T \rightarrow \infty} \sup_{t \in [1, \Lambda^T - 1]} \big( T {\mathcal Z}^T(t) \big) = \overline{\mathtt{f}}.
\end{equation}
Let us also recall that $X^{\infty}$ and $D^{\infty}$ are continuous and strictly increasing on $[0,\Lambda^{\infty}]$.

Before embarking on the proof of Theorem\ \ref{weakconvergethm1}, it will be useful to first establish the weak convergence of $X^T_{\alpha_1,\alpha_2}$ to $X^\infty_{\alpha_1,\alpha_2}$.
\begin{lemma}\label{rm-lemma4}
For $\mu \in (0,U), \alpha_1, \alpha_2 \geq 0$ s.t. $\alpha_1 < m_{\alpha_2}$, $X^T_{\alpha_1,\alpha_2}$ converges weakly to $X^\infty_{\alpha_1,\alpha_2}$.  Furthermore, 
for any bounded uniformly continuous function $H: [0,U] \rightarrow \bbr$,
\begin{equation}\label{rm-50}
\begin{aligned}
   \lim_{T\to \infty}  \sum_{k=\lfloor\alpha_1 T\rfloor+1}^{m_{\alpha_2}(T)} H\left(X^{T}_{k}\right){\mathcal Z}^T(k) =\  \int_{\alpha_1}^{m_{\alpha_2}} H\left( X^{\infty}_{z}\right)\mathtt{f}(z)dz.
\end{aligned}
\end{equation}
\end{lemma}
\proof{Proof:}
By the Portmanteau Theorem (cf. Theorem 2.1 in \citet{Billingsley}), to prove weak convergence, it suffices to demonstrate (\ref{rm-50}).  Let $\overline{H}$ denote an upper bound for $|H|$.  It follows from Lemma\ \ref{rm-lemma1}, (\ref{boundtz}), and the fact that $|a b - c d| \leq |a + c| |b - d| + |b + d| |a - c|$ that for any $\epsilon > 0$, there exists $T_{\epsilon} < \infty$ s.t. for all $T \geq T_{\epsilon}$, and $k \in [\lfloor\alpha_1 T\rfloor+1,m_{\alpha_2}(T)]$,
$$\big|H\left(X^{T}_{k}\right)\big(T {\mathcal Z}^T(k)\big) - H(X^{\infty}_{\frac{k}{T}}) \mathtt{f}(\frac{k}{T})\big|
\leq 
2 \overline{H} \epsilon + 4 \overline{\mathtt{f}} \sup_{k \in [\lfloor\alpha_1 T\rfloor+1,m_{\alpha_2}(T)]}|H\left(X^{T}_{k}\right) - H(X^{\infty}_{\frac{k}{T}})\big|.
$$
It follows from Lemma\ \ref{rm-lemma1} that there similarly exists $T_{1,\epsilon} < \infty$ s.t. for all $T  \geq T_{1,\epsilon}$,
$\sup_{k \in [\lfloor\alpha_1 T\rfloor+1,m_{\alpha_2}(T)]}|X^{T}_{k} - X^{\infty}_{\frac{k}{T}}\big| < \epsilon.$
Combining the above with the fact that $H$ is uniformly continuous implies that for any $\epsilon > 0$, there exists $T_{2,\epsilon} < \infty$ s.t. for all $T \geq T_{2,\epsilon}$,
$
\sup_{k \in [\lfloor\alpha_1 T\rfloor+1,m_{\alpha_2}(T)]}\big|H\left(X^{T}_{k}\right)\big(T {\mathcal Z}^T(k) \big) - H(X^{\infty}_{\frac{k}{T}}) \mathtt{f}(\frac{k}{T})\big| < \epsilon.
$
Thus for any $\epsilon > 0$, there exists $T_{3,\epsilon} < \infty$ s.t. for all $T \geq T_{3,\epsilon}$,
$
\bigg|T^{-1}\sum_{k=\lfloor\alpha_1 T\rfloor+1}^{m_{\alpha_2}(T)} H\left(X^{T}_{k}\right)\big(T {\mathcal Z}^T(k) \big) - T^{-1} \sum_{k=\lfloor\alpha_1 T\rfloor+1}^{m_{\alpha_2}(T)} H(X^{\infty}_{\frac{k}{T}}) \mathtt{f}(\frac{k}{T})\bigg| < \epsilon.
$
Noting that Riemann integrability implies $\lim_{T \rightarrow \infty} T^{-1} \sum_{k=\lfloor\alpha_1 T\rfloor+1}^{m_{\alpha_2}(T)} H(X^{\infty}_{\frac{k}{T}}) \mathtt{f}(\frac{k}{T}) = \int_{\alpha_1}^{m_{\alpha_2}} H\left( X^{\infty}_{z}\right)\mathtt{f}(z)dz$
completes the proof.  $\Halmos$
\endproof

We now complete the proof of Theorem \ref{weakconvergethm1}, and proceed by the standard path of demonstrating appropriate notions of tightness and convergence of finite-dimensional distributions.  We begin by providing some additional background regarding the space $D\left([0,1],\ \bbr^k\right)$ and associated notions of weak convergence.  Recall that the space $D\left([0,1],\ \bbr^k\right)$ contains all functions $x: [0,1]\to \bbr^k$ that are right-continuous in $[0,1)$ and have left-limits in $(0,1]$ (i.e. cadlag on $[0,1]$).  Given such a cadlag function $x = (x_1,\ldots,x_k): [0,1] \rightarrow \bbr^k$, let $||x_j|| \stackrel{\Delta}{=} \sup_{t\in[0,1]}|x_j(t)|$, and $||x||\stackrel{\Delta}{=} \sum_{j=1}^k ||x_j||$.
Given a set $\{\alpha_i\}_{i=0}^m$ satisfying $0=\alpha_0<\alpha_1<\ldots<\alpha_m=1$, we say that $\{\alpha_i\}_{i=0}^m$ is $\delta$-sparse if $\min_{1\leq i\leq m}(\alpha_i-\alpha_{i-1})\geq \delta$.  Let $S_{\delta,m}$ denote the collection of all such $\delta$-sparse sets consisting of exactly $m$ points.  Given $x$ which is cadlag on $[0,1]$ and $\delta \in (0,1)$, let
\[
   w_x(\delta)\stackrel{\Delta}{=}\ \inf_{m \in Z^+, \{\alpha_i\} \in S_{\delta,m}} \max_{1\leq i\leq m} \left(\sum_{j=1}^k\sup_{\alpha,\alpha'\in [\alpha_{i-1}, \alpha_i)}|x_j(\alpha)-x_j(\alpha')|\right),
\]
where the infimum extends over all $\delta$-sparse sets $\{\alpha_i\}$ of all possible finite lengths (i.e. values of $m$), where we note that one must only consider values of $m$ which are at most $\lceil \delta^{-1} \rceil + 1$.  Let us recall the following necessary and sufficient conditions for tightness in the space $D\left([0,1],\ \bbr^k\right)$.
\begin{lemma}[Theorem 13.2, \citet{Billingsley}]\label{bill1}
A sequence of probability measures $\{P_T\}_{T\geq 1}$ on the space $D\left([0,1],\ \bbr^k\right)$ under the $J_1$ topology is tight if and only if the following conditions hold:
\begin{equation}\label{rm-31}
   \lim_{a\to \infty}\limsup_{T \rightarrow \infty}P_T\left(||x||\geq a\right)=0;
\end{equation}
\begin{equation}\label{rm-32}
\textrm{For all}\ \epsilon > 0\ \ ,\ \ \lim_{\delta\to 0}\limsup_{T\rightarrow \infty}P_T\left(w_x(\delta)\geq \epsilon\right)=0.
\end{equation}
\end{lemma}
We now use Lemma\ \ref{bill1} to prove tightness of the sequence $\lbrace {\mathcal M}^T(\alpha)_{0 \leq \alpha \leq 1}, T \geq 1 \rbrace$, by verifying \eqref{rm-31} and \eqref{rm-32}.
\begin{lemma}\label{rm-lemmatight}
For $\mu \in (0,U)$, $\lbrace {\mathcal M}^T(\alpha)_{0 \leq \alpha \leq 1}, T \geq 1 \rbrace$ is tight in $D\left([0,1],\ \bbr^2\right)$.
\end{lemma}
\proof{Proof:}
\eqref{rm-31} follows immediately from the fact that ${\mathcal M}^{T}(\alpha)$ is supported on $[0, U]^2$ w.p.1 for all $\alpha\in [0,1]$, $T\geq 1$. Let us prove \eqref{rm-32}.  Recall from Corollary\ \ref{mainmartcorr2} that we may construct r.v.s $Z^T, Y^T$ on the same probability space as ${\mathcal M}^T$ s.t. the properties outlined in Corollary\ \ref{mainmartcorr2} hold.  It follows from Lemma\ \ref{rm-lemma1}, and the Riemann integrability of $\mathtt{f}$, that there exists $\delta_0 > 0$ s.t. for all $\delta \in (0,\delta_0)$, and $\delta' \in (0,\delta)$, there exists $T_{\delta'}$ s.t. for all $T \geq T_{\delta'}$,
\begin{equation}\label{weakimp1}
\sup_{t \in [1, \Lambda^T]} |X^T_t - X^{\infty}_{\frac{t}{T}}| < \delta'\ \ \ ,\ \ \ \sup_{t \in [1, \Lambda^T]} |D^T_t - D^{\infty}_{\frac{t}{T}}| < \delta'\ \ \ ,\ \ \ |\pr(\frac{Z^T}{T} \leq \delta) - \int_0^{\delta} \mathtt{f}(z)\ dz| < \delta'.
\end{equation}
Note that $X^{\infty}, D^{\infty}, \mathtt{f}$ are continuously differentiable on $[0,\Lambda^{\infty} + \epsilon_0]$ for some $\epsilon_0 > 0$.  Let $\overline{D} \stackrel{\Delta}{=} 32 \big(1 + \sup_{0 \leq x \leq \Lambda^{\infty} + \epsilon_0} \max(|\frac{d}{dx} D^{\infty}_x|, |\frac{d}{dx} X^{\infty}_x|, |\frac{d}{dx} \mathtt{f}|) \big)$.  Also note that Lemma\ \ref{rm-lemma1} implies that $\lbrace Y^T, T \geq 1 \rbrace$ converges weakly to the r.v. which equals $U$ w.p.1.

It thus follows from (\ref{weakimp1}) and the triangle inequality that for all $\delta \in (0,\delta_0)$, there exists $T_{\delta}$ s.t. for all $T \geq T_{\delta}$,
\begin{equation}\label{weakimp2}
\sup_{\substack{0 \leq \alpha_1 \leq \alpha_2 \leq 1\\ \lfloor \alpha_2 T \rfloor \leq \Lambda^T\\|\alpha_1 - \alpha_2| \leq 2 \delta}} \max\bigg( \big| X^{T}_{\lfloor\alpha_2 T\rfloor} - X^{T}_{\lfloor\alpha_1 T\rfloor} \big|, \big| D^{T}_{\lfloor\alpha_2 T\rfloor} - D^{T}_{\lfloor\alpha_1 T\rfloor} \big| \bigg) \leq \overline{D} \delta\ \ ,\ \ |X^{T}_{\Lambda^T} - U| \leq \overline{D} \delta,
\end{equation}
$$\pr\big( \big\lbrace \min(\frac{Z^T}{T}, 1 - \frac{Z^T}{T}) \leq 4 \delta \rbrace \bigcup \lbrace Y^T = 0 \rbrace \big) \leq \overline{D} \delta.
$$
As we will be referring to various sample-path constructions in the remainder of the proof, for a given $\omega \in \Omega$ (on the corresponding probability space), we will denote the dependence of various r.v.s on $\omega$ by adding $\omega$ as a subscript or superscript where appropriate.  Given $\delta \in (0,\delta_0)$ and $T \geq T_{\delta}$, let $\Omega_{\delta,T} \subseteq \Omega$ denote that subset of $\Omega$ s.t. $\omega \in \Omega_{\delta,T}$ iff $\min(\frac{Z^{T,\omega}}{T}, 1 - \frac{Z^{T,\omega}}{T}) > 4 \delta$ and $Y^{T,\omega} = U$.
Let us fix some $\epsilon > 0, \delta \in (0,\delta_0), T \geq T_{\delta}, \omega \in \Omega_{\delta,T}$.  For the given parameters, we now define a particular $\delta$-sparse mesh.  Note that for the given $\omega$, there exists $k_{\omega}$ s.t. $k_{\omega} \delta \leq \frac{Z^{T,\omega}}{T} < (k_{\omega}+1) \delta$, where $0 < (k_{\omega} - 2) \delta < (k_{\omega} + 2) \delta < 1$.  Also, let $n_{1,\omega} \stackrel{\Delta}{=} \sup\lbrace n: \frac{Z^{T,\omega}}{T} + n \delta < 1 - \delta \rbrace$, and $n_{\omega} \stackrel{\Delta}{=} k_{\omega} + n_{1,\omega} + 2$.  Then we let $\lbrace \alpha^{\omega}_i, i \in [1,n_{\omega}] \rbrace$ denote the following $\delta$-sparse mesh.  Let $\alpha^{\omega}_1 = 0$.  For $k \in [1, k_{\omega} - 1]$, let $\alpha^{\omega}_k = k \delta$.  Let $\alpha^{\omega}_{k_{\omega}} = \frac{Z^{T,\omega}}{T}$.  For $k \in [k_{\omega} + 1, k_{\omega} + 1 + n_{1,\omega}]$, let $\alpha^{\omega}_k = \frac{Z^{T,\omega}}{T} + (k - k_{\omega}) \delta$.  Finally, we let $\alpha^{\omega}_{n_{\omega}} = 1$.  
Note that 
\begin{equation}\label{smallalpha1}
\frac{Z^{T,\omega}}{T} = \alpha^{\omega}_{k_{\omega}} = \inf \lbrace \alpha \in [0,1]: \lfloor \alpha T \rfloor \geq Z^{T,\omega} \rbrace,
\end{equation} 
$$i \in [2, n_{\omega}+1], \alpha,\alpha'\in [\alpha^{\omega}_{i-1}, \alpha^{\omega}_i)\ \textrm{implies}\ |\alpha-\alpha'|\leq 2\delta.
$$
We treat two cases. First, suppose $Z^{T,\omega} \leq \Lambda^T-1$.  In this case, Corollary\ \ref{mainmartcorr2} implies that
\[
   {\mathcal M}^{T,\omega}(\alpha)= \begin{cases}
   \big( X^{T}_{\lfloor\alpha T\rfloor}, D^{T}_{\lfloor\alpha T\rfloor} \big) & \text{if}\ \lfloor\alpha T\rfloor\leq Z^{T,\omega}-1;\\
   \big( X^{T}_{Z^{T,\omega}}, 0\big) & \text{if}\ \lfloor\alpha T\rfloor\geq Z^{T,\omega}.
   \end{cases}
\]
It thus follows from (\ref{weakimp2}) and (\ref{smallalpha1}) that 
$$\max_{i \in [2,k_{\omega}]} \sum_{j=1}^2 \sup_{\alpha,\alpha'\in [\alpha_{i-1}, \alpha_i)} |{\mathcal M}^{T,\omega}_j(\alpha) - {\mathcal M}^{T,\omega}_j(\alpha')| \leq 2 \overline{D} \delta,$$
$$\max_{i \in [k_{\omega} + 1,n_{\omega} + 1]} \sum_{j=1}^2 \sup_{\alpha,\alpha'\in [\alpha_{i-1}, \alpha_i)} |{\mathcal M}^{T,\omega}_j(\alpha) - {\mathcal M}^{T,\omega}_j(\alpha')| = 0,$$
and thus in this case $w^{\omega}_x(\delta) \leq 2 \overline{D} \delta$.

Alternatively, suppose $Z^{T,\omega} = \Lambda^T$.  In this case, Corollary\ \ref{mainmartcorr2} implies that
\[
   {\mathcal M}^{T}(\alpha)= \begin{cases}
   \big( X^{T}_{\lfloor\alpha T\rfloor}, D^{T}_{\lfloor\alpha T\rfloor} \big) & \text{if}\ \lfloor\alpha T\rfloor\leq \Lambda^T-1;\\
   \big( X^{T}_{\Lambda^T}, U\big) & \text{if}\ \lfloor\alpha T\rfloor= \Lambda^T;\\
   \big( U, U\big) & \text{if}\ \lfloor\alpha T\rfloor\geq \Lambda^T+1.\\
   \end{cases}
\]
It again follows from (\ref{weakimp2}) and (\ref{smallalpha1}) that
$$\max_{i \in [2,k_{\omega}]} \sum_{j=1}^2 \sup_{\alpha,\alpha'\in [\alpha_{i-1}, \alpha_i)} |{\mathcal M}^{T,\omega}_j(\alpha) - {\mathcal M}^{T,\omega}_j(\alpha')| \leq 2 \overline{D} \delta.$$
Now, note that for $i \in [k_{\omega} + 1,n_{\omega} + 1]$ and $\alpha,\alpha' \in [\alpha_{i-1}, \alpha_i)$,
$|{\mathcal M}^{T,\omega}_2(\alpha) - {\mathcal M}^{T,\omega}_2(\alpha')| = 0$, and
 $|{\mathcal M}^{T,\omega}_1(\alpha) - {\mathcal M}^{T,\omega}_1(\alpha')| \leq |X^{T}_{\Lambda^T} - U| \leq \overline{D} \delta$, with the final inequality following from (\ref{weakimp2}).  Thus in this case, it again holds that $w^{\omega}_x(\delta) \leq 2 \overline{D} \delta$.

To complete the proof of (\ref{rm-32}), we need only carefully combine all of the above.  In particular, let us fix any $\epsilon > 0$.  It follows from (\ref{weakimp2}) that for all $\delta \in (0,\delta_0)$ and $T \geq T_{\delta}$, $\pr(\Omega_{\delta,T}) \geq 1 - 
\overline{D} \delta$, $w^{\omega}_x(\delta) \leq 2 \overline{D} \delta$ for all $\omega \in \Omega_{\delta,T}$, and thus $\pr\big(w^{\omega}_x(\delta) > 2 \overline{D} \delta\big) \leq \overline{D} \delta$.  We conclude that for all $\delta < \min(\delta_0, \frac{\epsilon}{3 \overline{D}})$
and $T \geq T_{\delta}$, $\pr\big(w^{\omega}_x(\delta) \geq \epsilon \big) \leq \overline{D} \delta$, from which (\ref{rm-32}) follows after taking limits in the appropriate order.  $\Halmos$
\endproof

Next, we prove the following weak convergence of finite-dimensional distributions.
\begin{lemma}\label{rm-lemmafinite}
For $\mu \in (0,U)$, and $0 = \alpha_1< \ldots < \alpha_{k_0} = \Lambda^{\infty} < \ldots < \alpha_n = 1$, $\left({\mathcal M}^{T}(\alpha_1), \ldots, {\mathcal M}^{T}(\alpha_n)\right)$ converges weakly to $\left({\mathcal M}^{\infty}(\alpha_1), \ldots, {\mathcal M}^{\infty}(\alpha_n)\right)$.
\end{lemma}

\proof{Proof:}
By the Portmanteau theorem, it suffices to prove that for any bounded uniformly continuous function $H:[0,U]^{2n}\to \bbr$,
\begin{equation}\label{rm-weak}
   \lim_{T\to \infty} \bbe\left[H\left({\mathcal M}^{T}(\alpha_1), \ldots, {\mathcal M}^{T}(\alpha_n)\right)\right]=\ \bbe\left[H\left({\mathcal M}^{\infty}(\alpha_1), \ldots, {\mathcal M}^{\infty}(\alpha_n)\right)\right].
\end{equation}
As a notational convenience, we define for $i \in [1,n]$,
\begin{eqnarray*}
   & H^\infty\left(i, x, y\right)&\stackrel{\Delta}{=}\   H\left(X^{\infty}_{\alpha_1}, D^{\infty}_{\alpha_1}, \ldots, X^{\infty}_{\alpha_i}, D^{\infty}_{\alpha_i},  x, y, x, y, \ldots, x, y\right),\\
   & H^T\left(i, x, y\right)&\stackrel{\Delta}{=}\   H\left(X^{T}_{\lfloor\alpha_1 T\rfloor}, D^{T}_{\lfloor\alpha_1 T\rfloor}, \ldots, X^{T}_{\lfloor\alpha_{i} T\rfloor}, D^{T}_{\lfloor\alpha_{i} T\rfloor},  x, y, x, y, \ldots, x, y\right),\\
   & \hat{H}^T\left(i, x,y\right)&\stackrel{\Delta}{=}\    H\left(X^{T}_{\lfloor\alpha_1 T\rfloor}, D^{T}_{\lfloor\alpha_1 T\rfloor}, \ldots, X^{T}_{\lfloor\alpha_{i} T\rfloor}, D^{T}_{\lfloor\alpha_{i} T\rfloor}, x, y, y, y, \ldots, y,y\right).
\end{eqnarray*}

Then it follows from the definition of ${\mathcal M}^{\infty}(\alpha)$ that
\[
   \bbe\left[H\left({\mathcal M}^{\infty}(\alpha_1), \ldots, {\mathcal M}^{\infty}(\alpha_n)\right)\right]=\ \sum_{i=1}^{k_0-1} \int_{\alpha_i}^{\alpha_{i+1}}  H^\infty\left(i, X^{\infty}_{z}, 0\right) \mathtt{f}(z) dz + \gamma   H^\infty\left(k_0, U,U\right).
\]
It follows from Lemma \ref{rm-lemma1} that there exists $T_0$ s.t. for all $T \geq T_0$, $\lfloor\alpha_{k_0-1} T\rfloor< \Lambda^T - 1 < \Lambda^T< \lfloor\alpha_{k_0+1} T\rfloor$.  From Corollary \ref{mainmartcorr2}, for all $T \geq T_0$,
\[
   \bbe\left[H\left({\mathcal M}^{T}(\alpha_1), \ldots, {\mathcal M}^{T}(\alpha_n)\right)\right]=\ W_1^T+W_2^T+W_3^T+W_4^T,
\]
where
\begin{eqnarray*}
   &W_{1}^T\stackrel{\Delta}{=}& \sum_{i=1}^{k_0-2} \sum_{j=\lfloor\alpha_i T\rfloor+1}^{\lfloor\alpha_{i+1} T\rfloor}  H^T\left(i, X_j^T, 0\right) \bbp(Z^T=j)  +\  \sum_{j=\lfloor\alpha_{k_0-1} T\rfloor+1}^{\min\{\lfloor\alpha_{k_0} T\rfloor, \Lambda^T-1\}}  H^T\left(k_0-1, X_j^T, 0\right)  \bbp(Z^T=j),\\
   &W_2^T\stackrel{\Delta}{=}&  I\left(\lfloor\alpha_{k_0} T\rfloor< \Lambda^T-1\right)   \sum_{j=\lfloor\alpha_{k_0} T\rfloor+1}^{\Lambda^T-1}  H^T\left(k_0, X_j^T, 0\right)\bbp(Z^T=j),\\
   &W_3^T\stackrel{\Delta}{=}&  \Bigg[ I\left(\lfloor\alpha_{k_0} T\rfloor= \Lambda^T\right)   \hat{H}^T\left(k_0-1, X^{T}_{\Lambda^T},U\right)     +\   I\left(\lfloor\alpha_{k_0} T\rfloor> \Lambda^T\right)   H^T\left(k_0-1, U,U\right)\\
   && \ \ \ \ \ \ +\  I\left(\lfloor\alpha_{k_0} T\rfloor< \Lambda^T\right)   H^T\left(k_0, U,U\right)  \Bigg]   \bbp\left(Z^T= \Lambda^T\right) \bbp\left(Y^T= U\right),\\
   &W_4^T\stackrel{\Delta}{=}& \Bigg[ I\left(\lfloor\alpha_{k_0} T\rfloor\geq \Lambda^T\right)  H^T\left(k_0-1, X^{T}_{\Lambda^T},0\right)\\
   && \ \ \ \ \ \ +\    I\left(\lfloor\alpha_{k_0} T\rfloor< \Lambda^T\right)  H^T\left(k_0, X^{T}_{\Lambda^T},0\right)
   \Bigg] \bbp\left(Z^T= \Lambda^T\right) \bbp\left(Y^T= 0\right).
\end{eqnarray*}
Lemma\ \ref{rm-lemma4} implies that $\lim_{T \rightarrow \infty} W^T_1 = \sum_{i=1}^{k_0-1} \int_{\alpha_i}^{\alpha_{i+1}}  H^\infty\left(i, X^{\infty}_{z}, 0\right) \mathtt{f}(z) dz$, and it follows from Lemma\ \ref{rm-lemma1} and the uniform continuity of $H$ that 
$\lim_{T \rightarrow \infty} W^T_4 = \gamma   H^\infty\left(k_0, U,U\right)$, and $\lim_{T \rightarrow \infty} W^T_2 = \lim_{T \rightarrow \infty} W^T_3 = 0$.  Combining the above completes the proof of \eqref{rm-weak}.   $\Halmos$
\endproof

\proof{Proof:}[Proof of Theorem\ \ref{weakconvergethm1}]
Combining Lemmas\ \ref{bill1},\ \ref{rm-lemmatight},\ and\ \ref{rm-lemmafinite} completes the proof of Theorem\ \ref{weakconvergethm1}.  $\Halmos$
\endproof

\section{Numerical experiments}\label{sec-nume}
In this section, we conduct numerical experiments to further demonstrate the benefits of the proposed policy.  More specifically, we aim to show the benefits of a robust policy which accounts for martingale structure over a robust policy which assumes that the demand is drawn independently over time.  To that end, we will compare the performance of the minimax optimal policy for the independent-demand model to the performance of the minimax optimal policy for the martingale-demand model, when the demand itself comes from a simple additive-MMFE model (and has the martingale property).  Indeed, in line with the numerical experiments conducted in \citet{Mamani17} (in which the joint distribution of demand was multivariate normal), we assume that 
the true demand process is the following martingale: $D_n=\mu+\sum_{t=1}^n\epsilon_t$, where $\mu$ is a constant and $\{\epsilon_t\}_{t\geq 1}$ is a sequence of i.i.d. r.v.s representing forecast adjustments, i.e. the demand comes from a simple additive-MMFE model.  We further assume that $\lbrace \epsilon_t, t \in [1,T] \rbrace$ are i.i.d., normally distributed with mean $0$ and variance $\sigma^2$.  We note that such an additive-MMFE model with normally distributed increments is common in the literature (e.g., \citet{HJ}, \citet{LSR}, \citet{WAK}), making it a good reference for numerical comparisons.  As a base-line we assume that the demand distribution has a per-period average of $10$, i.e., $\mu=10$, and the standard deviation $\sigma$ equals either $0.1\mu$ (relatively low) or $0.2\mu$ (relatively high).  We normalize the holding cost at $h =1$ and vary backorder costs in the specified range, namely, the service level $\frac{b}{b+h}$ is varied from $10\%$ (low) to $90\%$ (high).  The finite time horizon $T$ is varied from relatively short ($T=3$) to relatively long ($T=20$).  
\\\indent We note that although both policies we consider will be seeded by an upper-bound parameter U (as required in the definition of the policy),which we will vary from $U = 15$ to $U = 25$, the demand process itself in no cases depends on U, and may go above U (and may actually go below zero, in which case the demand actually adds to the inventory, although we note that such negative demand will be a relatively rare occurence under most of our parameter settings).  As any particular model assumptions may not hold in any given practical setting, it is in some sense natural to consider the performance when such a support assumption may fail.  Recall that the minimax optimal policy in the independent-demand setting (from Theorem\ \ref{thm-4.1}) is a base-stock policy with base-stock level $\chi_{\textbf{IND}}(\mu,U,b)$, which we note is independent of time or the realized demands, and in any case is either always 0 or always U (depending on the relation of $\mu,U$, and $b$).  Thus this policy can be implemented with no alterations even if demand leaves the interval [0,U].  Indeed, we let Policy$_\text{IND}$ denote exactly this policy, which we shall use in our numerical experiments.  
\\\indent Things are slightly more complicated in the martingale-demand setting.  Indeed, recall that the minimax optimal policy $\hat{\pi}^*$ in the martingale-demand setting (from Theorem\ \ref{mainmart1}) sets $\hat{x}^*_1 = \max\big(x_0, \chi^T_{\textbf{MAR}}(\mu,U,b)\big)$;  and For $t \in [1,T-1]$, sets $\hat{x}^*_{t+1}(\mathbf{d}_{[t]}) = \max\big(\hat{x}^*_t(\mathbf{d}_{[t-1]}) - d_t, \chi^{T-t}_{\textbf{MAR}}(d_t,U,b) \big).$  Formally, $\hat{\pi}^*$ is a state-and-time-dependent base-stock policy, ordering up to level $\chi^{T - t}_{\textbf{MAR}}(D_t,U,b)$ in period $t+1$ for $t \in [1,T-1]$, where we note that this level depends both on time and the most recent demand value.  Here, a difficulty arises if for some $t$ it holds that $D_t \notin [0,U]$, since in that case $\chi^{T - t}_{\textbf{MAR}}(D_t,U,b)$ is not defined.  We remedy this by defining a policy Policy$_\text{MAR}$ as follows.  In period 1, the policy orders up to level $\chi^T_{\textbf{MAR}}(\mu,U,b)$, consistent with $\hat{\pi}^*$.  For $t \in [1,T-1]$, if $D_t \in [0,U]$, then at the start of period $t + 1$ the policy orders up to level $\chi^{T - t}_{\textbf{MAR}}(D_t,U,b)$, consistent with $\hat{\pi}^*$.  However, if for some $t \in [1,T-1]$ it holds that $D_t < 0$, then at the start of period $t + 1$ the policy orders up to $0$.  Alternatively, if for some $t \in [1,T-1]$ it holds that $D_t > U$, then at the start of period $t + 1$ the policy orders up to $U$.  We note that this is in some sense the simplest and most natural way to reconcile this discrepancy, and it is this (slightly modified) policy Policy$_\text{MAR}$ that we use in our numerical experiments.
\\\indent In summary, the goal of our numerical experiment is to compare the performance of Policy$_\text{MAR}$ and Policy$_\text{IND}$ when the demand process is a Gaussian random walk (i.e. additive-MMFE).  For each parameter setting, we perform $10^6$ simulations, with our table displaying the average cost incurred under each of the two policies for each parameter setting.  In all cases we set the initial inventory level to 0.  Although (for clarity of exposition) we do not include a formal analysis of variance or confidence intervals, we note that our experiments suggested that the outcomes were generally quite stable.  The following table summarizes the different parameter settings which we will consider.

\begin{table}[h]
\centering
    \begin{tabular}{|c|c|c|c|c|c|c|c|c|}
      \hline
      Parameters & $h$ & $b$ & $T$ & $U$ & $\lbrace \epsilon_t, t \in [1,T] \rbrace$ & $\mu$ & $\sigma$ & $x_0$ \\
      \hline
      Values & 1 & $\frac{1}{9}$, $\frac{1}{4}$, 1, 4, 9 & 3, 10, 20 & 15, 20, 25 & i.i.d. Normal(0, $\sigma^2$) & 10 & 0.1$\mu$, 0.2$\mu$ & $0$\\
      \hline
    \end{tabular}
\caption{Summary of parameters in the numerical experiments.}\label{mart-table-1}
\end{table}

We report our numerical results in Table \ref{mart-table-2}. For each set of parameter values, we report the following three numbers: the average cost incurred by Policy$_\text{MAR}$ (denoted by C$_\text{MAR}$), the average cost incurred by Policy$_\text{IND}$ (denoted by C$_\text{IND}$), and the percentage of cost reduction defined as $\frac{C_\text{IND} - C_\text{MAR}}{C_\text{IND}}$.  Each cell contains the corresponding triplet (C$_\text{MAR}$, C$_\text{IND}$, $\frac{C_\text{IND} - C_\text{MAR}}{C_\text{IND}}$).

\begin{table}[h]
\centering
    \begin{tabular}{|c|c|c|c|}
      \hline
      $U$ = 15, $\sigma$ = 0.1$\mu$  & $T$ = 3 & $T$ = 10 & $T$ = 20 \\
      \hline
      $b$ = 1/9  & (3.333, 3.333, 0.00\%) & (10.19, 11.11, 8.28\%) & (18.66, 22.52, 17.1\%) \\
      \hline
      $b$ = 1/4  & (5.991, 7.500, 20.1\%) & (18.13, 25.00, 27.5\%) & (33.03, 50.31, 34.3\%) \\
      \hline
      $b$ = 1    &  (11.66, 15.00, 22.3\%) & (26.57, 50.49, 47.4\%) & (52.31, 104.6, 50.0\%) \\
      \hline
      $b$ = 4    &  (15.01, 15.01, 0.00\%) & (40.17, 51.22, 21.6\%) & (66.47, 111.1, 40.2\%) \\
      \hline
      $b$ = 9    &  (15.01, 15.01, 0.00\%) & (51.77, 52.44, 1.28\%) & (99.89, 122.0, 18.1\%) \\
      \hline
      \hline
      $U$ = 20, $\sigma$ = 0.1$\mu$  & $T$ = 3 & $T$ = 10 & $T$ = 20 \\
      \hline
      $b$ = 1/9  & (3.333, 3.333, 0.00\%) & (11.09, 11.11, 0.18\%) & (21.99, 22.52, 2.35\%) \\
      \hline
      $b$ = 1/4  & (7.500, 7.500, 0.00\%) & (23.67, 25.00, 5.32\%) & (44.32, 50.31, 11.9\%) \\
      \hline
      $b$ = 1    &  (18.88, 30.00, 37.1\%) & (43.82, 100.0, 56.2\%) & (84.20, 200.4, 58.0\%) \\
      \hline
      $b$ = 4    &  (30.00, 30.00, 0.00\%) & (53.39, 100.0, 46.6\%) & (84.69, 200.6, 57.8\%) \\
      \hline
      $b$ = 9    &  (30.00, 30.00, 0.00\%) & (87.81, 100.0, 12.2\%) & (112.7, 201.0, 43.9\%) \\
      \hline
      \hline
      $U$ = 25, $\sigma$ = 0.1$\mu$  & $T$ = 3 & $T$ = 10 & $T$ = 20 \\
      \hline
      $b$ = 1/9  & (3.333, 3.333, 0.00\%) & (11.11, 11.11, 0.00\%) & (22.48, 22.52, 0.18\%) \\
      \hline
      $b$ = 1/4  & (7.500, 7.500, 0.00\%) & (24.92, 25.00, 0.32\%) & (48.60, 50.31, 3.40\%) \\
      \hline
      $b$ = 1    &  (17.91, 30.00, 40.3\%) & (55.17, 100.0, 44.8\%) & (107.1, 200.4, 46.6\%) \\
      \hline
      $b$ = 4    &  (41.34, 45.00, 8.13\%) & (64.92, 150.0, 56.7\%) & (107.4, 300.2, 64.2\%) \\
      \hline
      $b$ = 9    &  (45.00, 45.00, 0.00\%) & (96.36, 150.0, 35.8\%) & (135.3, 300.2, 54.9\%) \\
      \hline
      \hline
      $U$ = 15, $\sigma$ = 0.2$\mu$  & $T$ = 3 & $T$ = 10 & $T$ = 20 \\
      \hline
      $b$ = 1/9  & (3.211, 3.336, 3.75\%) & (9.911, 12.35, 19.7\%) & (37.32, 44.60, 16.3\%) \\
      \hline
      $b$ = 1/4  & (5.476, 7.503, 27.0\%) & (16.76, 26.30, 36.3\%) & (49.90, 72.95, 31.6\%) \\
      \hline
      $b$ = 1    &  (10.84, 15.33, 29.3\%) & (29.30, 57.90, 49.4\%) & (80.55, 150.1, 46.3\%) \\
      \hline
      $b$ = 4    &  (15.77, 15.81, 0.25\%) & (53.94, 68.65, 21.4\%) & (151.1, 198.4, 23.8\%) \\
      \hline
      $b$ = 9    &  (16.63, 16.63, 0.00\%) & (80.62, 86.58, 6.88\%) & (251.2, 278.9, 9.93\%) \\
      \hline
      \hline
      $U$ = 20, $\sigma$ = 0.2$\mu$  & $T$ = 3 & $T$ = 10 & $T$ = 20 \\
      \hline
      $b$ = 1/9  & (3.335, 3.336, 0.03\%) & (11.65, 12.35, 5.67\%) & (40.92, 44.60, 8.25\%) \\
      \hline
      $b$ = 1/4  & (7.445, 7.503, 0.77\%) & (22.23, 26.30, 15.5\%) & (58.62, 72.95, 19.6\%) \\
      \hline
      $b$ = 1    &  (18.11, 30.00, 39.6\%) & (41.09, 101.6, 59.6\%) & (94.89, 226.0, 58.0\%) \\
      \hline
      $b$ = 4    &  (29.63, 30.01, 1.27\%) & (56.73, 103.2, 45.0\%) & (130.3, 239.7, 45.6\%) \\
      \hline
      $b$ = 9    &  (30.01, 30.02, 0.03\%) & (87.94, 105.6, 16.7\%) & (176.0, 261.5, 32.7\%) \\
      \hline
      \hline
      $U$ = 25, $\sigma$ = 0.2$\mu$  & $T$ = 3 & $T$ = 10 & $T$ = 20 \\
      \hline
      $b$ = 1/9  & (3.336, 3.336, 0.00\%) & (12.22, 12.35, 1.05\%) & (43.12, 44.60, 3.32\%) \\
      \hline
      $b$ = 1/4  & (7.503, 7.503, 0.00\%) & (24.89, 26.30, 5.36\%) & (65.22, 72.95, 10.6\%) \\
      \hline
      $b$ = 1    &  (18.80, 30.00, 37.3\%) & (51.96, 101.6, 48.9\%) & (113.6, 226.0, 49.7\%) \\
      \hline
      $b$ = 4    &  (40.48, 45.00, 10.0\%) & (68.18, 150.9, 54.8\%) & (139.8, 322.7, 56.7\%) \\
      \hline
      $b$ = 9    &  (44.98, 45.00, 0.04\%) & (96.32, 151.2, 36.3\%) & (175.6, 327.5, 46.4\%) \\
      \hline
   \end{tabular}
\caption{Costs incurred by Policy$_\text{MAR}$ and Policy$_\text{IND}$, and  percentages of cost reduction.}\label{mart-table-2}
\end{table}

Let us summarize our findings.  First, and most importantly, in all cases the cost incurred by Policy$_\text{MAR}$ is no greater than the cost incurred by Policy$_\text{IND}$, and the percentage of cost reduction can be as large as $64\%$.  This dominance, and (for some parameter settings) substantial improvement, suggest the importance of accounting for possible dependencies when implementing robust inventory control policies.  We note that in the deterministic robust optimization setting, related numerical results were presented in \cite{Mamani17}.
\\\indent Next, we comment on how the parameter choices impact the improvement.  First, we note that for several parameter settings, the two policies perform nearly identically.  This stems from the fact that for certain parameters, both policies will with high probability either always order up to U or always order up to 0.  Although Policy$_\text{IND}$ always exhibits such a behavior, a careful analysis of Theorem\ \ref{mainmart1} and Policy$_{\text{MAR}}$ shows that for certain values of $T, t, D_t, U$, and $b$, the thresholds of Policy$_{\text{MAR}}$ will exhibit the same behavior.  For example, we proved this equivalence for sufficiently large b in our previous Corollary\ \ref{largeb1}, in which case all relevant thresholds will be U (with high probability).  A similar behavior can be demonstrated as regards settings in which both policies set the relevant thresholds to 0 (again with high probability in the case of Policy$_{\text{MAR}}$).  We note that such a phenomena occurs almost exclusively when T = 3 and/or b takes either very small or very large values, which is consistent with the parameter regimes under which such agreement would occur (as indicated in Corollary\ \ref{largeb1}).  Indeed, Corollary\ \ref{largeb1} indicates that such agreement would occur (for example) if $b > (\frac{U}{D_t} - 1) \times (T - t)$ for all $t \in [0,T-1]$ (here letting $D_0 \stackrel{\Delta}{=} \mu$), which is more likely to occur when b is large and/or T is small.  Again, we note that a similar phenomena can occur for alternative parameter ranges in which case both policies order up to 0 instead of ordering up to U.  We note that when $T = 3$, this behavior is partially mitigated by the presence of a larger variance, i.e. when $T = 3$ the larger variance tends to yield a larger gap between the policies, as the larger variance can lead to realized demands which break the conditions needed for these degeneracies.  
\\\indent Second, we note that our results seem to be fairly insensitive to U, $\sigma$, and T, so long as $T \geq 10$.  Intuitively, for any such parameters the conditional expectation of demand will with high probability fluctuate considerably, leading to a significant advantage for Policy$_{\text{MAR}}$, which utilizes this information.  Indeed, for all such parameter settings, the gap fluctuates fairly consistently as a function of $b$, taking values near $50\%$ for values of b closer to 1, and taking smaller values for other settings of b.  We offer two possible explanations for this phenomena.  First, we note that for the parameters we consider, b being close to 1 tends to coincide with $\frac{U}{\mu \times (b + 1)}$ being close to 1.  We note that $\frac{U}{\mu \times (b + 1)}$ equals 1 exactly when Policy$_\text{IND}$ is indifferent between always ordering up to 0, and always ordering up to U.  Our results suggest that it is exactly in this regime, namely when $\frac{U}{\mu \times (b + 1)}$ is close to 1, that the relative gap between the two policies tends to be largest.
\\\indent However, a more refined analysis reveals that when T and U are both large, this effect seems to be mitigated on one side, as the gap remains large even when $\frac{U}{\mu \times (b + 1)}$ is small, although the gap shrinks when $\frac{U}{\mu \times (b + 1)}$ becomes large.  To help explain this refinement, we reason alternatively as follows, building on our earlier discussion about how certain parameter values can lead to the two policies degenerating and behaving similarly.  It follows from Theorem\ \ref{mainmart1} that in the first period Policy$_{\text{MAR}}$ orders up to 0 (i.e. does not order at all) if $\frac{\mu}{U} \in [0, \prod_{k=1}^T \frac{k}{b + k})$, and orders up to U if $\frac{\mu}{U} \in (\frac{T}{b + T}, 1]$.  Of course, Policy$_{\text{MAR}}$ may also order up to many intermediate values.  However, the above implies that if $\frac{\mu}{U} \in \big[0, \min(\prod_{k=1}^T \frac{k}{b + k} , \frac{1}{b + 1}) \big)$, or $\frac{\mu}{U} \in \big( \frac{T}{b + T}, 1 \big]$, then Policy$_{\text{MAR}}$ and Policy$_{\text{IND}}$ degenerate to a common policy (at least in the first period).  Thus if b is relatively large, since $\prod_{k=1}^T \frac{k}{b + k}$ will be extremely small, the primary manner for such a degeneracy/agreement to occur would be to have $\frac{\mu}{U} \in \big( \frac{T}{b + T}, 1 \big],$ equivalently $b \geq T \times (\frac{U}{\mu} - 1)$. The only parameter settings we consider for which this condition holds either require $T = 3$, or $T = 10, U = 15, b = 9$.  We note that indeed, the setting in which $T = 10, U = 15, b = 9$ is essentially the only setting in which T and b are moderately large, yet the gap between the policies remains very small.  
\\\indent Alternatively, if b is relatively small, since $\frac{T}{b + T}$ is then close to 1, the primary manner for such a degeneracy/agreement to occur would be to have $\frac{\mu}{U} \in [0, \prod_{k=1}^T \frac{k}{b + k})$ (since $\frac{1}{b+1}$ is very close to 1 for small values of b).  As a Taylor-series approximation shows that $\prod_{k=1}^T \frac{k}{b + k}$ behaves roughly like $T^{-b}$ for small b, this is (roughly) equivalent to requiring that $b \leq \frac{\log(\frac{U}{\mu})}{\log(T)}$.  However, this inequality holds for many parameter settings involving large T, specifically : $U = 15, T = 10, b = \frac{1}{9}; U = 15, T = 20, b = \frac{1}{9}; U = 20, T = 10, b = \frac{1}{9}; U = 20, T = 10, b = \frac{1}{4}; U = 20, T = 20, b = \frac{1}{9}; U = 25, T = 10, b = \frac{1}{9}; U = 25, T = 10, b = \frac{1}{4}; U = 25, T = 20, b = \frac{1}{9}; U = 25, T = 20, b = \frac{1}{4}$.  
\\\indent Thus we find that there is a substantial asymmetry between small and large values of b, with small values of b leading to more scenarios in which the two policies degenerate and behave identically, at least in the first period.  Intuitively, this should lead to smaller gaps for smaller values of b, while the presence of a large gap should be more robust for larger values of b, exactly as reflected in our numerical experiments.  Of course, the exact behavior of Policy$_{\text{MAR}}$ (and associated thresholds) over time depends on the realized demand (not just $\mu$), yet one might expect this intuition to (roughly) hold broadly over the entire time horizon.
\\\indent In summary, we find that Policy$_\text{MAR}$ substantially outperforms Policy$_\text{IND}$ for these settings, further demonstrating the utility of our model and the importance of taking dependencies into consideration when modeling robustness.  Our experiments also suggest many interesting and subtle phenomena, whose full exploration we leave as an interesting direction for future research.  
\section{Conclusion}\label{sec-conclusion}
In this paper, we formulated and solved a dynamic distributionally robust multi-stage newsvendor model which naturally unifies the analysis of several inventory models with demand forecasting, and
enables the optimizer to incorporate past realizations of demand into the structure of the uncertainty set going forwards.  We explicitly computed the minimax optimal policy (and associated worst-case distribution) in closed form.  Our main proof technique involved a non-trivial induction, combining ideas from convex analysis, probability, and DP.  By analyzing in-depth the interplay between the minimax optimal policy and associated worst-case distribution, we proved that at optimality the worst-case demand distribution corresponds to the setting in which inventory may become obsolete at a random time, a scenario of practical interest.  To gain further insight into our explicit solution, we computed the limiting dynamics as the time horizon grows large, by proving weak convergence to an appropriate limiting stochastic process, for which many quantities have a simple and intuitive closed form.  We also compared to the analogous setting in which demand is independent across periods, and made several comparisons between these two models.  Finally, we complemented our results by providing a targeted and concise numerical experiment further demonstrating the benefits of our model.
\\\indent Our work leaves many interesting directions for future research.  First, to make the model more practically applicable, it would be interesting to extend our results to more complicated dependency structures which robustify a larger class of forecasting models, e.g. that for which $E[D_{t+1} | D_{[t]}] = a_t D_t + b_t$ (i.e. general time-dependent linear conditional expectations).  One could also incorporate higher order dependencies, non-linearities, as well as higher-dimensional covariate and feature spaces such as those arising in factor models.  Similarly, one could consider inventory models with more realistic features, e.g. various notions of ordering costs, settings in which demand is lost, positive lead times, etc.  Of course, our general approach could also be applied to many models used by the Operations Research community where forecasting and model uncertainty non-trivially interact.  In all of these cases, it would be interesting to develop a more general theory of when the associated DP have a special structure which allows for a closed-form (or at least computationally tractable) solution.  We note that such questions may also shed further light on the search for optimal probability inequalities in a variety of settings, e.g. bounds for functionals of martingale sequences and their generalizations.  As suggested by our examples in Section\ \ref{nozerohere}, these generalizations will likely require several fundamentally new ideas and insights.
\\\indent Taking a broader view, perhaps the most interesting directions for future research involve exploring the relationship between different ways to model uncertainty in multi-stage optimization problems.  What is the precise relationship between static and dynamic models of uncertainty, and how should we think about comparing and selecting between different models?  This question involves notions of computational efficiency, model appropriateness in any given application, as well as questions related to the so-called price of correlations.  Of course, such issues are intimately related to various relevant questions of a statistical nature.  How precisely should our ability to fit and tune a model relate to our choice of uncertainty set?  How does the notion of model over-fitting come into play here?  Such questions become especially interesting in the context of multi-stage problems.  Another intriguing line of questioning involves the relation between classical probability and (distributionally) robust optimization.  For example, the obsolescence feature of the worst-case martingale in our setting can be interpreted as a ``robust" manifestation of the celebrated martingale convergence theorem.  A deeper understanding of the precise connection between the limit laws of classical probability and the solutions to distributionally robust optimization problems remains an intriguing open question, where we note that closely related questions have been the subject of several recent investigations (cf. \citet{abernethy13,bandi15a,bandi15b}).  Formalizing and analyzing these and related questions will likely require a combination of ideas and techniques from optimization, statistics, and probability, all of which bring to bear different perspectives with which to understand uncertainty.

\ACKNOWLEDGMENT{The authors would like to thank Erhan Bayraktar, Vishal Gupta, Garud Iyengar, Sandeep Juneja, Anton Kleywegt, Marcel Nutz, Alex Shapiro, and Alejandro Toriello for helpful discussions and insights.  David A. Goldberg also gratefully acknowledges support from NSF grant no. 1757394.}

\bibliographystyle{nonumber}

\newpage
\section{Technical Appendix}\label{rm-sec-appendix}
\subsection{Proof of Lemma\ \ref{lemma0u}}
\proof{Proof of Lemma\ \ref{lemma0u}:}
We first prove that one can w.l.o.g. restrict to policies always ordering up to at least 0.  For $s \in [1,T]$, let $\hat{\Pi}^{s,0}$ denote that subset of $\hat{\Pi}$ consisting of those policies $\hat{\pi}$ with the following property: for all $t \in [T - s + 1, T]$ and $\mathbf{d} = (d_1,\ldots,d_T) \in [0,U]_{{\mathcal Q}}^T$, $x^{\hat{\pi}}_t(\mathbf{d}_{[t-1]}) \geq 0$.  Also, let $\hat{\Pi}^{0,0} \stackrel{\Delta}{=} \hat{\Pi}$.  We now prove that for all $s \in [0,T-1]$, and any policy $\hat{\pi} \in \hat{\Pi}^{s,0}$, there exists a policy $\hat{\pi}' \in \hat{\Pi}^{s + 1,0}$ s.t. for all vectors $\mathbf{d} = (d_1,\ldots,d_T) \in [0,U]_{{\mathcal Q}}^T$ and $t \in [1,T]$, $C^{\hat{\pi}'}_t(\mathbf{d}_{[t]}) \leq C^{\hat{\pi}}_t(\mathbf{d}_{[t]})$, proving the desired result.  Let us fix any $s \in [0,T-1]$, and let $\hat{\pi}$ be any policy belonging to $\hat{\Pi}^{s,0}$.  We now define a sequence of functions $\lbrace x'_t, t \in [1,T] \rbrace$, where $x'_1$ is a real constant, and $x'_t$ is a map from $[0,U]_{{\mathcal Q}}^{t-1}$ to ${\mathcal Q}$ for $t \geq 2$.  We will then prove that these functions correspond to a policy belonging to $\hat{\Pi}^{s+1,0}$, which incurs cost no greater than $\hat{\pi}$ in every time period for every sample path of demand.
Let $x'_t = x^{\hat{\pi}}_t$ for all $t \in [1,T] \setminus \lbrace T - s \rbrace$, and $x'_{T - s} = \max(0, x^{\hat{\pi}}_{T- s})$.  We now verify that $\lbrace x'_t, t \in [1,T] \rbrace$ corresponds to some policy belonging to $\hat{\Pi}^{s+1,0}$, and begin by demonstrating that for all $t \in [0,T-1]$ and $\mathbf{d} = (d_1,\ldots,d_t) \in [0,U]_{{\mathcal Q}}^t$, 
$x'_{t+1}(\mathbf{d}) \geq x'_t(\mathbf{d}_{[t-1]}) - d_t$, where $x'_1(d_{[0]}) \stackrel{\Delta}{=} x'_1$, $x'_0(\mathbf{d}_{[-1]}) \stackrel{\Delta}{=} x_0$, and $d_0 \stackrel{\Delta}{=} 0$.  
The assertion is immediate for $t \in [1,T] \setminus \lbrace T - s \rbrace$, as in these cases the property is inherited from $\lbrace x^{\hat{\pi}}_t, t \in [1,T]\rbrace$, and the monotonicty of our construction.  We now verify the assertion for $t = T - s$, by demonstrating that $x'_{T - s + 1}(\mathbf{d}) \geq x'_{T - s}(\mathbf{d}_{[T-s-1]}) - d_{T-s}$, and proceed by a case analysis.  First, suppose $x^{\hat{\pi}}_{T-s}(\mathbf{d}_{[T-s-1]}) \geq 0$.  In this case, the desired property is trivially inherited from $x^{\hat{\pi}}_{T-s}$.  
If not, by construction $x'_{T - s}(\mathbf{d}_{[T-s-1]}) = 0$, and thus $x'_{T - s}(\mathbf{d}_{[T-s-1]}) - d_{T-s} \leq 0$.  However, $\hat{\pi} \in \hat{\Pi}^{s,0}$ implies that $x^{\hat{\pi}}_{T - s + 1}(\mathbf{d}) \geq 0$, demonstrating the desired property.  Combining 
the above with the fact that $x'_t \geq 0$ for all $t \in [T-s+1,T]$ by virtue of $x'_t$ equaling $x^{\hat{\pi}}_t$ for all $t \geq T - s + 1$ and $\hat{\pi}$ belonging to $\hat{\Pi}^{s,0}$, and the fact that $x'_{T-s} \geq 0$  by construction, completes the proof 
that for all $s \in [0,T-1]$, and any policy $\hat{\pi} \in \hat{\Pi}^{s,0}$, there exists a policy $\hat{\pi}' \in \hat{\Pi}^{s + 1,0}$ s.t. $x^{\hat{\pi}'}_t = x^{\hat{\pi}}_t$ for all $t \in [1,T] \setminus \lbrace T - s \rbrace$, and $x^{\hat{\pi}'}_{T - s} = \max(0,x^{\hat{\pi}}_{T - s})$.  Combining with the fact that for all $d \in [0,U]$ and $x < 0$,
$$C(x,d) = b (d + |x|) > b[d - 0]_+ + [0 - d]_+ = b d,$$
further demonstrates that for all vectors $\mathbf{d} = (d_1,\ldots,d_T) \in [0,U]_{{\mathcal Q}}^T$ and $t \in [1,T]$, $C^{\hat{\pi}'}_t(\mathbf{d}_{[t]}) \leq C^{\hat{\pi}}_t(\mathbf{d}_{[t]})$, completing the proof of the desired statement.   Applying the statement inductively, we conclude that for any $\hat{\pi} \in \hat{\Pi}$, there exists a policy $\hat{\pi}' \in \hat{\Pi}^{T,0}$ (i.e. a policy ordering up to at least 0 in every period for all sample paths) which incurs cost at most that incurred by $\hat{\pi}$.
\\\\
We next prove that one can restrict to policies always ordering up to at most $U$.  For $s \in [1,T]$, let $\hat{\Pi}^{s,U}$ denote that subset of $\hat{\Pi}^{T,0}$ consisting of those $\hat{\pi}$ with the following property: for all $t \in [T - s + 1, T]$ and $\mathbf{d} = (d_1,\ldots,d_T) \in [0,U]_{{\mathcal Q}}^T$, $x^{\hat{\pi}}_t(\mathbf{d}_{[t-1]}) \leq U$ whenever $y^{\hat{\pi}}_t(\mathbf{d}_{[t-1]}) \leq U$.  
Also, let $\hat{\Pi}^{0,U} \stackrel{\Delta}{=} \hat{\Pi}^{T,0}$.  We now prove that for all $s \in [0,T-1]$, and any policy $\hat{\pi} \in \hat{\Pi}^{s,U}$, there exists a policy $\hat{\pi}' \in \hat{\Pi}^{s + 1,U}$ s.t. for all vectors $\mathbf{d} = (d_1,\ldots,d_T) \in [0,U]_{{\mathcal Q}}^T$ and $t \in [1,T]$, $C^{\hat{\pi}'}_t(\mathbf{d}_{[t]}) \leq C^{\hat{\pi}}_t(\mathbf{d}_{[t]})$.  Let us fix any $s \in [0,T-1]$, and let $\hat{\pi}$ be any policy belonging to $\hat{\Pi}^{s,U}$.  We now define a sequence of functions $\lbrace x'_t, t \in [1,T] \rbrace$, where $x'_1$ is a real constant, and $x'_t$ is a map from $[0,U]_{{\mathcal Q}}^{t-1}$ to ${\mathcal Q}$ for $t \geq 2$.  We will then prove that these functions correspond to a policy belonging to $\hat{\Pi}^{s+1,U}$, which incurs cost no greater than $\hat{\pi}$ in every time period for every sample path of demand.  Let $x'_t = x^{\hat{\pi}}_t$ for all $t \in [1,T] \setminus \lbrace T - s \rbrace$; and for $\mathbf{d} \in [0,U]_{{\mathcal Q}}^{T-s-1}$ let 
$$
\begin{aligned}
   & x'_{T-s}(\mathbf{d}) \stackrel{\Delta}{=} \begin{cases}
\min\big(U ,  x^{\hat{\pi}}_{T- s}(\mathbf{d})\big) & \text{if}\ y^{\hat{\pi}}_{T- s}(\mathbf{d}) \leq U,\\
   x^{\hat{\pi}}_{T-s}(\mathbf{d}) & \text{else.}
   \end{cases}
\end{aligned}
$$
We now verify that $\lbrace x'_t, t \in [1,T] \rbrace$ corresponds to some policy belonging to $\hat{\Pi}^{s+1,U}$, and begin by demonstrating that for all $t \in [0,T-1]$ and $\mathbf{d} = (d_1,\ldots,d_t) \in [0,U]_{{\mathcal Q}}^t$, 
$x'_{t+1}(\mathbf{d}) \geq x'_t(\mathbf{d}_{[t-1]}) - d_t$.  The assertion is immediate for $t \in [1,T] \setminus \lbrace T - s - 1 \rbrace$, as in these cases the property is inherited from $\lbrace x^{\hat{\pi}}_t, t \in [1,T]\rbrace$, and the monotonicty of our construction.  We now verify the assertion for $t = T - s - 1$, by demonstrating that $x'_{T - s}(\mathbf{d}) \geq x'_{T - s - 1}(\mathbf{d}_{[T-s-1]}) - d_{T-s - 1}$, and proceed by a case analysis.  First, suppose that either $y^{\hat{\pi}}_{T-s}(\mathbf{d}_{[T-s-1]}) > U$, or 
$y^{\hat{\pi}}_{T-s}(\mathbf{d}_{[T-s-1]}) < U$ and $x^{\hat{\pi}}_{T- s}(\mathbf{d}_{[T-s]}) \leq U$.  In this case, the desired property is either trivially true, or inherited from $x^{\hat{\pi}}_{T-s}$.  
If not, it must hold that $y^{\hat{\pi}}_{T- s}(\mathbf{d}_{[T-s]}) \leq U$, and $x^{\hat{\pi}}_{T- s}(\mathbf{d}_{[T-s]}) > U$.  However, by construction this implies that $x'_{T - s - 1}(\mathbf{d}_{T-s-2}) - d_{T-s-1} \leq U$, while $x'_{T-s}(\mathbf{d}_{[T-s-1]}) = U$, demonstrating the desired property.  Combining the above with the fact that $x'_t(\mathbf{d}) \leq U$ whenever $x'_{t-1}(\mathbf{d}_{[t-2]})$ for all $\mathbf{d} \in [0,U]_{{\mathcal Q}}^{t-1}$ and $t \in [T-s+1,T]$ by virtue of $x'_t$ equaling $x^{\hat{\pi}}_t$ for all $t \geq T - s + 1$ and $\hat{\pi}$ belonging to $\hat{\Pi}^{s,U}$, and the fact that $x'_{T-s} \leq U$  by construction, completes the proof that for all $s \in [0,T-1]$, and any policy $\hat{\pi} \in \hat{\Pi}^{s,U}$, there exists a policy $\hat{\pi}' \in \hat{\Pi}^{s + 1,Y}$ s.t. $x^{\hat{\pi}'}_t = x^{\hat{\pi}}_t$ for all $t \in [1,T] \setminus \lbrace T - s \rbrace$, and $x^{\hat{\pi}'}_{T - s}$ is as defined above.  Combining with the fact that for all $d \in [0,U]$ and $x > U$, 
$$C(x,d) = x - d > b[d - U]_+ + [U - d]_+ = U - d,$$
further demonstrates that for all vectors $\mathbf{d} = (d_1,\ldots,d_T) \in [0,U]_{{\mathcal Q}}^T$ and $t \in [1,T]$, $C^{\hat{\pi}'}_t(\mathbf{d}_{[t]}) \leq C^{\hat{\pi}}_t(\mathbf{d}_{[t]})$, completing the proof of the desired statement.   Applying the statement inductively, we conclude that for any $\hat{\pi} \in \hat{\Pi}$, so long as $x_0 \leq U$, there exists a policy $\hat{\pi}' \in \hat{\Pi}^{T,U}$ (i.e. a policy ordering up to a level in $[0,U]$ in every period for all sample paths) which incurs cost at most that incurred by $\hat{\pi}$, completing the proof.  $\Halmos$.
\endproof
\subsection{Proof of Theorem\ \ref{thm-4.2}}
Our proof is quite similar, conceptually, to the proof in \citet{Iyengar} that a robust MDP satisfying the so-called rectangularity property can be solved by dynamic programming.  However, as our problem does not seem to fit precisely into that framework (e.g. since the conditional distribution of demand in a given time period can in principle depend on the entire history of demand up to that time), and as it is well-known that the rectangularity property and related notions can be quite subtle (cf. \citet{XGS}), we include a self-contained proof for completeness.  We will proceed  by a backwards induction, reasoning inductively that w.l.o.g. in period $T$ the conditional distribution of demand can be chosen to be consistent with $\hat{Q}^*$, and thus the policy in period $T$ can be chosen to be consistent with $\hat{\pi}^*$, and thus the conditional distribution of demand in period $T-1$ can be chosen to be consistent with $\hat{Q}^*$, etc.  It is worth noting that one could also proceed by using e.g. minimax theorems, strong duality, saddle-point theory, and/or general results from the theory of stochastic games, and we refer the interested reader to \citet{NE} for closely related results proved using such alternative approaches.  We begin by making several definitions.  For $\pi \in \Pi$, let $\textbf{MAR}^{\pi}_0 = \textbf{MAR}_T$.  For $t \in [1,T-1]$, let us define $\textbf{MAR}^{\pi}_{t}$ as follows.  $Q \in \textbf{MAR}^{\pi}_{t}$ iff $Q \in \textbf{MAR}^{\pi}_{t-1}$, and for all $\mathbf{q} \in \supp(Q_{[T-t]}), Q_{T-t+1 | \mathbf{q}} = \overline{Q}^t_{x^{\pi}_{T-t+1}(\mathbf{q}),q_{T-t}}$.  Also, let $\textbf{MAR}^{\pi}_T$ denote the singleton $\lbrace \hat{Q}^* \rbrace$.
Similarly, let $\Pi_0 = \Pi$, and for $t \in [1, T-1]$, let us define $\Pi_{t}$ as follows.  $\pi \in \Pi_{t}$ iff $\pi \in \Pi_{t-1}$, and for all $\mathbf{q} \in [0,U]_{{\mathcal Q}}^{T - t}$, $x^{\pi}_{T - t + 1}(\mathbf{q}) = \overline{\Phi}^t_{y^{\pi}_{T - t + 1}(\mathbf{q}), q_{T-t}}$.  Also, let $\Pi_T$ denote the singleton $\lbrace \hat{\pi}^* \rbrace$. 

We now state a lemma which follows immediately from our definitions, and formalizes the notion that for any probability measure $Q \in \textbf{MAR}$, and time $t$, there is a unique probability measure belonging to $\textbf{MAR}^{\pi}_t$ which agrees with $Q$ in periods $[1,t-1]$, and can be constructed naturally through an appropriate composition of marginal distributions; as well as the analogous statement for constructing policies belonging to $\Pi_t$ which agree with a given policy in periods $[1,t-1]$.  The associated measures (and policies) can in all cases be constructed through the natural composition of marginals (measurable functions).

\begin{lemma}\label{nonempty2}
Given $\pi \in \Pi$, $Q \in \textbf{MAR}_T$, and $t \in [1,T-1]$, there exists a unique probability measure $Q^{\pi,t} \in \textbf{MAR}^{\pi}_t$ satisfying $Q^{\pi,t}_{[T-t]} = Q_{[T-t]}$.  Similarly, given $\pi \in \Pi$ and $t \in [1,T-1]$, there exists a unique policy $\pi^t \in \Pi_t$ satisfying $x^{\pi}_{[T-t]} = x^{\pi^t}_{[T-t]}$. 
\end{lemma}

We now prove (by induction) a theorem which formalizes the intuitive backwards induction described above, and which will directly imply Theorem\ \ref{thm-4.2}.

\begin{theorem}\label{induct1}
For all $s \in [0,T-1]$, $\pi \in \Pi_s$, and $Q \in \textbf{MAR}^{\pi}_s$,
\begin{equation}\label{ver0}
\bbe_Q[\sum_{t = T - s}^T C^{\pi}_t] = \sum_{\mathbf{q} \in \supp(Q_{[T-s-1]})} \bigg( \sum_{z \in [0,U]_{{\mathcal Q}}} \hat{f}^{s+1}\big(x^{\pi}_{T-s}(\mathbf{q}),z\big) Q_{T-s | \mathbf{q}}(z) \bigg) Q_{[T-s-1]}(\mathbf{q}).
\end{equation}
For all $s \in [0,T-1]$, $\pi \in \Pi_s$, and $Q \in \textbf{MAR}^{\pi}_{s+1}$,
\begin{equation}\label{ver0p5}
\bbe_Q[\sum_{t = T - s}^T C^{\pi}_t] =  \sum_{\mathbf{q} \in \supp(Q_{[T-s-1]})} \hat{g}^{s+1}\big(x^{\pi}_{T-s}(\mathbf{q}),q_{T-s-1}\big) Q_{[T-s-1]}(\mathbf{q}).
\end{equation}
For all $s \in [0,T-1]$ and $\pi \in \Pi_s$, 
\begin{equation}\label{ver1}
\sup_{Q \in \textbf{MAR}} \bbe_Q[ \sum_{t=1}^T C^{\pi}_t ] = \sup_{Q \in \textbf{MAR}^{\pi}_{s+1}} \bbe_Q[ \sum_{t=1}^T C^{\pi}_t ].
\end{equation}
For all $s \in [0,T-1]$, $\pi \in \Pi_{s+1}$, and $Q \in \textbf{MAR}^{\pi}_{s+1}$,
\begin{equation}\label{ver2}
\bbe_Q[\sum_{t = T - s}^T C^{\pi}_t] =  \sum_{\mathbf{q} \in \supp(Q_{[T-s-1]})} \hat{V}^{s+1}\big(y^{\pi}_{T-s}(\mathbf{q}),q_{T-s-1}\big) Q_{[T-s-1]}(\mathbf{q}).
\end{equation}
For all $s \in [0,T-1]$,
\begin{equation}\label{ver3}
\inf_{\pi \in \Pi} \sup_{Q \in \textbf{MAR}} \bbe_Q[ \sum_{t=1}^T C^{\pi}_t ] = \inf_{\pi \in \Pi_{s+1}} \sup_{Q \in \textbf{MAR}} \bbe_Q[ \sum_{t=1}^T C^{\pi}_t ].
\end{equation}
\end{theorem}
\proof{Proof:} We begin with the base case $s = 0$.  (\ref{ver0}) follows immediately from the fact that $\hat{f}^1(x,d) = C(x,d)$.  (\ref{ver0p5}) then follows from the fact that for all $\pi \in \Pi_0$, $Q \in \textbf{MAR}^{\pi}_1$, and $\mathbf{q} \in \supp(Q_{[T-1]})$, 
$Q_{T | \mathbf{q}} = \overline{Q}^1_{x^{\pi}_T(\mathbf{q}),q_{T-1}} \in \hat{Q}^1(x^{\pi}_T(\mathbf{q}),q_{T-1})$, combined with the definition of $\hat{g}^1$ and $\hat{Q}^1$.

We now prove (\ref{ver1}).  Let us fix any $\pi \in \Pi_0$. Note that for any $\epsilon > 0$, there exists $\tilde{Q}_{\epsilon} \in \textbf{MAR}$ s.t. $\bbe[ \sum_{t=1}^T C^{\pi}_t\big(D^{\tilde{Q}_{\epsilon}}_{[t]}) ] > \sup_{Q \in \textbf{MAR}} \bbe_Q[ \sum_{t=1}^T C^{\pi}_t ] - \epsilon.$  Let us fix any such $\epsilon > 0$ and corresponding $\tilde{Q}_{\epsilon} \in \textbf{MAR}$, where we denote $\tilde{Q}_{\epsilon}$ by $\tilde{Q}$ for clarity of exposition.  We now prove that
\begin{equation}\label{ver10b}
\bbe[ \sum_{t=1}^T C^{\pi}_t\big(D^{\tilde{Q}^{\pi,1}}_{[t]}) ] \geq \bbe[ \sum_{t=1}^T C^{\pi}_t\big(D^{\tilde{Q}}_{[t]}) ].
\end{equation}
As $\tilde{Q}^{\pi,1}_{[T-1]} = \tilde{Q}_{[T-1]}$, it suffices by (\ref{ver0}) to demonstrate that for every $\mathbf{q} \in \supp(\tilde{Q}_{[T-1]})$, 
\begin{equation}\label{ver10d}
 \sum_{z \in [0,U]_{{\mathcal Q}}} \hat{f}^1\big( x^{\pi}_T(\mathbf{q}), z \big) \tilde{Q}^{\pi,1}_{T | \mathbf{q}}(z) \geq 
 \sum_{z \in [0,U]_{{\mathcal Q}}} \hat{f}^1\big( x^{\pi}_T(\mathbf{q}), z \big) \tilde{Q}_{T | \mathbf{q}}(z).
\end{equation}
By construction, $\tilde{Q}^{\pi,1} \in \textbf{MAR}^{\pi}_1$, and thus for all $\mathbf{q} \in \supp(\tilde{Q}_{[T-1]})$, $\tilde{Q}^{\pi,1}_{T | \mathbf{q}} = \overline{Q}^1_{x^{\pi}_T(\mathbf{q}),q_{T-1}}
\in \hat{Q}^1(x^{\pi}_T(\mathbf{q}),q_{T-1})$.  Thus by definition,   
the left-hand side of (\ref{ver10d}) equals
\begin{equation}\label{ver10e}
\sup_{Q \in \mathfrak{M}(q_{T-1})} \sum_{z \in [0,U]_{{\mathcal Q}}} \hat{f}^1\big( x^{\pi}_T(\mathbf{q}), z \big) Q(z)\ \ \ =\ \ \ \hat{g}^1\big(x^{\pi}_T(\mathbf{q}),q_{T-1}\big).
\end{equation}
Noting that the martingale property ensures $\tilde{Q}_{T | \mathbf{q}} \in \mathfrak{M}(q_{T-1})$ then completes the proof of (\ref{ver10d}), and (\ref{ver10b}).  Letting $\epsilon \downarrow 0$ completes the proof of (\ref{ver1}).

We now prove (\ref{ver2}).  As $\Pi_1 \subseteq \Pi_0$, it follows from (\ref{ver0p5}), combined with the fact that $\pi \in \Pi_1$ implies $x^{\pi}_{T}(\mathbf{q}) = \overline{\Phi}^1_{y^{\pi}_T(\mathbf{q}), q_{T-1}}$, that 
$$
\bbe_Q[C^{\pi}_T] =
\sum_{\mathbf{q} \in \supp(Q_{[T-1]})} \hat{g}^1\big(\overline{\Phi}^1_{y^{\pi}_T(\mathbf{q}), q_{T-1}},q_{T-1}\big) Q_{[T-1]}(\mathbf{q}).
$$
As $\overline{\Phi}^1_{y^{\pi}_T(\mathbf{q}), q_{T-1}} \in \hat{\Phi}^1\big(y^{\pi}_T(\mathbf{q}), q_{T-1}\big)$, (\ref{ver2}) then follows from definitions.

We now prove (\ref{ver3}).  For any $\epsilon > 0$, there exists $\tilde{\pi}_{\epsilon} \in \Pi_0$ s.t.
$
\inf_{\pi \in \Pi} \sup_{Q \in \textbf{MAR}} \bbe_Q[ \sum_{t=1}^T C^{\pi}_t ]
> 
\sup_{Q \in \textbf{MAR}} \bbe_Q[ \sum_{t=1}^T C^{\tilde{\pi}_{\epsilon}}_t ] - \epsilon
.$
Let us fix any such $\epsilon > 0$ and corresponding $\tilde{\pi}_{\epsilon} \in \Pi_0$, where we denote $\tilde{\pi}_{\epsilon}$ by $\tilde{\pi}$ for clarity of exposition.
We now prove that
\begin{equation}\label{ver30000a}
\sup_{Q \in \textbf{MAR}} \bbe_Q[ \sum_{t=1}^T C^{\tilde{\pi}}_t ]
\geq
\sup_{Q \in \textbf{MAR}} \bbe_Q[ \sum_{t=1}^T C^{\tilde{\pi}^1}_t ].
\end{equation}
As $\Pi_1 \subseteq \Pi_0$, by (\ref{ver1}) it suffices to demonstrate that 
\begin{equation}\label{ver30a}
\sup_{Q \in \textbf{MAR}^{\tilde{\pi}}_1} \bbe_Q[ \sum_{t=1}^T C^{\tilde{\pi}}_t ]
\geq
\sup_{Q \in \textbf{MAR}^{\tilde{\pi}^1}_1} \bbe_Q[ \sum_{t=1}^T C^{\tilde{\pi}^1}_t ].
\end{equation}
Note that $Q \in \textbf{MAR}^{\tilde{\pi}}_1$, combined with (\ref{ver0p5}), implies that
$$\bbe_Q[ \sum_{t=1}^T C^{\tilde{\pi}}_t ]
=
\bbe_Q[ \sum_{t=1}^{T-1} C^{\tilde{\pi}}_t ]
+
\sum_{\mathbf{q} \in \supp(Q_{[T-1]})} \hat{g}^{1}\big(x^{\tilde{\pi}}_{T}(\mathbf{q}),q_{T-1}\big) Q_{[T-1]}(\mathbf{q}).
$$
$\tilde{\pi} \in \Pi_0$ implies that for all $Q \in \textbf{MAR}^{\tilde{\pi}}_1$ and $\mathbf{q} \in [0,U]_{{\mathcal Q}}^{T-1}$, $x^{\tilde{\pi}}_{T}(\mathbf{q}) \geq y^{\tilde{\pi}}_{T}(\mathbf{q})$ and $x^{\tilde{\pi}}_{T}(\mathbf{q}) \in [0,U]_{{\mathcal Q}}$.  By combining the above, we conclude from definitions that
$Q \in \textbf{MAR}^{\tilde{\pi}}_1$ implies
\begin{eqnarray}
\bbe_Q[ \sum_{t=1}^T C^{\tilde{\pi}}_t ]
&\geq&
\bbe_Q[ \sum_{t=1}^{T-1} C^{\tilde{\pi}}_t ]
+
\sum_{\mathbf{q} \in \supp(Q_{[T-1]})} \hat{V}^{1}\big(y^{\tilde{\pi}}_{T}(\mathbf{q}),q_{T-1}\big) Q_{[T-1]}(\mathbf{q})\nonumber
\\&=&
\sum_{\mathbf{q} \in \supp(Q_{[T-1]})}
\bigg(
\sum_{t=1}^{T - 1} C\big(x^{\tilde{\pi}}_t(\mathbf{q}_{[t-1]}), q_t\big) + 
\hat{V}^{1}\big(y^{\tilde{\pi}}_{T}(\mathbf{q}),q_{T-1}\big) 
\bigg)
Q_{[T-1]}(\mathbf{q}).\label{ver300a}
\end{eqnarray}
Alternatively, if $Q \in \textbf{MAR}^{\tilde{\pi}^1}_1$, it follows from (\ref{ver2}), and the fact that $x^{\tilde{\pi}}_{[T-1]} = x^{\tilde{\pi}^1}_{[T-1]}$ and thus $y^{\tilde{\pi}^1}_T = y^{\tilde{\pi}}_T$, that
\begin{equation}\label{ver300b}
\bbe_Q[ \sum_{t=1}^T C^{\tilde{\pi}^1}_t ]
=
\sum_{\mathbf{q} \in \supp(Q_{[T-1]})}
\bigg(
\sum_{t=1}^{T - 1} C\big(x^{\tilde{\pi}}_t(\mathbf{q}_{[t-1]}), q_t\big) + 
\hat{V}^{1}\big(y^{\tilde{\pi}}_{T}(\mathbf{q}),q_{T-1}\big) 
\bigg)
Q_{[T-1]}(\mathbf{q}).
\end{equation}
For any $\epsilon_2 > 0$, there exists $\tilde{Q}_{\epsilon_2} \in \textbf{MAR}^{\tilde{\pi}^1}_1$ (which we denote simply as $\tilde{Q}$) s.t.
\begin{equation}\label{ver300c}
\bbe_{\tilde{Q}}[ \sum_{t=1}^T C^{\tilde{\pi}^1}_t ]
>
\sup_{Q \in \textbf{MAR}^{\tilde{\pi}^1}_1} \bbe_Q[ \sum_{t=1}^T C^{\tilde{\pi}^1}_t ]
- \epsilon_2.
\end{equation}
Lemma\ \ref{nonempty2} implies the existence of $\tilde{Q}^{\tilde{\pi},1} \in \textbf{MAR}^{\tilde{\pi}}_1$ s.t. $\tilde{Q}^{\tilde{\pi},1}_{[T-1]} = \tilde{Q}_{[T-1]}$.  Trivially,
$\sup_{Q \in \textbf{MAR}^{\tilde{\pi}}_1} \bbe_Q[ \sum_{t=1}^T C^{\tilde{\pi}}_t ] \geq \bbe_{\tilde{Q}^{\tilde{\pi},1}}[ \sum_{t=1}^T C^{\tilde{\pi}}_t]$, and it 
follows from (\ref{ver300a}), (\ref{ver300b}), and (\ref{ver300c}) that
\begin{eqnarray*}
\bbe_{\tilde{Q}^{\tilde{\pi},1}}[ \sum_{t=1}^T C^{\tilde{\pi}}_t] &\geq& \sum_{\mathbf{q} \in \supp(\tilde{Q}_{[T-1]})}
\bigg(
\sum_{t=1}^{T - 1} C\big(x^{\tilde{\pi}}_t(\mathbf{q}_{[t-1]}), q_t\big) + 
\hat{V}^{1}\big(y^{\tilde{\pi}}_{T}(\mathbf{q}),q_{T-1}\big) 
\bigg)
\tilde{Q}_{[T-1]}(\mathbf{q})
\\&=& \bbe_{\tilde{Q}}[ \sum_{t=1}^T C^{\tilde{\pi}^1}_t ]
\ \ \ >\ \ \ 
\sup_{Q \in \textbf{MAR}^{\tilde{\pi}^1}_1} \bbe_Q[ \sum_{t=1}^T C^{\tilde{\pi}^1}_t ]
- \epsilon_2.
\end{eqnarray*}
Letting $\epsilon_2 \downarrow 0$ completes the proof of (\ref{ver30a}) and (\ref{ver30000a}).  Letting $\epsilon \downarrow 0$ completes the proof of (\ref{ver3}). 
\\\\\indent We now proceed by induction.  Suppose that (\ref{ver0}) - (\ref{ver3}) hold for all $s' \in [0, s]$ for some $s \leq T - 2$.  We now prove that (\ref{ver0}) - (\ref{ver3}) also  hold for $s + 1$, and begin by verifying (\ref{ver0}).
Note that for all $\pi \in \Pi_{s+1}$ and $Q \in \textbf{MAR}^{\pi}_{s+1}$, 
\begin{equation}\label{ver0go1}
\bbe_Q[C^{\pi}_{T-s-1}] =
\sum_{\mathbf{q} \in \supp(Q_{[T-s-2]})} 
\bigg(\sum_{z \in [0,U]_{{\mathcal Q}}}
C(x^{\pi}_{T-s-1}(\mathbf{q}),z) 
Q_{T-s-1|\mathbf{q}}(z)\bigg)
Q_{[T-s-2]}(\mathbf{q}).
\end{equation}
Combining with (\ref{ver2}), applied to $s$, and the fact that 
 $y^{\pi}_{T-s}(\mathbf{q}:z) = x^{\pi}_{T-s-1}(\mathbf{q}) - z$, we conclude
\begin{equation}\label{ver0go2}
\bbe_Q[\sum_{t =T - s}^T C^{\pi}_t]
=
\sum_{\mathbf{q} \in \supp(Q_{[T-s-2]})}
\bigg(
\sum_{z \in [0,U]_{{\mathcal Q}}}
\hat{V}^{s+1}\big(x^{\pi}_{T-s-1}(\mathbf{q}) - z,z\big) 
Q_{T-s-1 | \mathbf{q}}(z)
\bigg)
Q_{[T-s-2]}(\mathbf{q}).
\end{equation}
Combining (\ref{ver0go1}) - (\ref{ver0go2}) with the definition of $\hat{f}^{s+2}$ completes the proof of (\ref{ver0}).

We now prove (\ref{ver0p5}).  It follows from (\ref{ver0}), applied to $s+1$, and the fact that $\textbf{MAR}^{\pi}_{s+2} \subseteq \textbf{MAR}^{\pi}_{s+1}$, that for all $\pi \in \Pi_{s+1}$ and $Q \in \textbf{MAR}^{\pi}_{s+2}$, 
$$\bbe_Q[\sum_{t = T-s-1}^T C^{\pi}_t] = \sum_{\mathbf{q} \in \supp(Q_{[T-s-2]})} \bigg( \sum_{z \in [0,U]_{{\mathcal Q}}} \hat{f}^{s+2}\big(x^{\pi}_{T-s-1}(\mathbf{q}),z\big) Q_{T-s-1 | \mathbf{q}}(z) \bigg) Q_{[T-s-2]}(\mathbf{q}).
$$
As $Q \in \textbf{MAR}^{\pi}_{s+2}$, $\mathbf{q} \in \supp(Q_{[T-s-2]})$ implies $Q_{T-s-1 | \mathbf{q}} = \overline{Q}^{s+2}_{x^{\pi}_{T-s-1}(\mathbf{q}),q_{T-s-2}} \in \hat{Q}^{s+2}(x^{\pi}_{T-s-1}(\mathbf{q}),q_{T-s-2}),$
(\ref{ver0p5}) then follows from the definitions of $\hat{g}^{s+2}$ and $\hat{Q}^{s+2}$.

We now prove (\ref{ver1}).  As $\Pi_{s+1} \subseteq \Pi_s$, it follows from (\ref{ver1}), applied to $s$, that it suffices to demonstrate that for all $\pi \in \Pi_{s+1}$, 
\begin{equation}\label{ver1induct1}
\sup_{Q \in \textbf{MAR}^{\pi}_{s+1}} \bbe_Q[ \sum_{t=1}^T C^{\pi}_t ] = \sup_{Q \in \textbf{MAR}^{\pi}_{s+2}} \bbe_Q[ \sum_{t=1}^T C^{\pi}_t ].
\end{equation}
Let us fix any $\pi \in \Pi_{s+1}$.  For any $\epsilon > 0$, there exists $\tilde{Q}_{\epsilon} \in \textbf{MAR}^{\pi}_{s+1}$ s.t. $\bbe[ \sum_{t=1}^T C^{\pi}_t\big(D^{\tilde{Q}_{\epsilon}}_{[t]}) ] > \sup_{Q \in \textbf{MAR}^{\pi}_{s+1}} \bbe_Q[ \sum_{t=1}^T C^{\pi}_t ] - \epsilon.$  Let us fix any such $\epsilon > 0$ and corresponding $\tilde{Q}_{\epsilon} \in \textbf{MAR}^{\pi}_{s+1}$, where we denote $\tilde{Q}_{\epsilon}$ by $\tilde{Q}$ for clarity of exposition.  We now prove that
\begin{equation}\label{ver1induct1b}
\bbe[ \sum_{t=1}^T C^{\pi}_t\big(D^{\tilde{Q}^{\pi,s+2}}_{[t]}) ] \geq \bbe[ \sum_{t=1}^T C^{\pi}_t\big(D^{\tilde{Q}}_{[t]}) ].
\end{equation}
As $\tilde{Q}^{\pi,s+2}_{[T-s-2]} = \tilde{Q}_{[T-s-2]}$, by (\ref{ver0}), applied to $s+1$, it suffices  to demonstrate that for every $\mathbf{q} \in \supp(\tilde{Q}_{[T-s-2]})$, 
\begin{equation}\label{ver1induct1d}
\sum_{z \in [0,U]_{{\mathcal Q}}} \hat{f}^{s+2}\big( x^{\pi}_{T-s-1}(\mathbf{q}), z \big) \tilde{Q}^{\pi,s+2}_{T-s-1 | \mathbf{q}}(z)
\geq 
\sum_{z \in [0,U]_{{\mathcal Q}}} \hat{f}^{s+2}\big( x^{\pi}_{T-s-1}(\mathbf{q}), z \big) \tilde{Q}_{T-s-1 | \mathbf{q}}(z).
\end{equation}
By construction, $\tilde{Q}^{\pi,s+2} \in \textbf{MAR}^{\pi}_{s+2}$, and thus for all $\mathbf{q} \in \supp(\tilde{Q}_{[T-s-2]})$, $\tilde{Q}^{\pi,s+2}_{T-s-1 | \mathbf{q}} 
= \overline{Q}^{s+2}_{x^{\pi}_{T-s-1}(\mathbf{q}),q_{T-s-2}} \in \hat{Q}^{s+2}
(x^{\pi}_{T-s-1}(\mathbf{q}),q_{T-s-2})$.  Thus by definition,   
the left-hand side of (\ref{ver1induct1d}) equals
\begin{equation}\label{ver1induct1e}
\sup_{Q \in \mathfrak{M}(q_{T-s-2})} \sum_{z \in [0,U]_{{\mathcal Q}}} \hat{f}^{s+2}\big( x^{\pi}_{T-s-1}(\mathbf{q}), z \big) Q(z)\ \ \ =\ \ \ \hat{g}^{s+2}\big(x^{\pi}_{T-s-1}(\mathbf{q}),q_{T-1}\big).
\end{equation}
Noting that the martingale property ensures $\tilde{Q}_{T-s-1 | \mathbf{q}} \in \mathfrak{M}(q_{T-s-2})$ then completes the proof of (\ref{ver1induct1d}), and (\ref{ver1induct1b}).  Letting $\epsilon \downarrow 0$ further completes the proof of (\ref{ver1}).

We now prove (\ref{ver2}).  As $\Pi_{s+2} \subseteq \Pi_{s+1}$, it follows from (\ref{ver0p5}), applied to $s+1$, and the fact that 
$\pi \in \Pi_{s+2}$ implies $x^{\pi}_{T-s-1}(\mathbf{q}) = \overline{\Phi}^{s+2}_{y^{\pi}_{T-s-1}(\mathbf{q}), q_{T-s-2}},$ that
$$
\bbe_Q[\sum_{t = T - s - 1}^T C^{\pi}_t] = 
\sum_{\mathbf{q} \in \supp(Q_{[T-s-2]})} \hat{g}^{s+2}\big(\overline{\Phi}^{s+2}_{y^{\pi}_{T-s-1}(\mathbf{q}),q_{T-s-2}}\big) Q_{[T-s-2]}(\mathbf{q}).
$$
As 
$\overline{\Phi}^{s+2}_{y^{\pi}_{T-s-1}(\mathbf{q}), q_{T-s-1}} \in \hat{\Phi}^{s+2}\big(y^{\pi}_{T-s-1}(\mathbf{q}), q_{T-s-1}\big)$, (\ref{ver2}) then follows from definitions.

We now prove (\ref{ver3}).  First, observe that by (\ref{ver3}), applied to $s$, it suffices to prove that 
\begin{equation}\label{ver3induct1a}
\inf_{\pi \in \Pi_{s+1}} \sup_{Q \in \textbf{MAR}} \bbe_Q[ \sum_{t=1}^T C^{\pi}_t ] = \inf_{\pi \in \Pi_{s+2}} \sup_{Q \in \textbf{MAR}} \bbe_Q[ \sum_{t=1}^T C^{\pi}_t ].
\end{equation}
For any $\epsilon > 0$, there exists $\tilde{\pi}_{\epsilon} \in \Pi_{s+1}$ s.t. $
\inf_{\pi \in \Pi_{s+1}} \sup_{Q \in \textbf{MAR}} \bbe_Q[ \sum_{t=1}^T C^{\pi}_t ]
> 
\sup_{Q \in \textbf{MAR}} \bbe_Q[ \sum_{t=1}^T C^{\tilde{\pi}_{\epsilon}}_t ] - \epsilon
.$  Let us fix any such $\epsilon > 0$ and corresponding $\tilde{\pi}_{\epsilon} \in \Pi_{s+1}$, where we denote $\tilde{\pi}_{\epsilon}$ by $\tilde{\pi}$ for clarity of exposition.
We now prove that
$$
\sup_{Q \in \textbf{MAR}} \bbe_Q[ \sum_{t=1}^T C^{\tilde{\pi}}_t ]
\geq
\sup_{Q \in \textbf{MAR}} \bbe_Q[ \sum_{t=1}^T C^{\tilde{\pi}^{s+2}}_t ].
$$
As $\Pi_{s+2} \subseteq \Pi_{s+1}$, by (\ref{ver1}), applied to $s+1$, it suffices to demonstrate that 
\begin{equation}\label{ver3induct1a1}
\sup_{Q \in \textbf{MAR}^{\tilde{\pi}}_{s+2}} \bbe_Q[ \sum_{t=1}^T C^{\tilde{\pi}}_t ]
\geq
\sup_{Q \in \textbf{MAR}^{\tilde{\pi}^{s+2}}_{s+2}} \bbe_Q[ \sum_{t=1}^T C^{\tilde{\pi}^{s+2}}_t ].
\end{equation}
Note that $Q \in \textbf{MAR}^{\tilde{\pi}}_{s+2}$, combined with (\ref{ver0p5}), applied to $s+1$, implies that
$$\bbe_Q[ \sum_{t=1}^T C^{\tilde{\pi}}_t ]
=
\bbe_Q[\sum_{t=1}^{T-s-2} C^{\tilde{\pi}}_t ]
+
\sum_{\mathbf{q} \in \supp(Q_{[T-s-2]})} \hat{g}^{s+2}\big(x^{\tilde{\pi}}_{T-s-1}(\mathbf{q}),q_{T-s-2}\big) Q_{[T-s-2]}(\mathbf{q}).
$$
$\tilde{\pi} \in \Pi_{s+1}$ implies that for all $Q \in \textbf{MAR}^{\tilde{\pi}}_{s+2}$ and $\mathbf{q} \in [0,U]^{T-s-2}_{{\mathcal Q}}$, $x^{\tilde{\pi}}_{T-s-1}(\mathbf{q}) \geq y^{\tilde{\pi}}_{T-s-1}(\mathbf{q})$ and $x^{\tilde{\pi}}_{T-s-1}(\mathbf{q}) \in [0,U]_{{\mathcal Q}}$.  By combining the above, we conclude from definitions that
$Q \in \textbf{MAR}^{\tilde{\pi}}_{s+2}$ implies $\bbe_Q[ \sum_{t=1}^T C^{\tilde{\pi}}_t ]$ is at least
$$\bbe_Q[ \sum_{t=1}^{T-s-2} C^{\tilde{\pi}}_t ]
+
\sum_{\mathbf{q} \in \supp(Q_{[T-s-2]})} \hat{V}^{s+2}\big(y^{\tilde{\pi}}_{T-s-1}(\mathbf{q}),q_{T-s-2}\big) Q_{[T-s-2]}(\mathbf{q}),
$$
which is itself equal to
\begin{equation}\label{ver3induct1a2}
\sum_{\mathbf{q} \in \supp(Q_{[T-s-2]})}
\bigg(
\sum_{t=1}^{T-s-2} C\big(x^{\tilde{\pi}}_t(\mathbf{q}_{[t-1]}), q_t\big) + 
\hat{V}^{s+2}\big(y^{\pi}_{T-s-1}(\mathbf{q}),q_{T-s-2}\big) 
\bigg) Q_{[T-s-2]}(\mathbf{q}).
\end{equation}
Alternatively, if $Q \in \textbf{MAR}^{\tilde{\pi}^{s+2}}_{s+2}$, it follows from (\ref{ver2}), applied to $s+1$, and the fact that $x^{\tilde{\pi}}_{[T-s-2]} = x^{\tilde{\pi}^{s+2}}_{[T-s-2]}$, that $\bbe_Q[ \sum_{t=1}^T C^{\tilde{\pi}^{s+2}}_t ]$ equals
$$\bbe_Q[ \sum_{t=1}^{T-s-2} C^{\tilde{\pi}}_t ]
+
\sum_{\mathbf{q} \in \supp(Q_{[T-s-2]})} \hat{V}^{s+2}\big(y^{\pi}_{T-s-1}(\mathbf{q}),q_{T-s-2}\big) Q_{[T-s-2]}(\mathbf{q}),$$
which is itself equal to
\begin{equation}\label{ver3induct1b2}
\sum_{\mathbf{q} \in \supp(Q_{[T-s-2]})}
\bigg(
\sum_{t=1}^{T-s-2} C\big(x^{\tilde{\pi}}_t(\mathbf{q}_{[t-1]}), q_t\big) + 
\hat{V}^{s+2}\big(y^{\pi}_{T-s-1}(\mathbf{q}),q_{T-s-2}\big) 
\bigg)
Q_{[T-s-2]}(\mathbf{q}).
\end{equation}
For any $\epsilon_2 > 0$, there exists $\tilde{Q}_{\epsilon_2} \in \textbf{MAR}^{\tilde{\pi}^{s+2}}_{s+2}$ (which we denote simply as $\tilde{Q}$) s.t.
\begin{equation}\label{ver3induct1c2}
\bbe_{\tilde{Q}}[ \sum_{t=1}^T C^{\tilde{\pi}^{s+2}}_t ]
>
\sup_{Q \in \textbf{MAR}^{\tilde{\pi}^{s+2}}_{s+2}} \bbe_Q[ \sum_{t=1}^T C^{\tilde{\pi}^{s+2}}_t ]
- \epsilon_2.
\end{equation}
Lemma\ \ref{nonempty2} implies the existence of $\tilde{Q}^{\tilde{\pi},s+2} \in \textbf{MAR}^{\tilde{\pi}}_{s+2}$ s.t. 
$\tilde{Q}^{\tilde{\pi},s+2}_{[T-s-2]} = \tilde{Q}_{[T-s-2]}$.
Trivially, $\sup_{Q \in \textbf{MAR}^{\tilde{\pi}}_{s+2}} \bbe_Q[ \sum_{t=1}^T C^{\tilde{\pi}}_t ] \geq \bbe_{\tilde{Q}^{\tilde{\pi},s+2}}[ \sum_{t=1}^T C^{\tilde{\pi}}_t]$, which by (\ref{ver3induct1a2}), (\ref{ver3induct1b2}), and (\ref{ver3induct1c2}) is itself at least
\begin{eqnarray*}
\ &\ &\ \ \ \sum_{\mathbf{q} \in \supp(\tilde{Q}^{\tilde{\pi},s+2}_{[T-s-2]})}
\bigg(
\sum_{t=1}^{T-s-2} C\big(x^{\tilde{\pi}}_t(\mathbf{q}_{[t-1]}), q_t\big) + 
\hat{V}^{s+2}\big(y^{\pi}_{T-s-1}(\mathbf{q}),q_{T-s-2}\big) 
\bigg)
\tilde{Q}^{\tilde{\pi},s+2}_{[T-s-2]}(\mathbf{q})
\\&\ &\ \ \ \ \ \ =\ \ \ 
\sum_{\mathbf{q} \in \supp(\tilde{Q}_{[T-s-2]})}
\bigg(
\sum_{t=1}^{T-s-2} C\big(x^{\tilde{\pi}}_t(\mathbf{q}_{[t-1]}), q_t\big) + 
\hat{V}^{s+2}\big(y^{\pi}_{T-s-1}(\mathbf{q}),q_{T-s-2}\big) 
\bigg)
\tilde{Q}_{[T-s-2]}(\mathbf{q})
\\&\ &\ \ \ \ \ \ =\ \ \ \bbe_{\tilde{Q}}[ \sum_{t=1}^T C^{\tilde{\pi}^{s+2}}_t ]
\ \ \ >\ \ \ 
\sup_{Q \in \textbf{MAR}^{\tilde{\pi}^{s+2}}_{s+2}} \bbe_Q[ \sum_{t=1}^T C^{\tilde{\pi}^{s+2}}_t ]
- \epsilon_2.
\end{eqnarray*}
Letting $\epsilon_2 \downarrow 0$ completes the proof of (\ref{ver3induct1a1}) and (\ref{ver3induct1a}).  Letting $\epsilon \downarrow 0$ completes the proof of (\ref{ver3}). 

Combining all of the above completes the desired induction and proof of the theorem.  $\Halmos$.
\endproof

\subsection{Proof of Observation\ \ref{thereiszero}}
\proof{Proof of Observation\ \ref{thereiszero} : }
For a general  $T \geq 1, \mu \in (0,U), x \in \big[\chi^{T}_{\textbf{MAR}}(\mu,U,b) , U \big)$, let $j'$ be the unique index s.t. $\mu \in (A^T_{j'}, A^T_{j' + 1}]$, and $k'$ to be the unique index s.t. $x \in [B^T_{k'}, B^T_{k' + 1})$, where existence and uniqueness follow from definitions and our assumptions.  From definitions (especially of ${\mathfrak q}^T_{x,\mu}$) and our assumptions (especially time-homogeneity of the relevant parameters and the accompanying self-reducibility of the DP equations, as well as the martingale property), to prove the desired result it suffices to prove (for such general $T,\mu,x$) that $x \geq \chi^{T}_{\textbf{MAR}}(\mu,U,b)$ imples : 1. $k' \geq j' + 1$ (i.e. that the first case in the definition of ${\mathfrak q}^T_{x,\mu}$ cannot occur); and 2. $A^T_{k'} \geq x.$  Thus suppose $x \geq \chi^{T}_{\textbf{MAR}}(\mu,U,b).$  Recall from definitions that $\chi^{T}_{\textbf{MAR}}(\mu,U,b) = B^T_{\Gamma^T_{\mu}}$, and $\mu \in \big( A^{T+1}_{\Gamma^T_{\mu} - 1}, A^{T+1}_{\Gamma^T_{\mu}} \big]$.  It follows that $\mu \leq A^{T+1}_{\Gamma^T_{\mu}} \leq A^T_{\Gamma^T_{\mu}}$ by the monotonicity (in T) of $A^T_j$.  Combining with the fact that $\mu \in (A^T_{j'}, A^T_{j' + 1}]$, and the monotonicity (in j) of $A^T_j$, it follows that
\begin{equation}\label{j'is}
\Gamma^T_{\mu} \geq j' + 1.
\end{equation}
Since by definition $\chi^{T}_{\textbf{MAR}}(\mu,U,b) = B^T_{\Gamma^T_{\mu}}$, the fact that $x \geq \chi^{T}_{\textbf{MAR}}(\mu,U,b)$ and $x \in [B^T_{k'}, B^T_{k' + 1})$, combined with the monotonicity (in k) of $B^T_k$, together imply that 
\begin{equation}\label{k'is}
k' \geq \Gamma^T_{\mu}.
\end{equation}
Combining (\ref{j'is}) and (\ref{k'is}) we conclude that $k' \geq j' + 1$, showing the first desired statement.  To complete the proof, it thus suffices to prove that (under the same assumptions)  $A^T_{k'} \geq x.$  Since by construction $x \in [B^T_{k'}, B^T_{k' + 1}),$ and it follows from (\ref{BAcompare}) that $A^T_{k'} \geq B^T_{k'+1}$, combining the above completes the proof.  $\qed$
\endproof

\subsection{Proof of Lemma\ \ref{proveme1}}
\proof{Proof:}[Proof of Lemma\ \ref{proveme1}]
We first prove continuity.  That ${\mathfrak g}^T(x,d)$ is a continuous function of $x$ on $(0,U) \setminus \bigcup_{i = \Gamma^{T-1}_{d}}^{T-1} \lbrace B^T_i \rbrace$, and a right-continuous function of $x$ on $[0,U] \setminus \lbrace U \rbrace$, follows from definitions, and the fact that $F^T_{i}(x,d)$ and $G^T_{i}(x,d)$ are continuous functions of $x$ on $[0,U]$ for all $i$.  It similarly follows that $\lim_{x \uparrow B^T_i} {\mathfrak g}^T(x,d)$ exists for all $i \in [\Gamma^{T-1}_d, T] \setminus \lbrace 0 \rbrace$.  It thus suffices to demonstrate that
$\lim_{x \uparrow B^T_i} {\mathfrak g}^T(x,d)$ equals ${\mathfrak g}^T(B^T_i,d)$ for all $i \in [\Gamma^{T-1}_d, T] \setminus \lbrace 0 \rbrace$.
We treat two cases: $i = \Gamma^{T-1}_d$, and $i \in [\Gamma^{T-1}_d + 1, T]$, and begin with the case $i = \Gamma^{T-1}_{d}$.  By assumption we preclude the case $i = 0$.  Thus suppose $i = \Gamma^{T-1}_d \in [1,T]$.
In this case,
\begin{eqnarray*}
\lim_{x \uparrow B^T_i} {\mathfrak g}^T(x,d) &=& F^T_{i - 1}(B^T_i,d)
\\&=& - b B^T_i + (b + T) B^T_i + \big( T b  - (b + 1) i \big) d.
\\&=& T B^T_i+ \big( T b  - (b + 1) i \big) d.
\end{eqnarray*}
Alternatively,
\begin{eqnarray*}
{\mathfrak g}^T(B^T_i,d) &=& G^T_i(B^T_i,d)
\\&=& (T - \frac{b+T}{A^T_i} d) B^T_i + (T - i) b d
\\&=& T B^T_i + \big( T b - (b + 1) i \big) d.
\end{eqnarray*}
Combining the above completes the proof for this case.  Next, suppose $i \in [\Gamma^{T-1}_d + 1, T]$.  In this case,
\begin{eqnarray*}
\lim_{x \uparrow B^T_i} {\mathfrak g}^T(x,d)&=& G^T_{i-1}(B^T_i,d)
\\&=& (T - \frac{b+T}{A^T_{i - 1}} d) B^T_i + \big(T - (i -1) \big) b d
\\&=& T B^T_i + \big( T b - (b + 1) i \big) d.
\end{eqnarray*}
Alternatively,
\begin{eqnarray*}
{\mathfrak g}^T(B^T_i,d) &=& G^T_i(B^T_i,d)
\\&=& (T - \frac{b+T}{A^T_i} d) B^T_i + (T - i) b d\ \ \ =\ \ \ T B^T_i + \big( T b - (b + 1) i \big) d,
\end{eqnarray*}
completing the proof.
\\\indent We now prove convexity.  As in our proof of continuity, that ${\mathfrak g}^T(x,d)$ is a right-differentiable function of $x$ on $(0,U) \setminus \bigcup_{i = \Gamma^{T-1}_d}^{T} \lbrace B^T_i \rbrace$, with non-decreasing right-derivative on the same set, follows from definitions and the fact that $F^T_{i}(x,d)$ and $G^T_{i}(x,d)$ are linear functions of $x$ on $[0,U]$ for all $i$.  It similarly follows that ${\mathfrak g}^T(x,d)$ is a right-differentiable function of $x$ on $[0,U] \setminus \lbrace U \rbrace$, and that $\lim_{x \uparrow B^T_i} \partial^+_x {\mathfrak g}^T(x,d)$ exists for all $i \in [\Gamma^{T-1}_d,T-1] \setminus \lbrace 0 \rbrace$.  It thus suffices to demonstrate that $\lim_{x \uparrow B^T_i} \partial^+_x {\mathfrak g}^T(x,d) \leq
\partial^+_x {\mathfrak g}^T(B^T_i,d)$ for all $i \in [\Gamma^{T-1}_d, T-1] \setminus \lbrace 0 \rbrace$.  We treat two cases: $i = \Gamma^{T-1}_d$, and $i \in [\Gamma^{T-1}_d + 1, T-1]$, and begin with the case $i = \Gamma^{T-1}_d$.  By assumption we preclude the cases $i = 0,T$.  Thus suppose $i = \Gamma^{T-1}_d \in [1,T-1]$.  Then
\begin{eqnarray*}
\lim_{x \uparrow B^T_i} \partial^+_x {\mathfrak g}^T(x,d) &=& \lim_{x \uparrow B^T_i} \partial^+_x  F^T_{i - 1}(x,d)
\\&=& - b.
\end{eqnarray*}
Alternatively,
\begin{eqnarray*}
\partial^+_x {\mathfrak g}^T(B^T_i,d) &=& \partial^+_x G^T_i(B^T_i,d)
\\&=&  T - \frac{b+T}{A^T_i} d.
\\&\geq& T - \frac{b+T}{A^T_i} A^T_i\ \ \ =\ \ \ - b.
\end{eqnarray*}
Combining the above completes the proof for this case.  Next, suppose $i \in [\Gamma^{T-1}_d + 1, T - 1]$.  Then
\begin{eqnarray*}
\lim_{x \uparrow B^T_i} \partial^+_x {\mathfrak g}^T(x,d) &=& \lim_{x \uparrow B^T_i} \partial^+_x  G^T_{i-1}(x,d)
\\&=& T - \frac{b + T}{A^T_{i - 1}} d.
\end{eqnarray*}
Alternatively,
\begin{eqnarray*}
\partial^+_x {\mathfrak g}^T(B^T_i,d)  &=& \partial^+_x G^T_{i}(B^T_i,d)
\\&=&  T - \frac{b + T}{A^T_i} d.
\end{eqnarray*}
The desired result then follows from the fact that $A^T_i$ is increasing in $i$ and $d$ is non-negative.  Combining the above completes the proof.  $\Halmos$
\endproof
\subsection{Proofs of Lemmas\ \ref{interlace1} and\ \ref{proveme2}}
\proof{Proof:}[Proof of Lemma\ \ref{interlace1}]
The desired result follows directly from the observation that $A^{T+1}_{j+1}$ is trivially strictly less than $A^T_{j+1}$ for all $j \in [0, T-2]$, and $\frac{A^T_j}{A^{T+1}_{j+1}} =  \frac{\frac{j+1}{b+j+1}}{\frac{T}{b+T}} < 1$ for all $j \in [0,T-2]$.  $\Halmos$
\endproof

\proof{Proof:}[Proof of Lemma\ \ref{proveme2}]
Let $i = \Gamma^T_d, j = \Gamma^{T-1}_d$.  From Lemma\ \ref{proveme1}, it suffices to prove that $\partial^+_x {\mathfrak g}^T(x,d) \leq 0$ for all $x < B^T_i$, and $\partial^+_x {\mathfrak g}^T\big(B^T_i,d\big) \geq 0$; or that $\partial^+_x {\mathfrak g}^T(x,d) \leq 0$ for all $x < U$ and $B^T_i = U$.  It follows from Lemma\ \ref{interlace1} that $i \in \lbrace j, j + 1 \rbrace$.  We now proceed by a case analysis.  First, suppose $i = j$.  In that case, $B^T_i = B^T_j$, $j \leq T - 1$, and $B^T_i < U$.  We conclude that for all $x < B^T_i$,
$$
\partial^+_x {\mathfrak g}^T(x,d)\ \ \ =\ \ \ \partial^+_x F^T_{j-1}(x,\mu)\ \ \  =\ \ \ - b.$$
Noting that
\begin{eqnarray*}
\partial^+_x {\mathfrak g}^T\big(B^T_i,d\big) &=& \partial^+_x G^T_i(B^T_i,d)
\\&=& T - \frac{b+T}{A^T_i} d\ \ \ \geq\ \ \ T - \frac{b+T}{A^T_i} A^{T+1}_i\ \ \ =\ \ \ 0
\end{eqnarray*}
completes the proof in this setting.
\\\indent Alternatively, suppose that $i = j + 1$.  In this case, $B^T_i > B^T_j$, and
\begin{eqnarray*}
\lim_{x \uparrow B^T_i} \partial^+_x {\mathfrak g}^T\big(x,d\big) &=& \lim_{x \uparrow B^T_i} \partial^+_x G^T_{i - 1}(x,d)
\\&=& T - \frac{b+T}{A^T_{i-1}} d\ \ \ \leq\ \ \ T - \frac{b+T}{A^T_{i-1}} A^{T+1}_{i-1}\ \ \ =\ \ \ 0.
\end{eqnarray*}
If $B^T_i = U$, the result follows from Lemma\ \ref{proveme1}.  Otherwise,
\begin{eqnarray*}
\partial^+_x {\mathfrak g}^T\big(B^T_i, d\big) &=& \partial^+_x G^T_i(B^T_i , d)
\\&=& T - \frac{b+T}{A^T_i} d\ \ \ \geq\ \ \ T - \frac{b+T}{A^T_i} A^{T+1}_i\ \ \ =\ \ \ 0.
\end{eqnarray*}
Combining the above completes the proof.  $\Halmos$
\endproof

\subsection{Proof of Lemma\ \ref{bareqbar2}}
\proof{Proof:}[Proof of Lemma\ \ref{bareqbar2}]
Let us proceed by showing that for each fixed $x \in [0,U]$, $\overline{{\mathfrak g}}^T(x,d) ={\mathfrak g}^T\bigg( \max\big( \beta^T_d , x - d \big) , d \bigg)$ for all $d \in [0,U]$.  As the equivalence is easily verified for the case $d = 0$, in which case $\overline{{\mathfrak g}}^T(x,d) = {\mathfrak g}^T\bigg( \max\big( \beta^T_d , x - d \big) , d \bigg) = T x$, suppose $d > 0$.
We proceed by a case analysis, beginning with the setting $d \in (0, z^T_x)$.
In this case, $\max\big( \beta^T_d , x - d \big) = x - d$, and
$$x\ \ \ >\ \ \ d + B^T_{\Gamma^T_d}\ \ \ \geq\ \ \ d+ B^T_{\Gamma^{T-1}_d},$$
where the second inequality follows from Lemma\ \ref{interlace1} and the monotonicity (in $i$) of $B^T_i$.
Combining with the easily verified fact that $\overline{F}^T_{j}(x,d) = G^T_{j}(x - d,d)$ for all $j$ completes the proof in this case.
\\\indent Alternatively, suppose $d \in [z^T_x, U]$, which implies that $\max\big( \beta^T_d , x - d \big) = \beta^T_d$, and $d \geq x - \beta^T_d$.  We proceed by a case analysis.  First, suppose $d \in (A^T_j, A^{T+1}_{j+1}]$ for some $j \in [-1, T-2]$.  In this case, Lemma\ \ref{interlace1} implies that $d \in (A^T_j, A^T_{j+1}] \bigcap (A^{T+1}_j, A^{T+1}_{j+1}]$, and
\begin{eqnarray*}
\overline{{\mathfrak g}}^T(x,d) &=& \overline{G}^T_{j}(d).
\\&=&  G^T_{j+1}(B^T_{j+1}, d)\ \ \ =\ \ \ {\mathfrak g}^T(\beta^T_d,d).
\end{eqnarray*}
Alternatively, suppose $d \in (A^{T+1}_{j+1}, A^T_{j+1}]$ for some $j \in [-1,T-2]$.  In this case, Lemma\ \ref{interlace1} implies that $d \in (A^T_{j}, A^T_{j +1}] \bigcap (A^{T+1}_{j+1}, A^{T+1}_{j+2}]$, and
\begin{eqnarray*}
\overline{{\mathfrak g}}^T(x,d) &=& \overline{G}^T_{j+1}(d).
\\&=& G^T_{j+2}(B^T_{j+2}, d)\ \ \ =\ \ \ {\mathfrak g}^T(\beta^T_d,d).
\end{eqnarray*}
Lemma\ \ref{interlace1} implies that this treats all cases.  Combining the above completes the proof.  $\Halmos$
\endproof

\subsection{Proof of Lemma\ \ref{convexconcave1}}
\proof{Proof:}[Proof of Lemma\ \ref{convexconcave1}]
First, let us treat the case $d \in [0, z^T_x)$, and begin by proving continuity.  Right-continuity at 0 when $x \neq 0$ follows from the fact that $\lim_{d \downarrow 0} \overline{F}^T_j(x,d) = T x$ for all $j$, and right-continuity at 0 when $x = 0$ follows from definitions.  That $\overline{{\mathfrak g}}^T(x,d)$ is a continuous function of $d$ on $(0, z^T_x) \setminus
\bigcup_{ j = 1}^{T - 1} \lbrace x - B^T_j \rbrace$, and a left-continuous function of $d$ on $(0, z^T_x)$, follows from the continuity (in $d$) of $\overline{F}^T_j(x,d)$ for all $j$.  It similarly follows that $\lim_{d \downarrow x - B^T_j} \overline{{\mathfrak g}}^T(x,d)$  exists for all $j$ s.t. $x - B^T_j \in (0,z^T_x)$, and it thus suffices to demonstrate that
$\lim_{d \downarrow x - B^T_j} \overline{{\mathfrak g}}^T(x,d)$ equals $\overline{{\mathfrak g}}^T(x, x - B^T_j)$ for all $j \in [1,T-1]$ s.t. $x - B^T_j \in (0, z^T_x)$.  Note that for all such $j$,
\begin{eqnarray*}
\lim_{d \downarrow x - B^T_j} \overline{{\mathfrak g}}^T(x,d) &=& \overline{F}^T_{j-1}(x,x - B^T_j)
\\&=& T x + \big( (b-1) T - b (j-1) - \frac{b+T}{A^T_{j-1}} x \big) (x - B^T_j) + \frac{b+T}{A^T_{j-1}} (x - B^T_j)^2.
\end{eqnarray*}
Alternatively,
\begin{eqnarray*}
\overline{{\mathfrak g}}^T(x, x - B^T_j) &=& \overline{F}^T_j(x,x - B^T_j)
\\&=& T x + \big( (b-1) T - b j - \frac{b+T}{A^T_j} x \big) (x - B^T_j) + \frac{b+T}{A^T_j} (x - B^T_j)^2.
\end{eqnarray*}
It follows that $\overline{{\mathfrak g}}^T(x, x - B^T_j) - \lim_{d \downarrow x - B^T_j} \overline{{\mathfrak g}}^T(x,d)$ equals
\begin{eqnarray*}
\ &\ &\ \big( - b + ( \frac{b+T}{A^T_{j-1}} - \frac{b+T}{A^T_j} ) x \big) (x - B^T_j) -  ( \frac{b+T}{A^T_{j-1}} - \frac{b+T}{A^T_j} )(x - B^T_j)^2
\\&\ &\ \ \ =\ (x - B^T_j) \big( - b + ( \frac{b+T}{A^T_{j-1}} - \frac{b+T}{A^T_j} ) x  -  ( \frac{b+T}{A^T_{j-1}} - \frac{b+T}{A^T_j} )(x - B^T_j) \big)
\\&\ &\ \ \ =\ (x - B^T_j) \big( - b + (b+T) B^T_j ( \frac{1}{A^T_{j-1}} - \frac{1}{A^T_{j}} ) \big)\ \ \ =\ \ \ 0,
\end{eqnarray*}
completing the proof of continuity.
\\\indent We now prove convexity.  Again applying Lemma\ \ref{suffcon1}, it suffices to demonstrate that
$\partial^+_d \overline{{\mathfrak g}}^T(x,d)$ exists and is non-decreasing on $(0,z^T_x\big)$.
Since $\overline{F}^T_j(x,d)$ is a convex quadratic function of $d$ for all $j$, we conclude that:
$\partial^+_d \overline{{\mathfrak g}}^T(x,d)$ exists on $(0,z^T_x)$;
$\partial^+_d \overline{{\mathfrak g}}^T(x,d)$ is non-decreasing on $(0,z^T_x) \setminus \bigcup_{j = 1}^{T - 1} \lbrace x - B^T_j \rbrace$; and
$\lim_{d \uparrow x - B^T_j}\partial^+_d \overline{{\mathfrak g}}^T(x,d)$ exists for all $j \in [1,T-1]$ s.t. $x - B^T_j \in (0,z^T_x)$.  It thus suffices to demonstrate that $\lim_{d \uparrow x - B^T_j}\partial^+_d \overline{{\mathfrak g}}^T(x,d) \leq
\partial^+_d \overline{{\mathfrak g}}^T(x , x - B^T_j)$ for all $j \in [1,T-1]$ s.t. $x - B^T_j \in (0, z^T_x)$.
Note that for any such $j$,
\begin{eqnarray*}
\lim_{d \uparrow x - B^T_j}\partial^+_d \overline{{\mathfrak g}}^T(x,d) &=& \partial^+_d \overline{F}^T_j(x,x - B^T_j)
\\&=& (b-1) T - b j - \frac{b+T}{A^T_j} x  + 2 \frac{b+T}{A^T_j} (x - B^T_j)
\\&=& (b-1) T - (b + 2) j + x \frac{b+T}{A^T_j}.
\end{eqnarray*}
Alternatively, it follows from continuity that
\begin{eqnarray*}
\partial^+_d \overline{{\mathfrak g}}^T(x,x - B^T_j) &=& \partial^+_d \overline{F}^T_{j-1}(x,x - B^T_j)
\\&=& (b-1) T - b (j - 1) - \frac{b+T}{A^T_{j-1}} x  + 2 \frac{b+T}{A^T_{j-1}} (x - B^T_j)
\\&=& (b-1) T - (b + 2) j - b + x \frac{b+T}{A^T_{j-1}}.
\end{eqnarray*}
It follows that
\begin{equation}\label{thediff1}
\partial^+_d \overline{{\mathfrak g}}^T(x,x - B^T_j)
- \lim_{d \uparrow x - B^T_j}\partial^+_d \overline{{\mathfrak g}}^T(x,d)\end{equation} equals
$$- b + x\frac{b+T}{A^T_{j-1}} -  x \frac{b+T}{A^T_j},$$
which will be the same sign as
\begin{equation}\label{showpos1}
- b A^T_{j-1} + x (b + T) - x (b + T) \frac{j}{b+j}\ \ \ =\ \ \
- b A^T_{j-1} + x (b + T) \frac{b}{b + j}.\end{equation}
Noting that $\frac{b+j}{b + T} A^T_{j-1} = B^T_j$, and multiplying through the right-hand side of (\ref{showpos1}) by
$\frac{b+j}{b (b + T)}$, we further conclude that (\ref{thediff1}) will be the same sign as
$x - B^T_j$.  Since by assumption $x - B^T_j \in (0, z^T_x)$, we conclude that $x - B^T_j \geq 0$, completing the proof of continuity and convexity for $d \in (0, z^T_x)$ and right-continuity at $0$.
\\\indent Next, let us treat the case $d \in (z^T_x, U]$, and begin by proving continuity.  That
$\overline{{\mathfrak g}}^T(x,d)$ is a continuous function of $d$ on $(z^T_x, U) \setminus
\bigcup_{ j = 0}^{T - 1} \lbrace A^{T+1}_j \rbrace$, and a left-continuous function of $d$ on $(z^T_x,U]$, follows from the continuity (in $d$) of $\overline{G}^T_j(d)$ for all $j$.  It similarly follows that $\lim_{d \downarrow A^{T+1}_j} \overline{{\mathfrak g}}^T(x,d)$ exist for all $j \in [0,T - 1]$ s.t. $A^{T+1}_j > z^T_x$.  It thus suffices to demonstrate that $\lim_{d \downarrow A^{T+1}_j} \overline{{\mathfrak g}}^T(x,d) = \overline{{\mathfrak g}}^T(x, A^{T+1}_j)$ for all $j \in [0,T - 1]$ s.t. $A^{T+1}_j > z^T_x$.  Note that for any such $j$,
\begin{eqnarray*}
\lim_{d \downarrow A^{T+1}_j} \overline{{\mathfrak g}}^T(x,d) &=& \overline{G}^T_j(A^{T+1}_j)
\\&=& T B^T_{j+1} + \big( T b - (b + 1) (j+1) \big) A^{T+1}_j.
\end{eqnarray*}
Alternatively,
\begin{eqnarray*}
\overline{{\mathfrak g}}^T(x,A^{T+1}_j) &=& \overline{G}^T_{j-1}(A^{T+1}_j)
\\&=& T B^T_j + \big( T b - (b + 1) j \big) A^{T+1}_j.
\end{eqnarray*}
Thus
\begin{eqnarray*}
\overline{{\mathfrak g}}^T(x,A^{T+1}_j) - \lim_{d \downarrow A^{T+1}_j} \overline{{\mathfrak g}}^T(x,d) &=& T (B^T_j - B^T_{j+1}) + (b+1) A^{T+1}_j
\\&=& \frac{T}{b+T}( j A^T_j - (j+1) A^T_{j+1} \big) + (b+1) \frac{T}{b+T} A^T_j
\\&=& \frac{T}{b+T}\big( ( b + j + 1) A^T_j - (j+1) A^T_{j+1} \big)\ \ \ =\ \ \ 0,
\end{eqnarray*}
completing the proof of continuity.
\\\indent We now prove concavity.  Again applying Lemma\ \ref{suffcon1}, it suffices to demonstrate that
$\partial^+_d \overline{{\mathfrak g}}^T(x,d)$ exists and is non-increasing on $(z^T_x, U)$.
Since $\overline{G}^T_j(d)$ is a linear function of $d$ for all $j$, it follows from the piece-wise definition of
$\overline{{\mathfrak g}}^T(x,d)$ that demonstrating the desired concavity is equivalent to showing that
$\partial^+_d \overline{G}^T_j(0)$ is non-increasing in $j$.  Noting that $\partial^+_d \overline{G}^T_j(0) = T b - (b + 1)(j + 1)$, which is trivially decreasing in $j$, completes the proof.
\\\indent Finally, let us prove continuity at $z^T_x$.  We consider two cases, depending on how $\eta_d \stackrel{\Delta}{=} x - d$ comes to go from lying above $\beta^T_d$ to lying below $\beta^T_d$.  This ``crossing" can occur in two ways.  In particular, either $\beta^T_d$ and $\eta_d$ actually intersect, or $z^T_x$ occurs at a jump discontinuity of $\beta^T_d$ and the two functions never truly intersect.  We proceed by a case analysis.  Let $i = \Gamma^T_{z^T_x}$.
\\\indent First, suppose that $\beta^T_d$ and $\eta_d$ actually intersect at $z^T_x$, namely $B^T_i = x - z^T_x$.  If $z^T_x \in \bigcup_{j = -1}^T \lbrace A^{T+1}_j \rbrace$, the proof of right-continuity follows identically to our previous proof of right-continuity at $A^{T+1}_j$ for all $j \in [0,T-1]$ s.t. $A^{T+1}_j > z^T_x$, and we omit the details.  Otherwise, right-continuity at $z^T_x$ follows from definitions.  Either way, we need only demonstrate left-continuity.  Note that $z^T_x \in (x - B^T_{i+1} , x - B^T_i]$, and
\begin{eqnarray*}
\lim_{d \uparrow z^T_x} \overline{{\mathfrak g}}^T(x,d) &=& \overline{F}^T_i(x,z^T_x)
\\&=&  T x + \big( (b-1) T - b i - \frac{b+T}{A^T_i} x \big) (x -
B^T_i) + \frac{b+T}{A^T_i}(x - B^T_i)^2
\\&=&  T B^T_i + \big( T b - (b + 1) i \big) (x - B^T_i).
\end{eqnarray*}
Alternatively,
\begin{eqnarray*}
\overline{{\mathfrak g}}^T(x,z^T_x) &=& \overline{G}^T_{i-1}(x - B^T_i)
\\&=& T B^T_i + \big( T b - (b + 1) i \big)(x - B^T_i),
\end{eqnarray*}
completing the proof of continuity in this case.
\\\indent Alternatively, suppose that $\beta^T_d$ and $\eta_d$ do not truly intersect at $z^T_x$.  In this case, $z^T_x = A^{T+1}_i \in (x - B^T_{i+1}, x - B^T_i]$, and
\begin{eqnarray*}
\lim_{d \uparrow z^T_x} \overline{{\mathfrak g}}^T(x,d) &=& \overline{F}^T_i(x,A^{T+1}_i)
\\&=&  T x + \big( (b-1) T - b i - \frac{b+T}{A^T_i} x \big) A^{T+1}_i + \frac{b+T}{A^T_i} (A^{T+1}_i)^2
\\&=& b (T - i) A^{T+1}_i,
\end{eqnarray*}
where the final equality follows from straightforward algebraic manipulations, the details of which we omit.
Alternatively,
\begin{eqnarray*}
\overline{{\mathfrak g}}^T\big(x,z^T_x) &=&
\overline{G}^T_{i-1}(A^{T+1}_i)
\\&=& T B^T_i + \big( T b - (b + 1) i \big) A^{T+1}_i
\\&=& b (T - i) A^{T+1}_i,
\end{eqnarray*}
where the final equality again follows from straightforward algebraic manipulations.  This completes the proof of left-continuity at $z^T_x$.  The proof of right-continuity in this case follows identically to our previous proof of right-continuity at $A^{T+1}_j$ for all $j \in [0,T-1]$ s.t. $A^{T+1}_j > z^T_x$, and we omit the details.  Combining the above completes the proof of the lemma.  $\Halmos$
\endproof
\subsection{Proof of Lemmas\ \ref{piececon1} and\ \ref{linedom1}}
\proof{Proof:}[Proof of Lemma\ \ref{piececon1}]
The statements regarding continuity follow from Lemma\ \ref{convexconcave1}.  Noting that $d > \alpha^T_x$ implies $d > x$, concavity on $(\alpha^T_x, U)$ follows from (\ref{ztcompx}) and Lemma\ \ref{convexconcave1}.  Let $i = \zeta^{T-1}_x$.  Convexity on $(0, z^{T-1}_x)$ follows from (\ref{ztcompx}) and Lemma\ \ref{convexconcave1}.
Supposing $z^{T-1}_x \notin \lbrace 0, U \rbrace$, it follows from definitions that $z^{T-1}_x \in (A^T_{i-1}, A^T_i]$, where ${\mathcal A}^T_x = A^T_i$.  Combining with the convexity of $\overline{F}^{T-1}_j(x,d)$ and $\overline{G}^{T-1}_j(d)$ for all $j$, to prove the lemma, it suffices to demonstrate that: 1.\ $z^{T-1}_x \notin \bigcup_{j = -1}^{T-1} \lbrace A^T_j \rbrace$ implies that
\begin{equation}\label{show1}
\lim_{d \uparrow z^{T-1}_x} \partial^+_d {\mathfrak f}^T(x, d) \leq \partial^+_d {\mathfrak f}^T( x, z^{T-1}_x );
\end{equation}
and 2.\ $x \notin \bigcup_{j = -1}^{T-1} \lbrace A^T_j \rbrace$ implies
\begin{equation}\label{show2}
\lim_{d \uparrow x} \partial^+_d {\mathfrak f}^T(x, d) \leq \partial^+_d {\mathfrak f}^T(x, x).
\end{equation}
We treat several cases, and assume throughout that $z^{T-1}_x, x \notin \bigcup_{j = -1}^{T-1} \lbrace A^T_j \rbrace$.  First, suppose $z^{T-1}_x = x$.  In this case, it follows from a straightforward contradiction argument that
$\Gamma^T_x = i = 0$, $x \in (0, A^T_0)$, and
\begin{eqnarray*}
\lim_{d \uparrow z^{T-1}_x} \partial^+_d {\mathfrak f}^T(x, d) &=& -1 + \lim_{d \uparrow z^{T-1}_x} \partial^+_d \overline{F}^{T-1}_{0}(x,d)
\\&=& - 1 + (b - 1) (T - 1) - \frac{b + T - 1}{A^{T-1}_0} x + 2 \frac{b + T - 1}{A^{T-1}_0} z^{T-1}_x
\\&=& -1 + (b - 1) (T - 1) + \frac{b + T - 1}{A^{T-1}_0} x
\\&\leq& -1 + (b - 1) (T - 1) + \frac{b + T - 1}{A^{T-1}_0} A^T_0\ \ \ =\ \ \ b T - (b + 1).
\end{eqnarray*}
Alternatively,
\begin{eqnarray*}
\partial^+_d {\mathfrak f}^T( x, z^{T-1}_x )
&=& b + \partial^+_d \overline{G}^{T-1}_{-1}(x)
\\&=& b + (T - 1) b \ \ \ =\ \ \ b T.
\end{eqnarray*}
Combining the above completes the proof of both (\ref{show1}) and (\ref{show2}) in this case.
\\\indent Next, suppose $z^{T-1}_x < x$.  We again consider two cases, depending on how $\eta_d \stackrel{\Delta}{=} x - d$ comes to go from lying above $\beta^{T-1}_d$ to lying below $\beta^{T-1}_d$.  In the case that the two actually intersect, it may be easily verified that $z^{T-1}_x = x - B^{T-1}_i$, and
\begin{eqnarray*}
\lim_{d \uparrow z^{T-1}_x} \partial^+_d {\mathfrak f}^T(x, d)
&=& - 1 + \lim_{d \uparrow z^{T-1}_x} \partial^+_d \overline{F}^{T-1}_{i}(x,d)
\\&=& - 1 + (b - 1)(T - 1) - b i - \frac{b+T-1}{A^{T-1}_i} (z^{T-1}_x + B^{T-1}_i ) + 2 \frac{b + T - 1}{A^{T-1}_i} z^{T-1}_x
\\&=& - 1 + (b - 1)(T - 1) - (b + 1) i + \frac{b+T-1}{A^{T-1}_i} z^{T-1}_x.
\end{eqnarray*}
Alternatively,
\begin{eqnarray*}
\partial^+_d {\mathfrak f}^T\big(x, z^{T-1}_x\big)
&=& - 1 + \partial^+_d \overline{G}^{T-1}_{i - 1}(z^{T-1}_x)
\\&=& - 1 +  (T - 1)b - (b+1) i.
\end{eqnarray*}
It follows that
$$\partial^+_d {\mathfrak f}^T(x, z^{T-1}_x) - \lim_{d \uparrow z^{T-1}_x} \partial^+_d {\mathfrak f}^T(x, d)\ \ \ =\ \ \   T-1 -  \frac{b+T-1}{A^{T-1}_i} z^{T-1}_x,$$
which will be the same sign as
$$A^T_i - z^{T-1}_x\ \ \ \geq\ \ \ A^T_i - A^T_i\ \ \ =\ \ \ 0.
$$
Combining the above completes the proof of (\ref{show1}) in this case.  Furthermore, the proof of (\ref{show2}) follows from the fact that $x \notin \bigcup_{j = -1}^{T-1} \lbrace A^T_j \rbrace$, from which it follows that
$$\partial^+_d {\mathfrak f}^T(x, x) - \lim_{d \uparrow x} \partial^+_d {\mathfrak f}^T(x, d) = b + 1.$$
\\\indent Finally, suppose that $z^{T-1}_x < x$ and $z^{T-1}_x \notin \bigcup_{j = 1}^{T - 2} \lbrace x - B^{T-1}_j \rbrace$.  As $x \notin \bigcup_{j = -1}^{T-1} \lbrace A^T_j \rbrace$, this final case follows nearly identically to the proof of (\ref{show2}) in the previous case, and we omit the details.  Combining all of the above cases completes the proof of the lemma.  $\Halmos$
\endproof

\proof{Proof:}[Proof of Lemma\ \ref{linedom1}]
Note that for any measure $Q \in \mathfrak{M}(\mu)$, it is true that
$$\bbe_Q[f(D)] \leq \bbe_Q[\eta(D)] =  \frac{ f(R) - f(L) }{R - L} \mu + \frac{ R f(L) - L f(R) }{R - L}.$$
But since $q$ has support only on points $d \in [0,U]$ s.t. $f(d) = \eta(d)$, it follows that
$$\bbe_{q}[f(D)] =\bbe_{q}[\eta(D)] = \frac{ f(R) - f(L) }{R - L} \mu + \frac{ R f(L) - L f(R) }{R - L}.$$
Combining the above completes the proof.  $\Halmos$
\endproof

\subsection{Proof of Lemma\ \ref{wearefamily1}}
\proof{Proof:}[Proof of Lemma\ \ref{wearefamily1}]
We first prove that ${\mathcal K}^T(x,d) \geq {\mathfrak f}^T(x,d)$ for all $d \in [0,U]$.  Let $i = \zeta^{T-1}_x, j = \Gamma^{T-1}_x$, and $k = \Upsilon^{T-1}_x$.  That ${\mathcal K}^T(x,0) \geq {\mathfrak f}^T(x,0)$ follows from definitions.   We now prove that ${\mathcal K}^T(x,A^T_l) \geq {\mathfrak f}^T(x,A^T_l)$ for all $l \in [i, j - 1]$.  First, it will be useful to rewrite ${\mathcal K}$ in a more convenient form.  Noting that $\aleph^T_x = A^T_k$,
\begin{eqnarray*}
{\mathfrak f}^T( x, A^T_k ) &=& b ( A^T_k - x ) + \overline{G}^{T-1}_{k - 1} (A^T_k)
\\&=&  b ( A^T_k - x ) + (T - 1) B^{T-1}_k + \big( (T - 1) b - (b + 1) k \big) A^T_k,  \end{eqnarray*}
and
\begin{eqnarray}
{\mathcal K}^T(x,d) &=& \frac{ {\mathfrak f}^T(x,A^T_k) - T x}{A^T_k} d + T x \nonumber
\\&=& \big( Tb - (b +1) k + \frac{(T-1) B^{T-1}_k - (b + T) x}{A^T_k} \big) d + T x \nonumber
\\&=& \big( b (T - k) - \frac{ (b + T) x }{A^T_k} \big) d + T x.\label{zetaiswhat1}
\end{eqnarray}
Combining with the fact that for all $l \in [i,j-1]$ one has $A^T_l \in [z^{T-1}_x, x]$, proving the desired statement is
equivalent to proving that
$$
\big( b (T - k) - \frac{ (b + T) x }{A^T_k} \big) A^T_l + T x
\geq x - A^T_l + (T - 1) B^{T-1}_l + \big( (T - 1) b - (b + 1) l \big) A^T_l,
$$
which is itself equivalent to demonstrating that
\begin{equation}\label{showme2a}
\big( b (l + 1 - k) + 1 - \frac{ (b + T) x }{A^T_k} \big) A^T_l + (T - 1) x \geq 0.
\end{equation}
First, it will be useful to prove that the left-hand side of (\ref{showme2a}),
$$\eta(l) \stackrel{\Delta}{=} \big( b (l + 1 - k) + 1 - \frac{ (b + T) x }{A^T_k} \big) A^T_l + (T - 1) x,$$
is decreasing in $l$, for $l \in [0, j - 1]$.  Indeed, after simplifying, we find that
$$
\eta(l + 1) - \eta(l) = \bigg( b + \frac{b}{b + l + 1}\big( b (l + 1 - k) + 1 - \frac{ (b + T) x }{A^T_k} \big) \bigg) A^T_{l+1},$$
which will be the same sign as
\begin{eqnarray*}
\ &\ & \big( b (l + 2 - k) + l + 2 \big) A^T_k - (b + T) x
\\&\ &\ \ \ \leq \big( b (l + 2 - k) + l + 2 \big) A^T_k - (b + T) B^T_k
\\&\ &\ \ \ = (b + 1) (l + 2 - k) A^T_k.
\end{eqnarray*}
Noting that $\ell+1\leq j-1$ and (\ref{BAcompare2}) implies $j \leq k $ completes the proof of monotonicity, which we now use to complete the proof of (\ref{showme2a}).  In particular, the above monotonicity implies that to prove (\ref{showme2a}), it suffices to prove that $\eta(k - 1) \geq 0$.  Note that
$$\eta(k- 1) = A^T_{k - 1} + \frac{ b (T - 1 - k) - k }{b + k} x.
$$
If $b (T - 1 - k) - k \geq 0$, then trivially $\eta(k - 1) \geq 0$.  Thus suppose
$b (T - 1 - k) - k < 0$.  In this case, $\eta(k - 1) \geq 0$ iff
$$
x \leq \frac{ (b + k) A^T_{k - 1} }{k - b(T - 1 - k)}.
$$
As $x \leq B^T_{k + 1}$, it thus suffices to prove that
$$
B^T_{k + 1} \leq \frac{b + k}{k - b(T - 1 - k)} A^T_{k - 1},
$$
which, dividing both sides by $A^T_{k - 1} (b+k)$ and simplifying, is itself equivalent to proving that
$$
(b + T) k \geq \big(k - b(T - 1 - k)\big)(b + k + 1).$$
Noting that $k \leq T - 1$, and thus $k - b(T - 1 - k) \leq k$, thus completes the desired proof that ${\mathcal K}^T(x,A^T_l) \geq {\mathfrak f}^T(x,A^T_l)$ for all $l \in [i, j - 1]$.
\\\indent We now prove that ${\mathcal K}^T(x,A^T_l) \geq {\mathfrak f}^T(x,A^T_l)$ for all $l \in [j,k]$.
By construction,
\[
{\mathcal K}^T(x,A^T_k) = {\mathfrak f}^T(x,A^T_k),
\]
and thus it suffices to prove the desired claim for $l \in [j,k-1]$.  Note that the degenerate case for which $x = A^T_k$ can be ignored, as in that case  $j = k$.  Thus suppose $x < A^T_k$.  In this case, Lemma\ \ref{piececon1} implies that ${\mathfrak f}^T(x,d)$ is a continuous, concave, piecewise linear function of $d$ on $[A^T_j, A^T_k]$.  As ${\mathcal K}^T(x,d)$ is a linear function of $d$, it follows from the basic properties of concave functions that to prove the desired claim, it suffices to demonstrate that
\begin{equation}\label{user222}
\partial^+_d {\mathcal K}(x, A^T_k) \leq \lim_{d \uparrow A^T_k} \partial^+_d {\mathfrak f}^T(x,d),
\end{equation}
which is equivalent to proving that
\begin{equation}\label{comp00a}
b (T - k) - \frac{ (b + T) x}{A^T_k} \geq
b + \lim_{d \uparrow A^T_k} \partial^+_d \overline{G}^{T-1}_{k - 1}(d).
\end{equation}
It follows from definitions that
$$A^T_k\ \ \ =\ \ \ \frac{b + T}{k} B^T_k\ \ \ \leq\ \ \ \frac{b + T}{k} x.$$
Combining with (\ref{comp00a}), we find that to prove the desired claim, it suffices to demonstrate that
\begin{equation}\label{comp00b}
b (T - k) - k \geq b +  (T - 1) b - (b + 1) k.
\end{equation}
Noting that both sides of (\ref{comp00b}) are equivalent completes the proof.
\\\indent Finally, let us prove that ${\mathcal K}^T(x,A^T_l) \geq {\mathfrak f}^T(x,A^T_l)$ for all $l \in [k + 1, T - 1]$, which will complete the proof.  It again follows from Lemma\ \ref{piececon1} and the basic properties of concave functions that in this case it suffices to demonstrate that
$\partial^+_d {\mathcal K}^T(x,A^T_k) \geq \partial^+_d{\mathfrak f}^T(x,A^T_k)$,
which is equivalent to proving that
\begin{equation}\label{comp01a}
b (T - k) - \frac{ (b + T) x }{A^T_k} \geq
b + \partial^+_d \overline{G}^{T-1}_k(A^T_k).
\end{equation}
It follows from definitions that
\begin{eqnarray*}
A^T_k &=& \frac{k + 1}{b + k + 1} A^T_{k + 1}
\\&=& \frac{b + T}{b + k + 1} B^T_{k + 1}\ \ \ \geq\ \ \ \frac{b + T}{b + k + 1} x.
\end{eqnarray*}
Combining with (\ref{comp01a}), we find that in this case it suffices to demonstrate that
\begin{equation}\label{comp01b}
b (T - k) - (b + k + 1) \geq  b +  (T - 1) b - (b + 1)( k + 1).
\end{equation}
Noting that both sides of (\ref{comp01b}) are equivalent completes the proof.  Combining all of the above with
the piece-wise convexity guaranteed by Lemma\ \ref{piececon1} completes the proof that ${\mathcal K}^T(x,d) \geq {\mathfrak f}^T(x,d)$ for all $d \in [0,U]$.
\\\indent We now prove that for all $l \in [k, T - 2]$, ${\mathcal L}^T_l(x,d) \geq {\mathfrak f}^T(x,d)$ for all $d \in [0,U]$.  Note that ${\mathfrak f}^T(x,d)$ is a concave function on $[A^T_k,U]$, and by construction ${\mathcal L}^T_l(x,d)$ is a line tangent to ${\mathfrak f}^T(x,d)$ at $A^T_l$.  It follows from the basic properties of concave functions that
\[
{\mathcal L}^T_l(x,d) \geq {\mathfrak f}^T(x,d)\ \ \ \textrm{for all}\ \ \ d \in [A^T_k,U].
\]
Combining those same properties with (\ref{user222}) and the basic properties of linear functions, it follows that ${\mathcal L}^T_l(x,d) \geq {\mathcal K}^T(x,d)$ for all $d \in [0,A^T_k].$
Combining all of the above completes the proof.  $\Halmos$
\endproof

\subsection{Proof of Theorems\ \ref{compare1theorem} and\ \ref{compare2theorem}, and Observation\ \ref{compare1theorembb}}
\proof{Proof:}[Proof of Theorem\ \ref{compare1theorem}]
Let $\pi_1$ denote the base-stock policy which always orders up to 0, and $\pi_2$ the base-stock policy which always orders up to $U$.  Note that if $x_0 = 0$, then for any $Q \in \mathcal{M}_{\textbf{GEN}}$, w.p.1 $C^{\pi_1}_t = b d_t$, and $C^{\pi_2}_t = U - d_t$, $t \in [1,T]$, which implies that $\text{Opt}^T_{\textbf{GEN}}(\mu,U,b) \leq \min\left\{T b\mu,\ T (U - \mu) \right\} = \text{Opt}^T_{\textbf{IND}}(\mu,U,b)$.  $\Halmos$ 
\endproof
\proof{Proof:}[Proof of Theorem \ref{compare2theorem}]
From Theorems \ref{thm-4.1} and \ref{mainmart1},
\[
\begin{aligned}
   \frac{\text{Opt}^T_{\textbf{MAR}}(\mu,U,b)}{\text{Opt}^T_{\textbf{IND}}(\mu,U,b)}=\ \ \frac{G_{\Gamma^{T}_{\mu}}^T(\beta^T_{\mu},\mu)}{\min\left\{T b\mu,\ T (U - \mu) \right\}}=\ \ \frac{\Gamma^{T}_{\mu}A^{T+1}_{\Gamma^{T}_{\mu}}+ \left(Tb-(b+1)\Gamma^{T}_{\mu}\right)\mu}{\min\left\{T b\mu,\ T (U - \mu) \right\}}.
\end{aligned}
\]
Combining the above with Lemma\ \ref{rm-lemma1} completes the proof.    $\Halmos$
\endproof
\proof{Proof:}[Proof of Observation\ \ref{compare1theorembb}]
The results of \cite{S-12} (i.e. Theorem\ \ref{thm-4.1}) imply that if $x_0 = U$ and $\frac{\mu}{U} \leq \frac{1}{b+1}$, then the dynamics of Problem\ \ref{dynamic0-indep} at optimality are as follows.  The initial inventory level equals $U$.  In each period, with probability $\frac{\mu}{U}$, the demand equals $U$, and with probability $1 - \frac{\mu}{U}$ equals 0.  Up until the first time that the demand equals $U$, the inventory level will equal $U$, and at the end of each period a holding cost equal to $U$ is incurred.  The first time the demand equals $U$, the inventory level drops to 0, and in that period no cost is incurred.  In all later periods, the inventory level is raised to 0, and if in that period the demand equals $U$, a cost of $b \mu$ is incurred - otherwise no cost is incurred.  It follows that 
the value of Problem\ \ref{dynamic0-indep} equals 
$$\sum_{k=1}^T (1 - \frac{\mu}{U})^{k-1} \frac{\mu}{U} \big( (k-1) U + (T - k) b \mu \big) + (1 - \frac{\mu}{U})^T U T,$$ 
which is itself equal to
$$b \mu T + \frac{U^2}{\mu} - (b + 1)U + (1 - \frac{\mu}{U})^T (1 + b - \frac{U}{\mu}) U.$$
Alternatively, it follows from Theorems\ \ref{thm-4.2} and\ \ref{mainmart1}, combined with the fact that ${\mathfrak q}^T_{U,\mu}(U) = \frac{\mu}{U}, {\mathfrak q}^T_{U,\mu}(0) = 1 - \frac{\mu}{U}$, and the martingale property, that the dynamics of Problem\ \ref{dynamic0-mar} at optimality are as follows.  The initial inventory level equals $U$.  With probability $\frac{\mu}{U}$, the demand in every period is $U$, and the inventory level is raised to $U$ at the start of each period, so no cost is incurred over the entire time horizon.  With probability $1 - \frac{\mu}{U}$, the demand in every period is $0$, and the inventory level in every period is $U$, implying that a cost of $T U$ is incurred over the time horizon.  It follows that the value of Problem\ \ref{dynamic0-mar} equals $(1 - \frac{\mu}{U}) T U = (U - \mu) T$.  Combining the above with a straightforward limiting argument, and the fact that $\frac{\mu}{U} < \frac{1}{b+1}$ implies $\frac{U - \mu}{b \mu} > 1$, completes the proof.  $\Halmos$
\endproof

\subsection{Further interpretation of Observation\ \ref{explicitt2}}\label{furtherobssec}
To further help interpret Observation\ \ref{explicitt2} and the explicit forms for the various quantities, we now show how (for the case $T = 2$) one may derive the same quantities from a certain ``heuristic relaxation" of our original problem.  In doing so, we will be able to see precisely where and why certain quantities arise.  In light of Observation\ \ref{thereiszero}, let us consider the following simplified version of  Problem\ \ref{dynamic0-mar}.  To set up this simplified problem, we reason as follows.  Observation\ \ref{thereiszero} suggests that at optimality a worst-case distribution always has support on 2 points, one on 0 and one on some value which clears the inventory.  Thus we create a simplified problem by ``enforcing" that in the first period, the adversary selects such a distribution.  Namely, if the post-ordering inventory level (in the first period) equals $x$, then in the first period the adversary must select (as a function of $x$) some value $d \in \big[ \max(x,\mu) , U \big]$ and set the demand (in the first period) equal to $d$ w.p. $\frac{\mu}{d}$, and equal to 0 w.p. $1 - \frac{\mu}{d}$.  Note that if the adversary acts in this manner in the first period then the dynamics in the second period will (since the policy-maker can order up to any level as their inventory will have been cleared in the first period) be equivalent to the dynamics of a single-period problem in which the policy-maker may select any starting level in [0,U] and the demand must have mean $D_1$ (either d or 0 depending on the random realized demand in the first period), a minimax problem whose explicit solution is given in Theorem\ \ref{thm-4.1}.  Combining the above, we are led to the following simplified version of Problem\ \ref{dynamic0-mar} : 
$$
\inf_{x \in [0,U]} \sup_{d \in \big[ \max(x,\mu), U \big]} \bigg( (1 - \frac{\mu}{d}) \times (2 x) + \frac{\mu}{d} \times \big( b \times (d - x) + \min(b d , U - d) \big) \bigg).$$
Further relaxing the problem to allow d to take any value in [0,U] (i.e. relaxing the constraint that $d \geq \max(x,\mu)$, which will yield the same insights through a simpler analysis), we are led to the following further-simplified problem : 
\begin{equation}\label{simpleman1}
\inf_{x \in [0,U]} \sup_{d \in [0,U]} \bigg( (1 - \frac{\mu}{d}) \times (2 x) + \frac{\mu}{d} \times \big( b \times (d - x) + \min(b d , U - d) \big) \bigg).
\end{equation}
For a fixed $x \in [0,U]$, let us consider the inner maximization problem.  Noting that $\min(b d , U - d) = b d$ exactly when $d \leq \frac{U}{b + 1}$, we find (after some simple algebra) that the inner problem has optimal value
\begin{equation}\label{simpleman2}
\max\bigg( \sup_{d \in [0, \frac{U}{b+1}]}\big( 2 x + 2 \mu b - \frac{(b + 2) \mu x }{d} \big) , \sup_{d \in [\frac{U}{b+1}, U]}\big( 2 x + (b - 1) \mu + \mu \times \frac{U - (b + 2) x }{d} \big) \bigg).
\end{equation}
Noting that $\sup_{d \in [0, \frac{U}{b+1}]}\big( 2 x + 2 \mu b - \frac{(b + 2) \mu x }{d} \big)$ attains its maximum at $d = \frac{U}{b+1}$, which is equivalent to 
$2 x + (b - 1) \mu + \mu \frac{U - (b + 2) x }{\frac{U}{b+1}},$ we find that the first term within the maximum is superfluous in (\ref{simpleman2}), and that (\ref{simpleman2}) equals
\begin{equation}\label{simpleman3}
\sup_{d \in [\frac{U}{b+1}, U]}\big( 2 x + (b - 1) \mu + \mu \times \frac{U - (b + 2) x }{d} \big).
\end{equation}
Now, we find that the optimal value for $d$ in (\ref{simpleman3}) is determined by the sign of $U - (b + 2) x$, i.e. whether $x \leq \frac{U}{b+2}$ or not.  We have thus so far identified two critical values : $\frac{U}{b+1}$ as determining the value and behavior in the final period (in line with the single-period problem), and $\frac{U}{b+2}$ as determining the sign of a critical coefficient which drives whether a worst-case distribution selects the non-zero value for demand to be ``small" or ``large".  Combining the above with some simple algebra, we conclude the following.  For $x \in [0, \frac{U}{b+2})$, the inner maximization has a maximizer at $d = \frac{U}{b+1}$, and has value (at this maximizer) $2 \mu b + \big(2 - (b+1)(b+2)\frac{\mu}{U}\big) x$.  For $x \in [\frac{U}{b+2}, U]$, the inner maximization has a maximizer at $d = U$, and has value (at this maximizer) $\mu b + \big(2 - (b+2)\frac{\mu}{U}\big) x$.  Thus the optimal choice of $x$ is driven by the sign of $2 - (b+1)(b+2)\frac{\mu}{U}$ and $2 - (b+2)\frac{\mu}{U}$, which coincides exactly with whether $\mu \in [0, \frac{2 U}{(b + 1)(b+2)}]$, $\mu \in (\frac{2 U}{(b + 1)(b+2)}, \frac{2 U}{b + 2}]$, or $\mu \in (\frac{2 U}{b + 2}, U]$.  
\\\indent Indeed, first suppose $\mu \in [0, \frac{2 U}{(b + 1)(b+2)}]$.  Then $\inf_{x \in [0,\frac{U}{b+2}]} \bigg( 2 \mu b + \big(2 - (b+1)(b+2)\frac{\mu}{U}\big) x \bigg) = 2 \mu b$ (attained at $x = 0$), while $\inf_{x \in [\frac{U}{b+2}, U]} \bigg( \mu b + \big(2 - (b+2)\frac{\mu}{U}\big) x \bigg) = \mu b+ \big(2 - (b+2)\frac{\mu}{U}\big) \times \frac{U}{b+2}$ (attained at $x = \frac{U}{b+2}$).  As it follows from some straightforward algebra (the details of which we omit) that $\mu b+ \big(2 - (b+2)\frac{\mu}{U}\big) \times \frac{U}{b+2} \geq 2 \mu b$ for all $\mu \in [0, \frac{2 U}{(b + 1)(b+2)}]$, in this case we conclude that the optimal x equals 0, the optimal d (for that choice of x) equals $\frac{U}{b+1}$, and the optimal value equals $2 \mu b$.
\\\indent Next, suppose $\mu \in (\frac{2 U}{(b + 1)(b+2)}, \frac{2 U}{b + 2}]$.  Then $\inf_{x \in [0,\frac{U}{b+2}]} \bigg( 2 \mu b + \big(2 - (b+1)(b+2)\frac{\mu}{U}\big) x \bigg) = 
2 \mu b + \big(2 - (b+1)(b+2)\frac{\mu}{U}\big) \times \frac{U}{b+2}$ (attained at $x = \frac{U}{b+2}$), while $\inf_{x \in [\frac{U}{b+2}, U]} \bigg( \mu b + \big(2 - (b+2)\frac{\mu}{U}\big) x \bigg) = \mu b+ \big(2 - (b+2)\frac{\mu}{U}\big) \times \frac{U}{b+2}$ (attained at $x = \frac{U}{b+2}$).  In this case we conclude that the optimal x equals $\frac{U}{b+2}$, the optimal d (for that choice of x) equals $U$, and the optimal value equals $\mu (b - 1) + \frac{2 U}{b + 2}$.
\\\indent Finally, suppose $\mu \in (\frac{2 U}{b + 2}, U]$.  Then $\inf_{x \in [0,\frac{U}{b+2}]} \bigg( 2 \mu b + \big(2 - (b+1)(b+2)\frac{\mu}{U}\big) x \bigg) = 
2 \mu b + \big(2 - (b+1)(b+2)\frac{\mu}{U}\big) \times \frac{U}{b+2}$ (attained at $x = \frac{U}{b+2}$), while $\inf_{x \in [\frac{U}{b+2}, U]} \bigg( \mu b + \big(2 - (b+2)\frac{\mu}{U}\big) x \bigg) = \mu b+ \big(2 - (b+2)\frac{\mu}{U}\big) \times U$ (attained at $x = U$).  In this case, it again follows from some straightforward algebra (the details of which we omit) that the optimal x equals U, the optimal d (for that choice of x) equals $U$, and the optimal value equals $2 (U - \mu)$.
\\\indent We note that the above cases and values are perfectly consistent with the cases and values for all quantities appearing in Observation\ \ref{explicitt2}, including 
$\chi^2_{\textbf{MAR}} , \text{Opt}^2_{\textbf{MAR}}, D^2_1$, and $X^2_1$.  In summary, the break-point $\frac{U}{b+1}$ comes from the two possible behaviors in the single-period problem (corresponding to the final period), the break-point $\frac{U}{b+2}$ (for x) determines whether a critical coefficient in nature's inner maximization is positive or negative, hence dictating whether nature selects a ``small" or ``large" value for its non-zero support point.  Finally, the break-points $\frac{2 U}{(b + 1)(b+2)}$ and $\frac{2 U}{b + 2}$ (for $\mu$) determine which coefficients in the policy-maker's outer minimization are positive, also acting as important determinants for the optimal policy.
\\\indent We leave as an interesting direction for future research understanding more formally and generally the connection between such a simplified problem and our true problem.
\subsection{Different behaviors when our assumptions do not hold}\label{nozerohere}
As noted in Section\ \ref{mainsubsec}, the ``obsolescence phenomena", i.e. the property that conditional on the past demand realizations and the current inventory level (under an optimal policy), there always exists a worst-case distribution that assigns a strictly positive probability to zero (and otherwise clears the inventory), is a feature of our particular modeling assumptions, and may not hold under different assumptions.  Here we provide several examples showing that this feature need not hold if one relaxes our modeling assumptions.  The first three examples collectively demonstrate that: 1. if one allows for time-dependent costs, then a worst-case distribution may not put positive probability at zero (first example); 2. if one removes the upper bound on the support, then a worst-case distribution may not even exist (second example); and 3. if one imposes a lower bound on the support, then a worst-case distribution may not put positive probability on this lower bound, and furthermore may not clear the inventory (third example).  Collectively, these findings suggest that extending our framework to more complex models will require several fundamentally new ideas, as much of the structure which allowed for our explicit analysis may no longer hold.  However, our fourth example shows there is indeed hope of extending our results to more general settings, by showing that the obsolescence phenomena again manifests even if holding costs can vary arbitrarily over time, if one enforces (the admittedly strong assumption) that all backlogging costs are 0.
\ \\
\begin{example}  \textbf{Non-stationary and possibly zero cost parameters.}  Here we show that if one relaxes the constraint that the cost parameters are stationary and strictly positive, different behavior arises.  Indeed, let us consider the modification of Problem\ \ref{dynamic0-mar} in which $x_0 = 0, T = 2, h_1 = 1, b_1 = 0, h_2 = 1.$  We let $b_2 > 0, U > 0$ be general, and set $\mu = \frac{U}{b_2 + 1}$.  Namely, the problem is identical to Problem\ \ref{dynamic0-mar}, except that we allow for non-stationary and possibly zero costs.  First, we claim that in this setting, there exists a minimax optimal policy $\pi^* = (x^*_1, x^*_2)$ satisfying $x^*_1 = 0$, and $x^*_2(d_1) = \chi_{\textbf{IND}}(d_1,U,b_2)$.  Indeed, to see this, note that we can derive a lower bound on the value of the associated minimax problem by relaxing the problem s.t. at the start of the second period, the policy-maker may (at zero cost) select any inventory level in $[0,U]$, i.e. the policy-maker may order or dispose of inventory at that time (in the actual problem inventory disposal is infeasible).  It follows from Theorem\ \ref{thm-4.1}, non-negativity of all costs, and the fact that the dynamics in the final period coincide with those of a corresponding single-period problem (with $\mu$ equal to the realized value of $D_1$), that $\pi^*$ is actually optimal for this modified problem.  As $\pi^*$ is also feasible for the original problem, optimality for the original problem immediately follows.  Again applying Theorem\ \ref{thm-4.1}, it follows that there exists a worst-case measure $\hat{Q}^*$ (against $\pi^*$) s.t. for all $d_1 \in [0,U]$, $\hat{Q}^*_{2|d_1}(0) = 1 - \frac{d_1}{U}, \hat{Q}^*_{2|d_1}(U) = \frac{d_1}{U}$.  Furthermore, again applying Theorem\ \ref{thm-4.1}, we find that $\hat{Q}^*_1$ can be taken to be any measure belonging to 
$$\argmax_{Q \in \mathfrak{M}(\mu)} E_Q\big[\min(b_2 \times D_1, U - D_1)\big].$$
A simple application of Jensen's inequality shows that the measure $Q$ which assigns probability 1 to the $d$ maximizing $\min(b_2 \times d, U - d)$ is feasible and optimal, i.e. one can take $\hat{Q}^*_1$ to be the measure s.t. $\hat{Q}^*_1(\frac{U}{b_2 + 1}) = 1$.  Namely, here $\hat{Q}^*_1$ does not put strictly positive probability at 0, i.e. such an obsolescence phenomena does not occur (at least not in the first period).  Note that non-stationarity was critical to this conclusion, since if $b_2$ also equalled zero then the problem would become degenerate with the policy which never orders optimal, and every feasible martingale inducing zero cost.
\end{example}
\ \\
\begin{example}  \textbf{No upper bound on support.}
Here we show that if one relaxes the constraint that the demand must be at most U, different behavior arises.  Let us consider the modification of Problem\ \ref{dynamic0-mar} in which $x_0 = 1, T = 1, h_1 = 1, b_1 = 1, \mu = 1$, but there is no upper bound on the support (equivalently $U = \infty$).  Let $\mathfrak{M}^\infty(\mu)$ denote the set of all probability measures with non-negative support and mean $\mu$.  In this case, for any fixed choice of $x_1 \geq x_0$, the inner maximization (over measures) is equivalent to
\begin{equation}\label{ex2noup}
\sup_{Q\in \mathfrak{M}^\infty(1)}\mathbb{E}_{Q}\left[ [D_1 - x_1]_+ + [x_1 - D_1]_+ \right].
\end{equation}
Note that for any strictly positive real number $y$ and non-negative r.v. $Z$ s.t. $P(Z = 0) < 1,$ it is easily verified that $P\big( [Z - y]_+ + [y - Z]_+ < Z + y \big) > 0$.  Since trivially $[Z - y]_+ + [y - Z]_+ \leq Z + y$ w.p.1, it follows that for any fixed $x_1 \geq x_0$ and $Q \in \mathfrak{M}^\infty(1)$, $\mathbb{E}_{Q}\left[ [D_1 - x_1]_+ + [x_1 - D_1]_+ \right] < 1 + x_1.$  However, note that if (for $M \geq x_1 + 2$) we let $Q^M$ denote the measure s.t. $Q^M(\frac{1}{M}) = 1 - \frac{1}{M}, Q^M(M - 1 + \frac{1}{M}) = \frac{1}{M}$, we find that $Q^M \in \mathfrak{M}^\infty(1)$, and $\mathbb{E}_{Q^M}\left[ [D_1 - x_1]_+ + [x_1 - D_1]_+ \right] = (1 - \frac{1}{M}) \times (x_1 - \frac{1}{M}) + 
\frac{1}{M} \times (M - 1 + \frac{1}{M} - x_1) \geq 1 + x_1 - \frac{2 (x_1 + 1)}{M}$, and thus $\lim_{M \rightarrow \infty} \mathbb{E}_{Q^M}\left[ [D_1 - x_1]_+ + [x_1 - D_1]_+ \right] = 1 + x_1$.  Combining the above, we conclude that for any feasible policy, there does not exist any worst-case distribution.  Furthermore, the unique optimal policy is that which orders nothing in the first period, and there exists a sequence of probability measures which (in the limit) yield a minimax value of 2, such that no distribution in the sequence puts strictly positive probability at 0. 
\end{example}
\ \\
\begin{example}  \textbf{Positive lower bound on support.}
Here we show that if one imposes an additional lower bound constraint on the support of the demand, different behavior arises.  Let us consider the modification of Problem\ \ref{dynamic0-mar} in which $x_0 = U, T = 2, h_1 = b_1  = h_2 = b_2 = 1$.  We let $U > 0, L \in (0,U),$ be general so long as they satisfy the constraints $U < 2 \times L$, and set $\mu = \frac{U + L}{2}.$  Let $\mathfrak{M}^L(\mu)$ denote the set of all probability measures with support a subset of $[L,U]$ and mean $\mu$, and $\textbf{MAR}^L$ denote the collection of all probability measures corresponding to discrete-time martingale sequences $(D_1,D_2)$ s.t. for all $t$, the marginal distribution of $D_t$ belongs to $\mathfrak{M}^L(\mu)$.  Let $\Pi^L$ denote the family of (appropriately adapted) policies that both : 1. never order up to more than U, and 2. always order up to at least L.  Restricting to policies in $\Pi^L$ is without loss of generality, where the proof is analogous to that of Lemma\ \ref{lemma0u}, and we omit the details.
Then we consider the problem $\inf_{\pi \in \Pi^L} \sup_{Q \in \textbf{MAR}^L} \bbe_Q[ \sum_{t=1}^T C^{\pi}_t ].$  We begin by defining some auxiliary functions, which correspond to the solution to the corresponding 1-period problem.  It follows from a straightforward convexity argument (analogous to that used in the proof of Theorem\ \ref{thm-4.1}) and some algebra that
$$
   g^{1,L}(x, \hat{\mu})\stackrel{\Delta}{=}\ \sup_{Q \in \mathfrak{M}^L(\hat{\mu})} \mathbb{E}_{Q}\left[\left|x-D_1\right| \right]=\ \frac{(U - x) \times (\hat{\mu} - L) + (x - L) \times (U - \hat{\mu})}{U - L},$$
where one can always take as a worst-case distribution the measure $Q$ s.t. $Q(U) = \frac{\hat{\mu} - L}{U - L}, Q(L) = 1 - \frac{\hat{\mu} - L}{U - L}$.
Furthermore, $\inf_{x \in [L,U]} g^{1,L}(x, \hat{\mu}) = \min(\hat{\mu} - L, U - \hat{\mu})$, and 
$$x^{*,L}(\hat{\mu}) \stackrel{\Delta}{=} \argmin_{x \in [L,U]}\ g^1(x, \hat{\mu}) =\ \begin{cases}
U & \text{if $\hat{\mu} > \frac{U+L}{2}$};\\
L & \text{if $\hat{\mu} <\frac{U+L}{2}$};\\
[L,U] & \text{if $\hat{\mu} = \frac{U+L}{2}$}.
\end{cases}
$$
Furthermore, using the facts that: 1. $x_0 = U$, 2. $U < 2 \times L$ implies $x_1 - D_1 \leq L$ w.p.1, and 3. the dynamics in the final period coincide with those of a corresponding single-period problem (whose solution is given by $g^{1,L}),$ we find that $\inf_{\pi \in \Pi^L} \sup_{Q \in \textbf{MAR}^L} \bbe_Q[ \sum_{t=1}^T C^{\pi}_t ]$ equals
\begin{equation}\label{test-12}
U - \mu + \sup_{Q \in \mathfrak{M}^L(\mu)}\mathbb{E}_{Q}\left[ \min(D_1 - L, U - D_1) \right],
\end{equation}  
where the policy $\pi^*$ that orders nothing in period 1, and orders up to $x^{*,L}(D_1)$ in period 2 (where any consistent choice may be selected from $[L,U]$ if $D_1 = \frac{U+L}{2}$) will be optimal.  A simple application of Jensen's inequality shows that the measure $Q$ which assigns probability 1 to the $d$ maximizing $\min(d - L, U - d)$ is feasible and optimal, i.e. one can take $Q$ to be the measure s.t. $Q(\frac{U + L}{2}) = 1$.  Summarizing the above, we conclude that a worst-case measure (for $\pi^*$) is the martingale $\hat{Q}^*$ s.t. $\hat{Q}^*_1(\frac{U+L}{2}) = 1$, while $\hat{Q}^*_2(L) = \hat{Q}^*_2(U) = \frac{1}{2}$.
Hence the worst-case measure not only puts no probability at 0 (which would be impossible due to the lower bound), but actually puts no probability even at the lower bound itself (at least in the first round).  The intuition here is that due to the lower bound, it is less appealing for the adversary to leave the policy-maker stuck holding inventory, as some of this inventory is always reduced in later rounds (due to the lower bound).  Interestingly, we note that here, $D_1$ is not large enough to clear the initial inventory, instead only bringing the inventory below the lower bound (in contrast to the behavior when no such lower bound is imposed, e.g. in Observation\ \ref{thereiszero}).
\end{example}
\ \\
\begin{example}  \textbf{Non-stationary cost parameters revisited: a setting where the obsolescence phenomena again manifests.}
Here we show that there is a simple (albeit not completely trivial) setting in which costs can be non-stationary, yet our results (and general insights) regarding obsolescence continue to hold.  Let us consider the modification of Problem\ \ref{dynamic0-mar} in which $x_0 \in (0,U]$ is strictly positive (but otherwise general), $T$ is general, $\lbrace h_i, i = 1,\ldots,T \rbrace$ are strictly positive (but otherwise general), $\mu \in (0,U)$ is general, and the major restriction is that $b_i = 0$ for all $i$.  Namely, there is no back-logging cost whatsoever, and the only incurred costs are holding costs, which again may vary in a general manner over time.  It follows from a straightforward contradiction argument (similar in spirit to the proof of Lemma\ \ref{lemma0u}, the details of which we omit) that in this setting the policy $\pi^*$ which never orders is optimal.  The associated inner maximization (over measures) is thus equivalent to
$$
\sup_{Q \in \textbf{MAR}} \sum_{t=1}^T h_t E_Q[\max(0, x_0 - \sum_{i=1}^t D_i)].$$
For $t \in [1,T]$, let $\mathfrak{M}^{t U}(t \mu)$ denote the set of all probability measures Q with support on $[0, t U]$ such that $E[D^Q] = t \mu$.
As for all $t \in [1,T]$, $\sum_{i=1}^t D_i$ has mean $t \mu$ and support a subset of $[0, t U]$, we conclude (after applying convexity and Theorem\ \ref{thm-4.1}) that for all $Q \in \textbf{MAR}$ and $t \in [1,T]$, 
$$
E_Q[\max(0, x_0 - \sum_{i=1}^t D_i)] \leq \sup_{Q \in \mathfrak{M}^{t U}(t \mu)}\mathbb{E}_{Q}\big[ \max(0, x_0 - D_1) \big] = (1 - \frac{\mu}{U}) \times x_0,$$
where the measure $Q^*$ such that $Q^*(0) = 1 - \frac{\mu t}{U t} = 1 - \frac{\mu}{U}, Q^*(t U) = \frac{\mu}{U}$ belongs to  $\argmax_{Q \in \mathfrak{M}^{t U}(t \mu)}\mathbb{E}_{Q}\big[ \max(0, x_0 - D_1) \big]$.
Noting that under the martingale measure $\hat{Q}^* \in \textbf{MAR}$ such that $P\big(D^{\hat{Q}^*}_i = U\ \textrm{for all}\ i \in [1,T]\big) = \frac{\mu}{U}, P\big(D^{\hat{Q}^*}_i = 0\ \textrm{for all}\ i \in [1,T]\big) = 1 - \frac{\mu}{U}$ it holds (for all t) that $P(\sum_{i=1}^t D_i = 0) = 1 - \frac{\mu}{U}, P(\sum_{i=1}^t D_i = t U) = \frac{\mu}{U}$, we conclude that 
(under optimal policy $\pi^*$) $\hat{Q}^* \in \argmax_{Q \in \textbf{MAR}} \sum_{t=1}^T E_Q[\max(0, x_0 - \sum_{i=1}^t D_i)]$.  We note that the measure $\hat{Q}^*$ corresponds to a strong manifestation of the obsolescence phenomena, as even in the first period the demand is either 0 or U (with the demand stuck at its initial value in all subsequent periods).  The intuition here is that when the only costs possible are holding costs, all costs are derived from being ``stuck with inventory", and the obsolescence setting (and the correlations over time in demand it induces) are (in an appropriate sense) worst-case along these lines. 
\\\indent We note that in the symmetric setting in which all holding costs are 0 while the back-logging costs may be general, the inner maximization becomes degenerate, as the policy which always orders up to U will be optimal, and no cost will be incurred under any distribution for demand (under this policy).
\end{example}
\end{document}